\documentclass[a4paper,twocolumn,final,natbib]{article}                     
\usepackage{graphicx}
\usepackage{mathptmx}      
\usepackage{epstopdf}
\usepackage{amsmath}
\usepackage{amssymb}

\usepackage{wasysym} 

\usepackage{subfigure} 

\usepackage[pagewise]{lineno}

\usepackage{ulem} 
\usepackage{dcolumn} 
\newcolumntype{.}{D{.}{.}{-1}}
\newcolumntype{d}[1]{D{.}{.}{#1}}

\usepackage{fancyhdr}

\cleardoublepage

\begin{document}

\sloppy

\title{Is the astronomical forcing a reliable and unique pacemaker for Climate?
}

%

\author{by \\ \sc  Bernard De Saedeleer\textsuperscript{\dag} 
\\ \sc 
Michel Crucifix\textsuperscript{\dag}
\\ \sc 
Sebastian Wieczorek\textsuperscript{\ddag}
    \\[0.5cm]    
    \textsuperscript{\dag} 
 \it     Georges Lema\^itre Centre for Earth and Climate Research, \\
  \it Universit\'e catholique de Louvain, Louvain-la-Neuve, Belgium                     
    \\
    \textsuperscript{\ddag} 
 \it     University of Exeter, Exeter, UK
}

%
%





\date{}

\begin{center}
  \fbox{\parbox{16cm}{{\textbf{Important note:} This is the originally submitted version.
  Readers are advised to refer to the final version, available free of charge
  at \url{http://dx.doi.org/10.1007/s00382-012-1316-1} with the following reference: }
  B. De Saedeleer, M. Crucifix, S. Wieczorek, Is the astronomical forcing a reliable and unique pacemaker for climate? A conceptual model study, Climate Dynamics,  (2013)\, 40:273--294.}}
\end{center}
 
\maketitle

\begin{abstract}

There is  evidence that ice age cycles are paced by astronomical forcing, suggesting some kind of synchronization phenomenon. Here, we identify the type of such synchronization and explore systematically its uniqueness and robustness  using a simple paleoclimate model akin to the Van der Pol relaxation oscillator and dynamical system theory. As the insolation is quite a complex quasiperiodic signal, the traditional concepts of phase- or frequency-locking used to define synchronization to periodic forcing are inadequate. Instead, we explore a different concept of {\it generalized synchronization}  in terms of (coexisting) synchronized solutions for the forced system, their basins of attraction and instabilities.

We propose a clustering technique to compute the number of synchronized solutions, each of which corresponds to a different paleoclimate history. In this way, we uncover multistable synchronization (reminiscent of phase- or frequency-locking to individual periodic components of astronomical forcing) at low forcing strength, and monostable or unique synchronization at stronger forcing. In the multistable regime, different initial conditions may lead to different paleoclimate histories. To study their robustness, we analyze Lyapunov exponents that quantify the rate of convergence towards each synchronized solution (local stability), and basins of attraction that indicate critical levels of external perturbations (global stability). We find that  even though synchronized solutions are stable on a long term, there exist short episodes of desynchronization where nearby climate trajectories diverge temporarily (for about 50 kyr). We also show that when a synchronized solution approaches the boundary of its basin of attraction,  small perturbations may cause a jump to a different (coexisting) paeleoclimate history.

Our study  brings new insight into paleoclimate dynamics and reveals a possibility for the climate system  to wander throughout different climatic histories related to preferential synchronization regimes on obliquity, precession or combinations of both, as environmental parameters varied throughout the history of the Pleistocene.



\end{abstract}

\thispagestyle{fancy}



\section{Introduction}
\label{Sec:Introduction}

\begin{figure*}
\centering
\mbox{
\subfigure[{The LR04 stack \cite{Lisiecki:2005} of 57 benthic $\delta^{18}O$ [\permil] records; the $\delta^{18}O$ is a proxy for the global volume of ice. 
High values of  $\delta^{18}O$ correspond to a colder climate (glacial state).}]
{\includegraphics[width=0.485\textwidth]{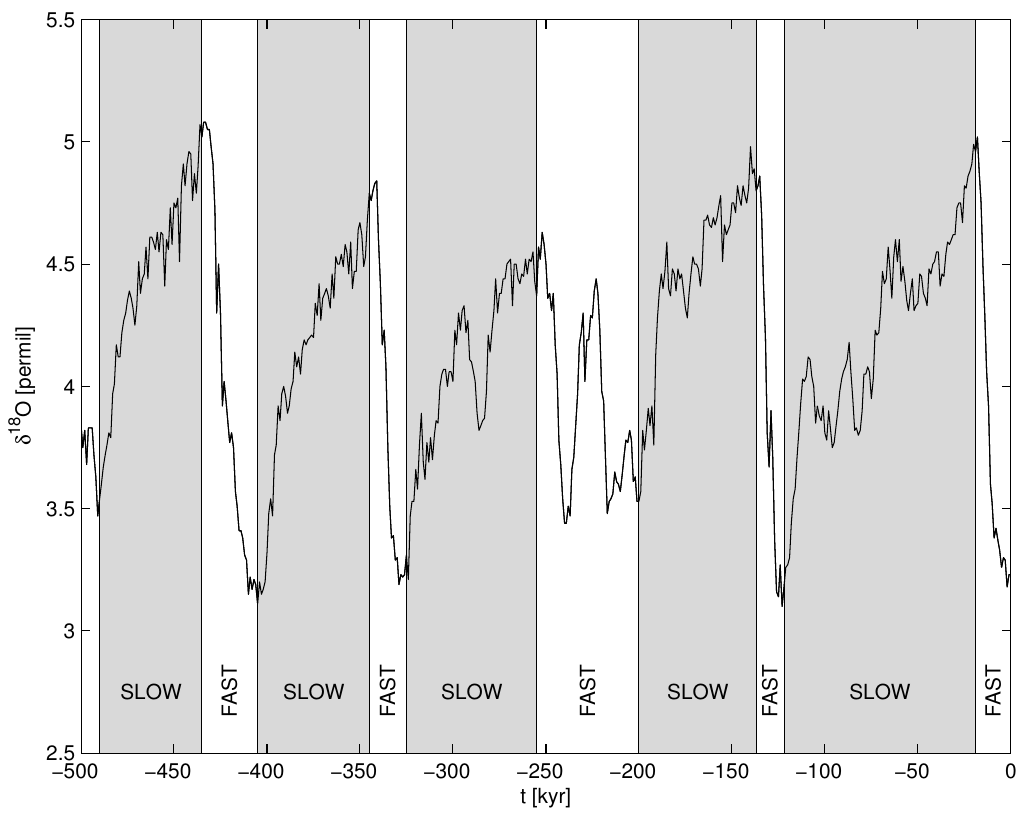}
\label{Fig:Fig_LR04_Slow_Fast}} \quad
\subfigure[{$CO_2$ composite record [ppmv] \cite{Luethi:2008}.
High values of  $CO_2$ correspond to a warmer climate (interglacial state).}]
{\includegraphics[width=0.485\textwidth]{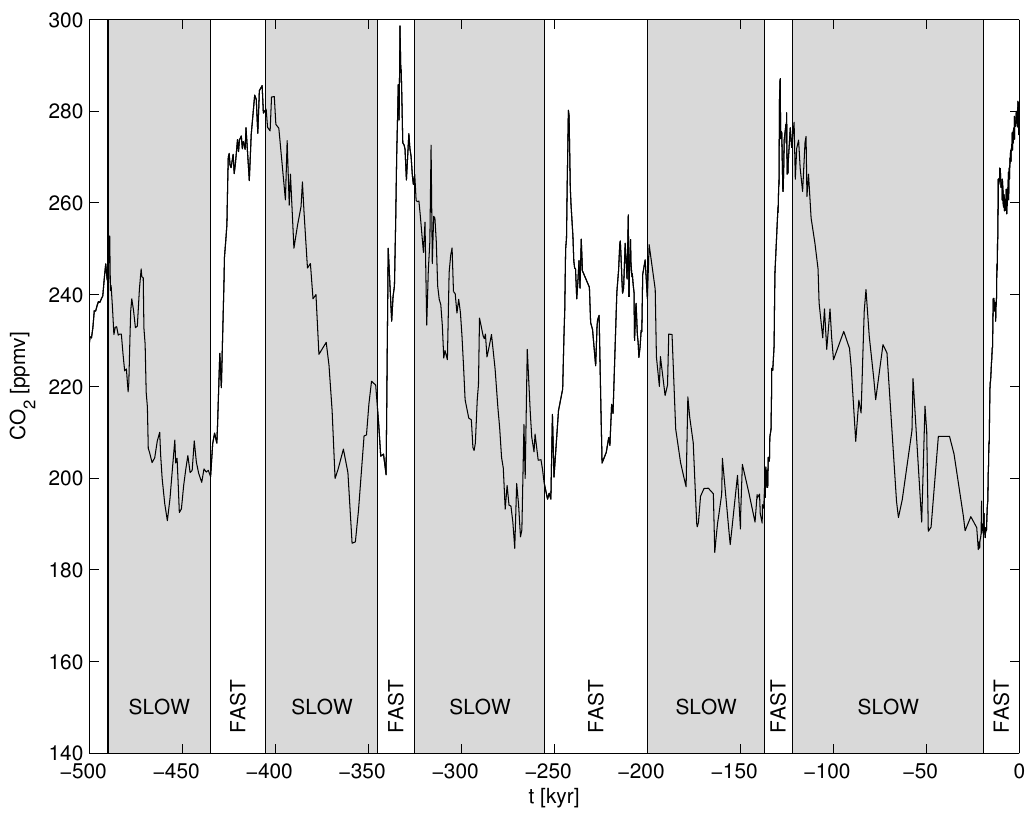}
\label{Fig:Fig_CO2_Slow_Fast}}}
\parbox{\textwidth}{
\caption{The long-term climatic signals reveal a slow-fast dynamics (only the most recent 500 kyr of the data are displayed here). The areas in grey (width  $\tau_{slow}$) are wider than the areas in white (width $\tau_{fast}$): while glaciation is a slow process of ice build-up, deglaciation occurs much more rapidly ($\tau_{slow}  > \tau_{fast}$). One also recognizes the last deglaciation which started some 20 kyr ago, up to the present time ($t=0$).}
\label{Fig:LR04_CO2}}
\end{figure*}

This article is a contribution to the field of paleoclimate dynamics theory, which has experienced many developments in terms of ice age models since many years notably by  \cite{Le-Treut:1983, Saltzman:1990}, and many others, and remains an active research field.
Paleoclimate modeling is a complex problem, hence an uncomfortable situation for a scientist. On one hand, the data are scarce and marred by uncertainties. On the other hand, there is not a single well established model, the problem is non autonomous, the forcing is aperiodic, and stochastic effects are present.

Here, we focus on the slow variations of climate over the few last million years, which include the phenomenon of ice ages
\cite{Hays:1976}, that is, the repeated growth and decay of ice sheets in the Northern Hemisphere of a total mass as big as modern Antarctica's.
When
examining long-term climatic signals like the 5.3 Myr-long stack produced in \cite{Lisiecki:2005}, or the 800 kyr-long EPICA Dome C Ice Core from \cite{Luethi:2008},
plotted respectively
in Fig. \ref{Fig:Fig_LR04_Slow_Fast}
and  Fig. \ref{Fig:Fig_CO2_Slow_Fast} for the last 500 kyr,
one immediately identifies three clearly visible features of the climatic time series:
\renewcommand{\theenumi}{\roman{enumi}}
\renewcommand{\labelenumi}{(\theenumi)}
\begin{enumerate}
\item \emph{oscillations}: the signal oscillates between higher and lower values of ice volume corresponding to the glacial and interglacial states,
\item \emph{asymmetry}: in  Fig. \ref{Fig:Fig_LR04_Slow_Fast} typical transitions from a minimum to a maximum
take much longer than transitions from a maximum to a minimum: deglaciations occur much more rapidly ($\tau_{fast} \approx 10$ kyr) than glaciations ($\tau_{slow} \approx 80$ kyr), giving a distinctive saw-tooth structure in the glacial/interglacial (G/I) cycles, especially pronounced over the last 500 kyrs,
\item \emph{100-kyr dominant period}: this has been identified by many authors since
\cite{Broecker:1970}. Note that  the G/I cycles are not  periodic.
\end{enumerate}

The asymmetry in the oscillations has been studied by many authors.
In order to reproduce it,
some authors use  underlying physical principles to build phenomenological models that exhibit slow-fast dynamics reasonably mimicking the climatic proxies  \cite{Saltzman:2002}.
Others assume this asymmetry by explicitly defining 2 different parameters such as the time intervals $\tau_{up} = \tau_{slow}$ and $\tau_{down} =\tau_{fast}$ \cite{Ashkenazy:2006}
or time constants ($\tau_{R}$ and $\tau_{F}$ in \cite{Paillard:1998} and
$T_{w}$ and $T_{c}$ in \cite{Imbrie:1980}).
Whatever the model, it has to ultimately exhibit asymmetric oscillations under the effect of the forcing,
as it is aimed to mimic the oscillations between G/I states. Relaxation oscillators are therefore very
straightforward natural candidates of ice age models.
In this article, we will consider a slightly modified van der Pol oscillator model to illustrate the new contributions
of our synchronization concepts. 

In this paper,
we concentrate on the influence of the astronomical forcing on Earth's climate.
This forcing
is induced by the slow variations in the spatial and seasonal distributions of incoming solar radiation (insolation) at the top of the atmosphere, associated with the slow variations of the Earth's astronomical elements: eccentricity ($e$), true solar longitude of the perihelion measured with respect to the moving vernal equinox ($\varpi$), and Earth obliquity ($\epsilon_E$). These quantities are now accurately known over several tens of millions of years \cite{Laskar:2004}, but analytical approximations of $e$,  $e\sin \varpi$, and $\epsilon_E$ valid back to one million years have been known since \cite{Berger:1978}. They take the form of d'Alembert series ($\sum A_i \sin[\omega_i t + \phi_i]$).

\begin{figure}
\begin{center}
\includegraphics[width=0.48\textwidth]{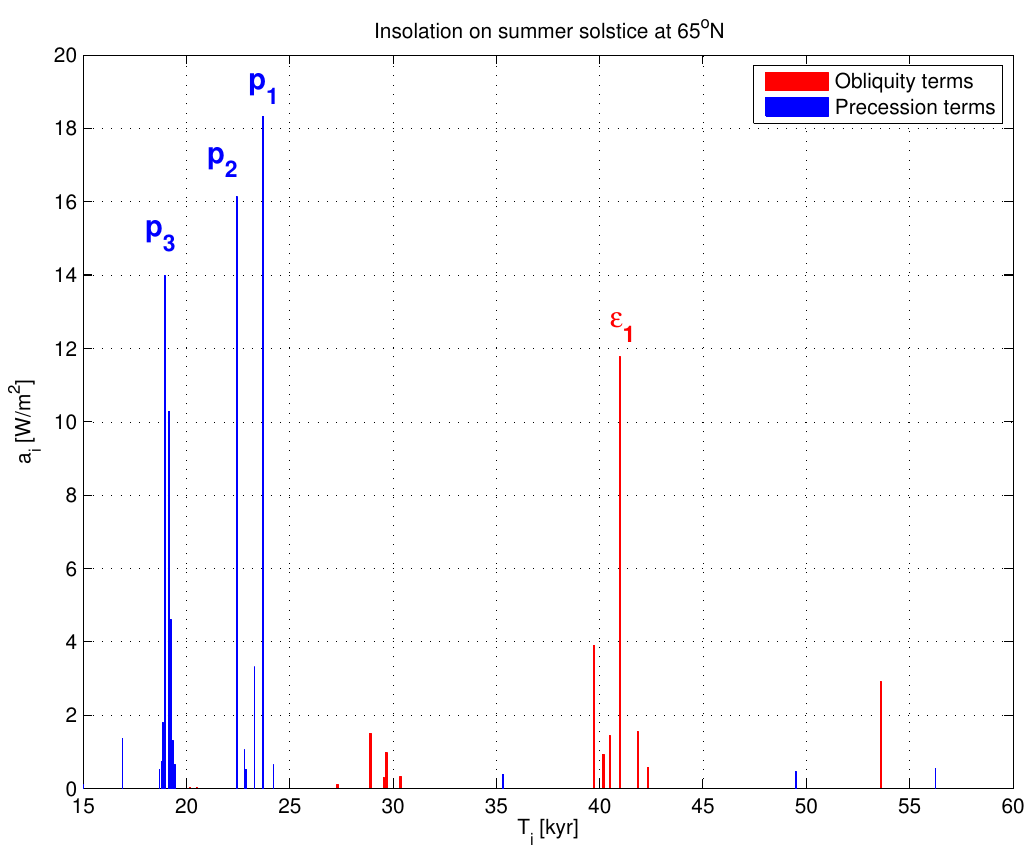}
\caption{Spectrum decomposition of the
incoming solar radiation (insolation)  at 65$^\circ$N latitude on the day of the summer solstice
($F(t)$ given in Eq. (\ref{Eq:Insolation_35_terms})).
This is a graphical representation of Table~\ref{Tab:Insolation35} in Appendix \ref{Annex:Insolation},
with $a_i^2=s_i^2 + c_i^2$ and $T_i=2\pi/\omega_i$.
The eight largest parameters $a_i$, which represent already $80\%$ of the signal, are the major components of the insolation; they clearly come from the precession (19 and 23 kyr), and from the obliquity (41 kyr) associated series. This insolation is the forcing used to construct all Figures subject to the quasiperiodic astronomical forcing.
The main harmonic $\epsilon_1$ associated with obliquity has an angular velocity
$\omega_{\epsilon_1} = 0.1532$ rad/kyr
($T_{\epsilon_1} = 41.0$ kyr)
and an amplitude of $a_{\epsilon_1} = 11.77 ~\text{W/m}^2$.
The three main harmonics associated with precession are denoted $p_1$, $p_2$ and $p_3$.
}
\label{Fig:Fig_spectre_insol}
\end{center}
\end{figure}

The external forcing $F(t)$ used throughout this article is the insolation at 65$^\circ$N latitude on the day of the summer solstice.
That specific insolation quantity is commonly related to the Milankovitch theory and can be thought of as a measure of how much ice may melt over summer.
It can be written under the following compact form:
\begin{equation}
F(t) = \sum_{i=1}^{35} [ s_i \,\sin(\omega_i t) + c_i \,\cos(\omega_i t) ]
\label{Eq:Insolation_35_terms}
\end{equation}
where the value of the $3\times 35$ parameters (including $\omega_i, s_i,$ and  $c_i $) are given in the Table~\ref{Tab:Insolation35} of Appendix \ref{Annex:Insolation}.
The coefficients were extracted from \cite{Berger:1978}
by performing a linear regression of the insolation on the $\omega_i$.
The  validity range of this  approximation is [-1 Myr, 0 Myr], and
its mean error 
(mean of the absolute value of the difference, compared to \cite{Laskar:2004})
is 6.7 W/m$^2$ with peaks at 27.5 W/m$^2$.
Note that the mean value (494.2447 W/m$^2$) has been removed;
the theoretical framework that allows to work with anomalies was justified by
\cite{Saltzman:1991}. In short, this can be done as we are interested in oscillations,
and not in the mean values themselves.
The \emph{quasiperiodic}\footnote{
A \emph{quasiperiodic} signal is the superposition of several periodic signals with uncommensurate periods.
}
nature of the  insolation forcing $F(t)$ is illustrated in its spectrum decomposition
in Fig. \ref{Fig:Fig_spectre_insol}.  Precession is dominated by two harmonics around 19 and 23 kyrs (1 kyr = 1,000 years) and obliquity is dominated by an harmonic with a period of 41 kyrs but it bears periods as long as 1,200 kyrs.

\paragraph{Synchronization\\[0.25cm]}

There is ample evidence that the astronomical forcing influences the climate system. The phrase 'pacemaker of ice ages' was coined in a seminal paper  \cite{Hays:1976} to express the idea that the timing of ice ages is controlled by the astronomical forcing, while the ice age cycle itself is shaped by internal system dynamics. The paradigm has prevailed since then and is it still supported by the most recent analyses of palaeoclimate records \cite{Lisiecki:2007,Huybers:2007}.
The notion of 'pacemaker' naturally evokes some sort of synchronization.
However, despite some attempts, the actual
type of synchronization has not been clearly identified or demonstrated to date.
For example, \cite{Ashkenazy:2006,Tziperman:2006} speak of "nonlinear phase-locking" although they
do not define  suitable ''phase variables'' that can be used to demonstrate a fixed-in-time relationship between
phases of the forcing and the oscillator response.



Synchronization, as a universal nonlinear phenomenon, is a pervasive process in Nature, as it is associated with rhythmic processes. It is therefore not surprising to have synchronization also in Paleoclimatic Sciences.
Depending on the forcing type (periodic, chaotic, stochastic), one can distinguish  many types of synchronization including
complete, lag, phase, frequency, identical, generalized \cite{Rulkov:1995}, achronal and isochronous \cite{Wu:2006},
and even noise synchronization.
For a review, the reader is referred to \cite{Balanov:2009, Pikovsky:2001}, among others.
While some terminology is still debated,   \cite{Brown:2000} proposed an unified definition of synchronization
for dynamical systems---there is synchronization if there exists a relationship $h$ between the measured properties of the forcing, $g(\vec{u})$, and those of the oscillator, $g(\vec{v})$: \begin{eqnarray}
\label{Eq:synch}
h(g(\vec{u}), g(\vec{v}))=0,
\end{eqnarray}
that is fixed-in-time, meaning that $h$ is time independent.
Because we are interested in synchronization that is stable,
for arbitrary initial conditions $\vec{u}(0)$ and $\vec{v}(0)$ that do not satisfy  (\ref{Eq:synch}),
we require  that \cite{Brown:2000}:
\begin{eqnarray}
\label{Eq:synchstab}
\lim_{t \rightarrow \infty}\; h(g(\vec{u}(t)), g(\vec{v}(t))) = 0.
\end{eqnarray}
For example, if $g(\vec{u})=\vec{u}$,  $g(\vec{v})=\vec{v}$, $\vec{u}$ and $\vec{v}$ have the same dimension, and
(\ref{Eq:synch}) can be written as $\vec{u}=\vec{v}$, we speak of {\it identical synchronization}.
More generally, if vectors $\vec{u}$ and $\vec{v}$ have different dimensions and (\ref{Eq:synch}) cannot be reduced
to more than a functional relationship $\vec{u}=H(\vec{v})$, we speak of {\it generalized synchronization}; see also
\cite{Abarbanel:1996,Rulkov:1995,Pikovsky:2001}.
Note that  the relationship (\ref{Eq:synch}) need not be unique. If there are two or more relationships
(\ref{Eq:synch}) for the same parameter settings, we speak of {\it multistable synchronization}
\cite[Ch.15.3.2]{Pikovsky:2001}. Then, which of the relationships the system settles to will depend
on initial conditions.

In this paper, we use a simple van der Pol oscillator model to identify and illustrate
 for the first time the phenomenon of generalized
synchronization between ice age cycles and astronomical forcing. The dynamical systems approach
outlined in the next section (i) allows for stability analysis of such synchronization,
(ii) uncovers interesting effects relating to the robustness of the synchronization with respect to
external perturbations, and (iii) uncovers the phenomenon of multistable synchronization that
has been overlooked by previous studies. We show that, in contrast to claims in  \cite{Tziperman:2006},
synchronization needs not be unique.

The article is structured as follows. Section \ref{Sec:VDP} introduces a slightly modified version of the van der Pol oscillator as a suitable model for studying synchronization of ice ages to astronomical forcing.
In Section \ref{Sec:oscillators}, we analyse synchronization to periodic forcing and quasiperiodic astronomical forcing
in terms of largest Lyapunov exponents.
Section \ref{Sec:Clustering}
is dedicated to the study of multistable synchronization in terms of attracting trajectories in the phase
space of the forced system, and the associated basins of attraction.
In Section \ref{Sec:beta}, we investigate effects of the symmetry-breaking parameter $\beta$ for the van der Pol oscillator model.
Section \ref{Sec:Robustness} is concerned with the robustness of the synchronization and focuses on two aspects relating to predictability. Firstly, it shows that the local stability can be lost temporarily causing divergence of nearby climatic trajectories. Secondly, it demonstrates that in the multistable regime external perturbations (such as noise) may cause jumps between coexisting synchronized solutions when these solutions come close to their basin boundary.  To be clear, all the treatment below is deterministic, except for Figs.
\ref{Fig:Spread} and \ref{Fig:stochastic}.


This article requires some basics of Dynamical Systems theory (dynamical systems, nonlinear oscillations,  limit cycles, bifurcations of vector fields, etc.), for which
we refer the reader to \cite{Guckenheimer:1983, Arnold:1983, Strogatz:1994},
to \cite{Saltzman:2002} for dynamical paleoclimatology,
and also to \cite{Savi:2005} for a review of some useful concepts.

\section{Generic ice age model: a modified van der Pol relaxation oscillator}
\label{Sec:VDP}

The hypothesis at the basis of the work by Milankovitch \cite{Milankovitch:1941} is that changes in total amount of continental ice (say: $x$) are driven by summer insolation $F(t)$ already described in Eq.~(\ref{Eq:Insolation_35_terms}).
One straightforward interpretation of this hypothesis is a simple differential equation
$\dot x = -d\Psi(x)/dx - \gamma ~F(t)$, where $d\Psi(x)/dx$ is the derivative of a climatic potential and $\gamma$ is the \emph{forcing efficiency}. However, models of this form fail in practice to correctly capture the rapid deglaciation phenomenon. We therefore propose to model the paleoclimatic dynamical system with a dissipative self-sustained oscillator resembling the classical van der Pol oscillator\footnote{
We give in the Appendix \ref{Annex:vdp} a summary description of the classical van der Pol oscillator  
and of its dynamical behaviour.
}:
\begin{subequations}
\label{Eq:System}
\begin{eqnarray}
\tau ~\dot x &=& -  [\; y + \beta - \gamma   \; F(t) \;] \label{Eq:x}   \\
\tau ~\dot y &=&  - \alpha \; [\; \Phi'(y) - x  \;]  \label{Eq:y}
\end{eqnarray}
\end{subequations}
where $\Phi'(y) = y^3 / 3 -y$. 
Note that this system is nonautonomous because the right-hand side depends explicitly on time.

The  physical interpretation of the model is as follows. Ice volume $x$  integrates the external forcing $F(t)$ over time but with a drift $y + \beta$. Assuming $\alpha \gg 1$, $y$ is the faster variable whose dynamics is controlled by a two-well potential $\Phi(y)$. For example, there are arguments that the dynamics of the Atlantic ocean circulation may be approximated by an equation similar to Eq. \ref{Eq:y}
\cite{Rahmstorf:2005,Dijkstra:2003}.
Further interpretation and discussion of the fast variable
can be found in \cite{Saltzman:1984,Tziperman:2003,Paillard:2004,Tziperman:2006,Crucifix:2011a}.
The parameter $\tau$ sets the slow time scale.
The coupled system Eq. (\ref{Eq:System}) has one stable equilibrium solution for $| \beta | > 1$ and a stable periodic orbit for $| \beta | < 1$. The ratio of time spent near the two stable branches of the slow manifold given by $\Phi'(y) = x$ depends on  $\beta$ (see the paragraph "Time Spent" in  Appendix \ref{Annex:vdp}).
We use $T_{ULC} $ to denote the period of the stable periodic orbit
and $\omega_{ULC} = 2\pi/T_{ULC} $ to denote the corresponding  angular velocity.

\begin{figure}
\begin{center}
\includegraphics[width=0.48\textwidth]{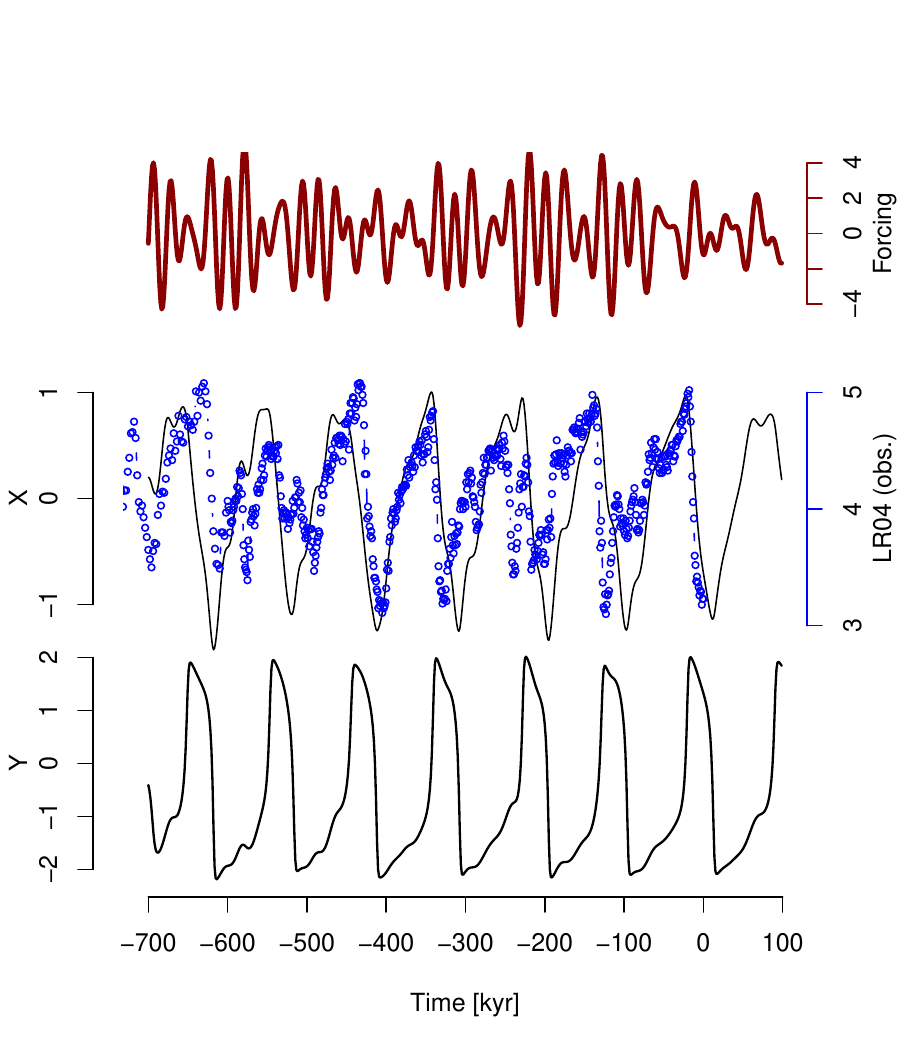}
\caption{
\emph{Top}: the insolation forcing $F(t)$, scaled by $a_{\epsilon_1}$ in order to work dimensionless. \emph{Bottom}: the $x$ and $y$ climatic trajectories obtained using system  Eq. (\ref{Eq:System})
with $\alpha = 11.11$, $\beta=0.25$, $\gamma = 0.75$ and $\tau = 35.09$. With these parameters, $\omega_{\epsilon 1} = 2.5 \omega_{ULC}$, where  $\omega_{\epsilon 1}$ is the angular velocity associated with the dominant harmonic of obliquity and $\omega_{ULC}$ is the angular velocity associated with the unforced system's periodic orbit. Blue dots correspond to the Lisiecki and Raymo stack (LR04) described in Fig.~\ref{Fig:Fig_LR04_Slow_Fast}. Time $t= 0$ corresponds conventionally to the year 1950.
The model is in good agreement with the ice volume proxy.
}
\label{Fig:Fit_VDP}
\end{center}
\end{figure}

Relaxation oscillators have been proposed previously to study ice ages \cite{Saltzman:1984,Tziperman:2003,Paillard:2004} although, to our knowledge, in a less general form than here.
We adopted this form\footnote{Note that the van der Pol oscillator model is also used as a reference in \cite[page 101]{Saltzman:2002} for a coupled ocean/sea-ice model.}
because it is very close to the well-studied van der Pol oscillator, and a good agreement 
(timing of glaciations and deglaciations, and their amplitude)
with ice volume proxies was easily found for well chosen values of $\alpha, \beta, \gamma$ and $\tau$ (Fig. \ref{Fig:Fit_VDP}).
We note, though, that small changes in parameters or additive fluctuations may easily shift the timings of ice-age terminations for reasons that will be clarified later in the paper.

The definition of synchronization can be applied to our model Eq. (\ref{Eq:System}) as follows. The astronomical forcing $F(t)$ corresponds to $\vec{u}(t)$, and the state vector whose two components are the slowly-varying ice volume $x$ and the faster variable $y$ corresponds to $\vec{v}(t)$. For nonperiodic forcing, relationship (\ref{Eq:synch}) can be very complicated (non-functional or even fractal-like) and hence difficult to detect. Therefore, other methods of detecting (\ref{Eq:synch}) had to be developed. As suggested by the auxiliary system approach \cite{Abarbanel:1996}, relationships (\ref{Eq:synch}) and (\ref{Eq:synchstab}) are implied by an (invariant) \emph{attracting trajectory}  in the $(x,y,t)$ phase space of the nonautonomous forced system (\ref{Eq:System}) \cite{Wieczorek:2011}. In the remainder of the paper,  such an attracting trajectory  is denoted with $AT$ and referred to as an \emph{attracting climatic trajectory} or {\it synchronized solution}. All other solutions to Eq. (\ref{Eq:System}) will be referred to as \emph{climatic trajectories}.

Previous approaches to nonlinear dynamics of quasiperiodically forced oscillators focused on discrete-time  mappings and two-frequency forcing \cite{Glendinning:1999,Osinga:2000,Belogortsev:1992,Broer:1998}. They uncovered interesting dynamics including Arnol'd or mode-locked tongues consisting of `interlocking' bubbles and open regions of multistability, nonsmooth bifurcations, and strange nonchaotic attractors. Here, we consider quasiperiodic forcing with 35 frequency components and focus on the regions of mode locking. Our  approach is based on instabilities of attracting trajectories in the $(x,y,t)$ phase space of the continous-time forced system  because they relate directly to the concept of generalized synchronization. We can provide a systematic study of generalized synchronization to astronomical forcing by demonstrating existence of such trajectories and exploring their {\it local} and {\it global} stability properties. More specifically, we perform three types of calculations. Firstly, a clustering detection technique uncovers parameter regions with monostable (unique) and multistable (non-unique) synchronization. Secondly, the largest Lyapunov exponent along $AT$ quantifies its long- and short-term local (linear) stability. Thirdly, a basin of attraction of $AT$  quantifies its global (nonlinear) stability. Finally, we remark that in the theory of nonautonomous dynamical systems,  attracting trajectories in the $(x,y,t)$ phase space are linked to a modern and more general concept of a \emph{pullback attractor} \cite{Kloeden:2000,Langa:2002,Kosmidis:2003,Wiggins:2003}.

\begin{figure*}
\centering
\mbox{
\subfigure[Time series : $x$ (the ice volume) is the slow variable, while $y$ exhibits slow-fast dynamics. The same color convention white/grey as in Fig. \ref{Fig:LR04_CO2} has been used in order to highlight the slow and fast episodes.]
{\includegraphics[width=0.485\textwidth]{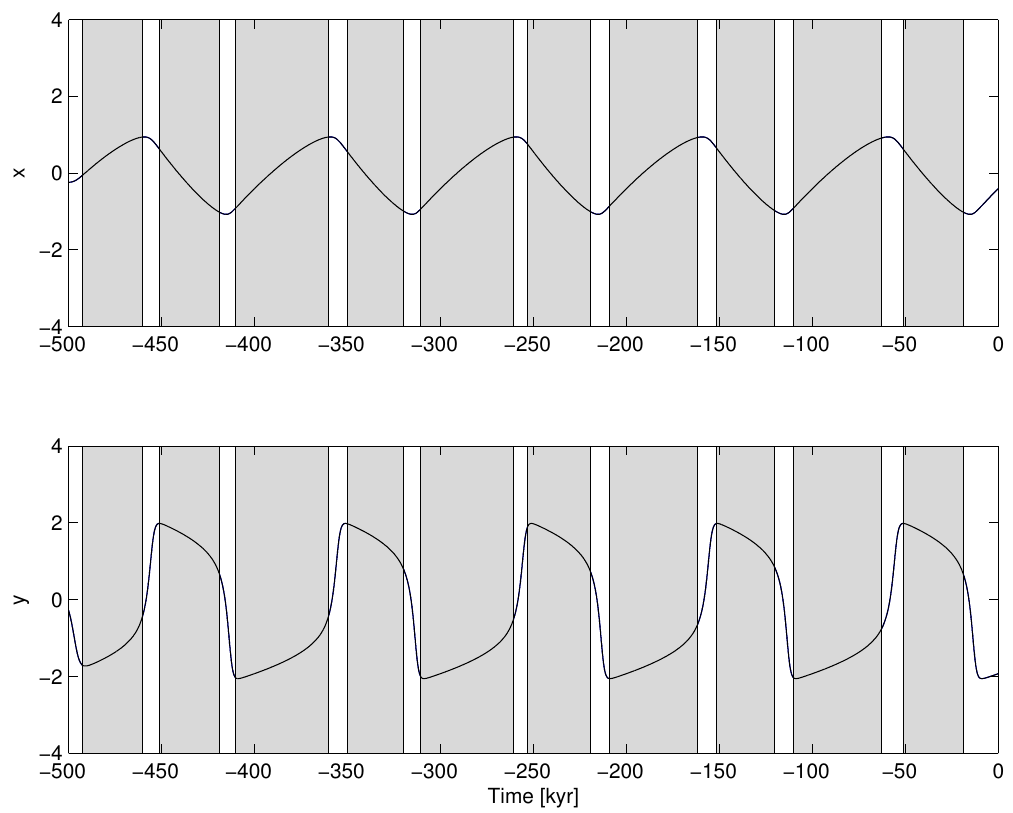}
\label{Fig:vdp_eps_0_1_phi3_cycle_flow}} \quad
\subfigure[Phase space portrait  in the  two-dimensional $(x,y)$ phase space  of the autonomous system: dynamical flow, limit cycle and residence plot (circles are spaced at every 0.5 kyr). More time is spent along the slow manifold (red dashed curve, corresponding to the function $\Phi'(y) =  y^3 / 3 -y = x$). 
The trajectory converges to the limit cycle (thin curve) with slow-fast dynamics.]
{\includegraphics[width=0.485\textwidth]{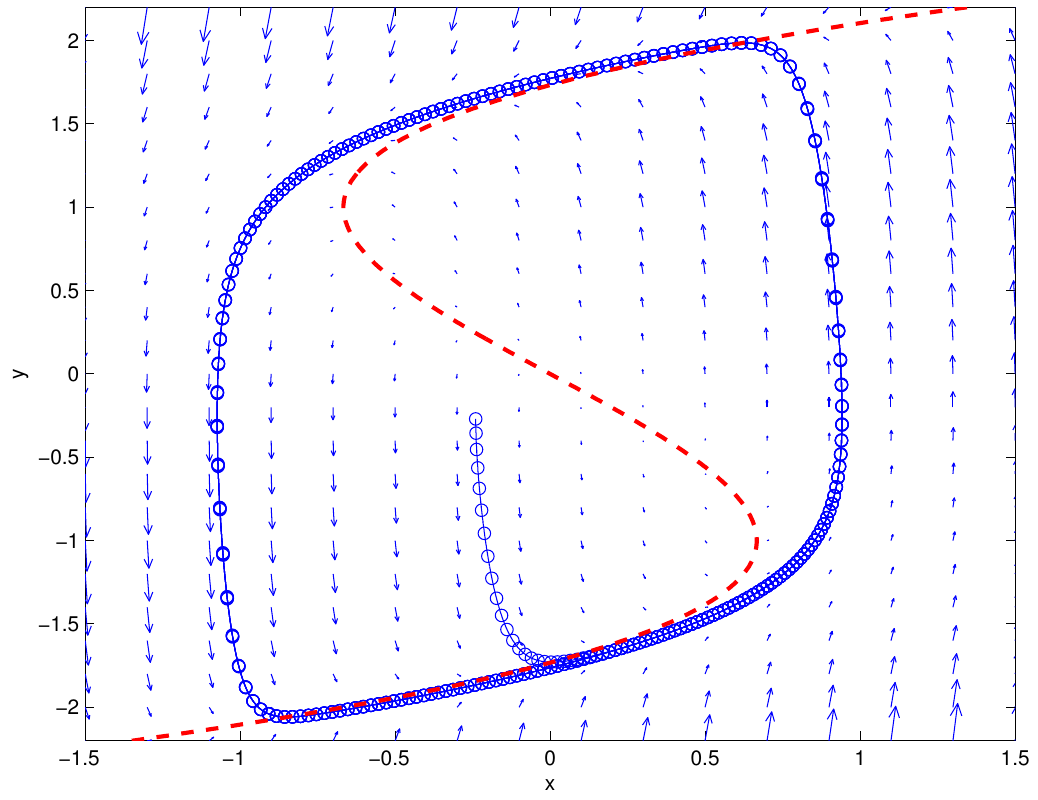}
\label{Fig:vdp_eps_0_1_phi3_xy}}}
\parbox{\textwidth}{
\caption{Dynamics of the unforced ice ages model Eq. (\ref{Eq:System}) with $\alpha = 11.11$, $\beta=0.25$, $\gamma = 0$ and $\tau = 35.09$. The slow-fast variable is $y$, while $x$ is always slow. Example of a typical climatic trajectory with initial conditions ($x_{-500},y_{-500}) = (-0.24, -0.27)$.}
\label{Fig:vdp_03}
}
\end{figure*}




\section{Synchronization of the paleoclimatic system to the insolation forcing}
\label{Sec:oscillators}

\begin{figure*}[!hp]
\centering

\mbox{
\subfigure[without any forcing ($\gamma = 0$, $\tau = 35.09 \rightarrow T_{ULC} \approx$ 100 kyr): 
all trajectories converge to a cylinder---the attracting set in the $(x,y,t)$ space.]
{\includegraphics[width=0.48\textwidth]{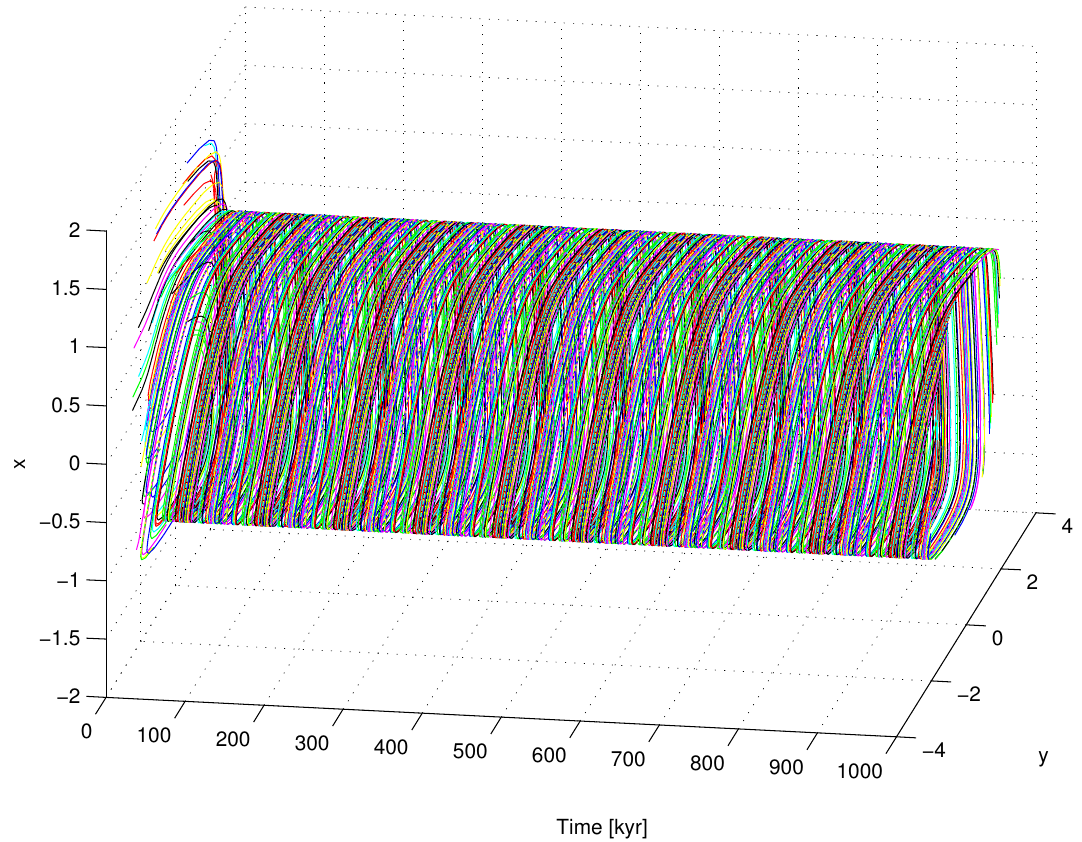}
\label{Fig:No_Insol}} \quad
\subfigure[{Section of Fig. \ref{Fig:No_Insol} at $t=550$ kyr: no clear dense cluster of climatic trajectories is identified: no attracting trajectory does exist.}]
{\includegraphics[width=0.48\textwidth]{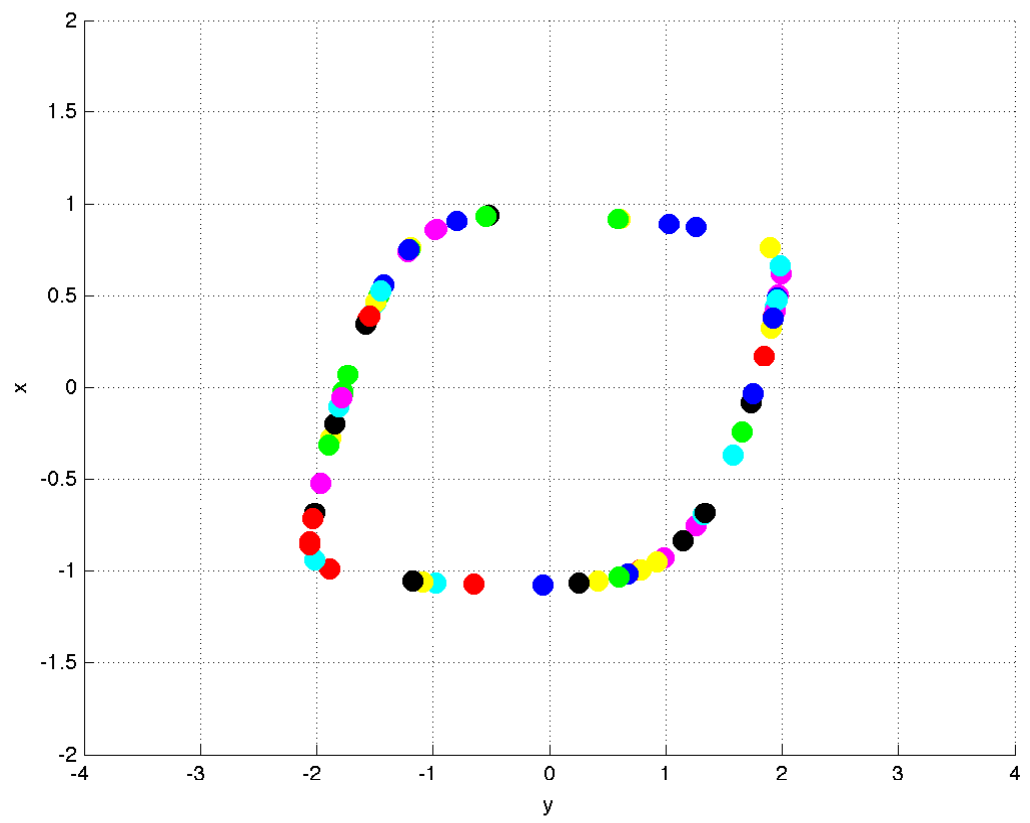}
\label{Fig:No_Insol_section}}}

\mbox{
\subfigure[with a 
purely periodic forcing with a period  of $T_{F} =$ 41 kyr 
  ($\gamma = 3.33$, $\tau = 35.09 \rightarrow T_{ULC} \approx$ 100 kyr): 
 the trajectories converge to two attracting trajectories of period 82 kyr:  
  there is a frequency-locking 2:1.]
{\includegraphics[width=0.48\textwidth]{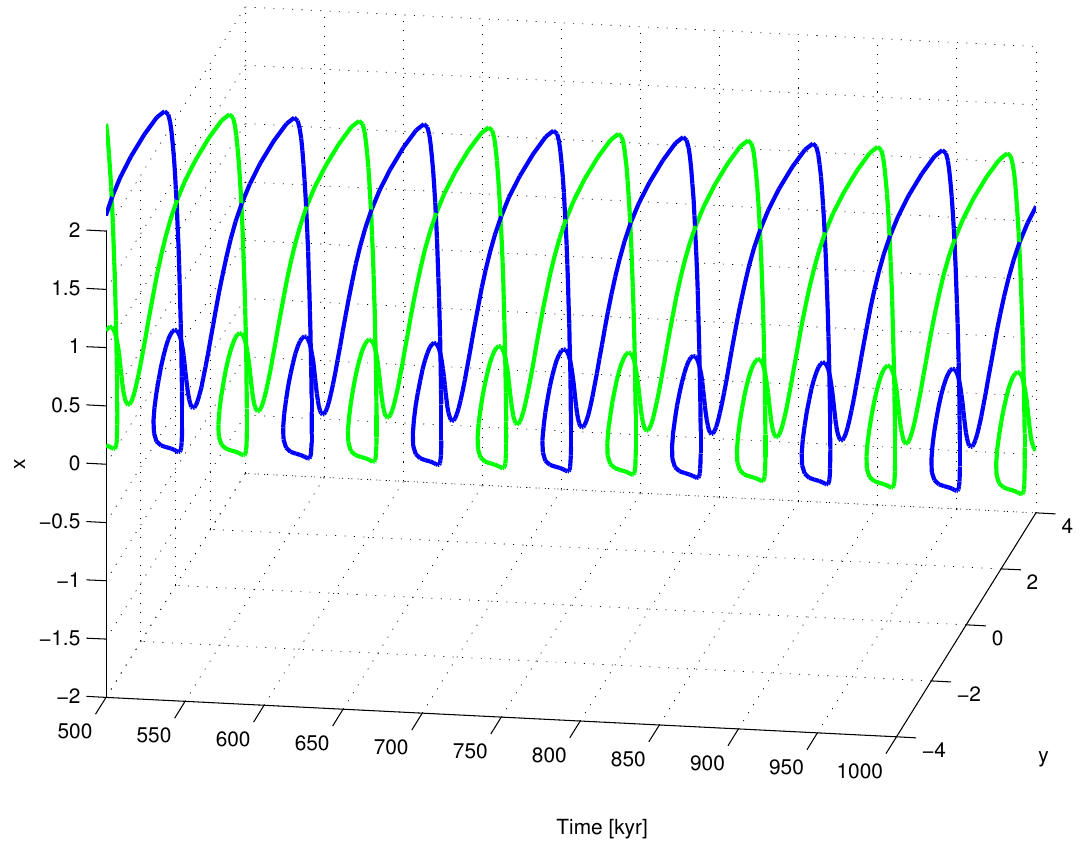}
\label{Fig:Insol_41k}} \quad
\subfigure[{Section of Fig.~\ref{Fig:Insol_41k} at $t=550$ kyr: two clusters are identified, corresponding to the two attracting trajectories born from frequency-locking 2:1.}]
{\includegraphics[width=0.48\textwidth]{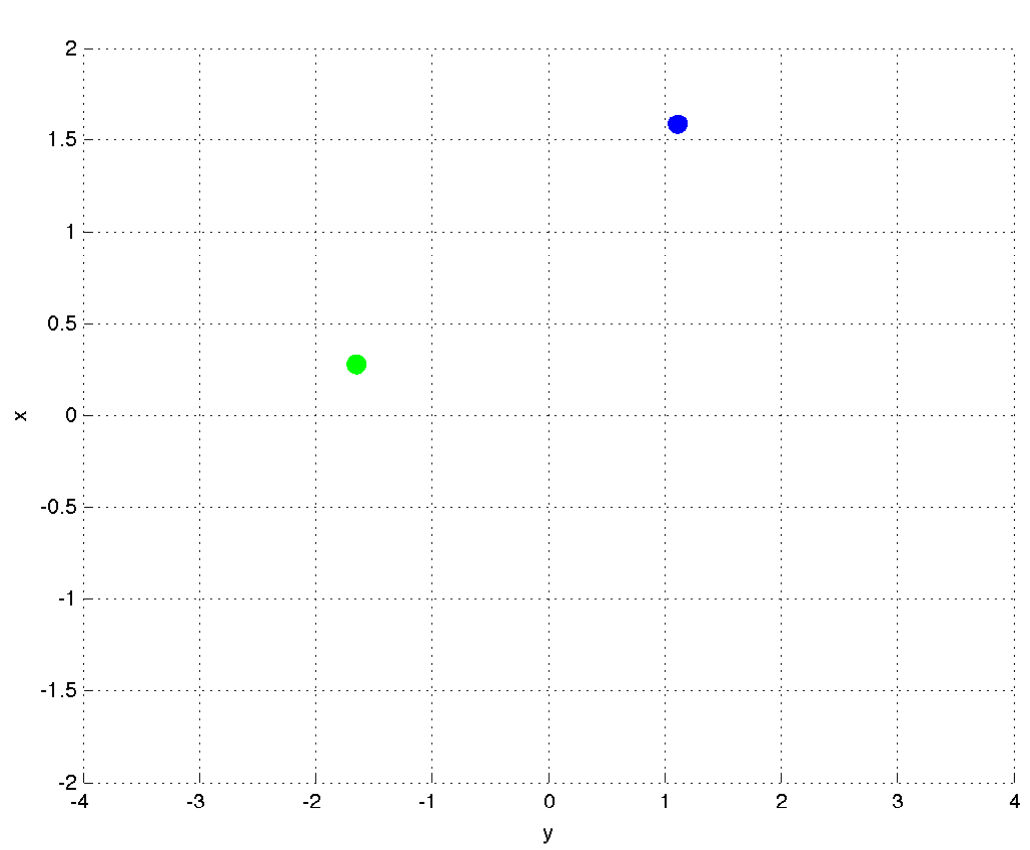}
\label{Fig:Insol_41k_section}}}

\mbox{
\subfigure[with the quasiperiodic insolation ($\gamma = 0.75$, $\tau = 43.86 \rightarrow T_{ULC} \approx$ 125 kyr) given by Eq.~(\ref{Eq:Insolation_35_terms}): the trajectories converge now to three attracting trajectories for a long time, revealing a multistable synchronization.]
{\includegraphics[width=0.48\textwidth]{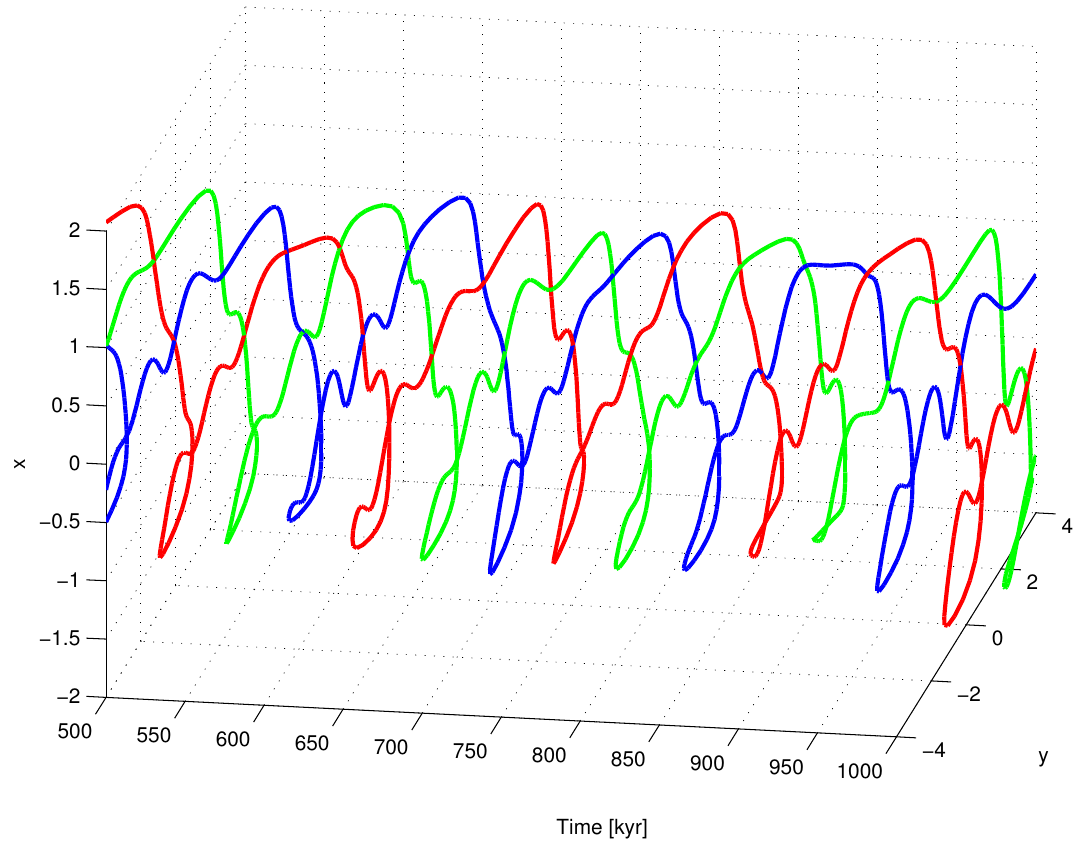}
\label{Fig:Insol}} \quad
\subfigure[{Section of Fig.~\ref{Fig:Insol} at $t=550$ kyr: three clusters of trajectories are identified, corresponding to the three attracting trajectories born from multistable synchronization.}]
{\includegraphics[width=0.48\textwidth]{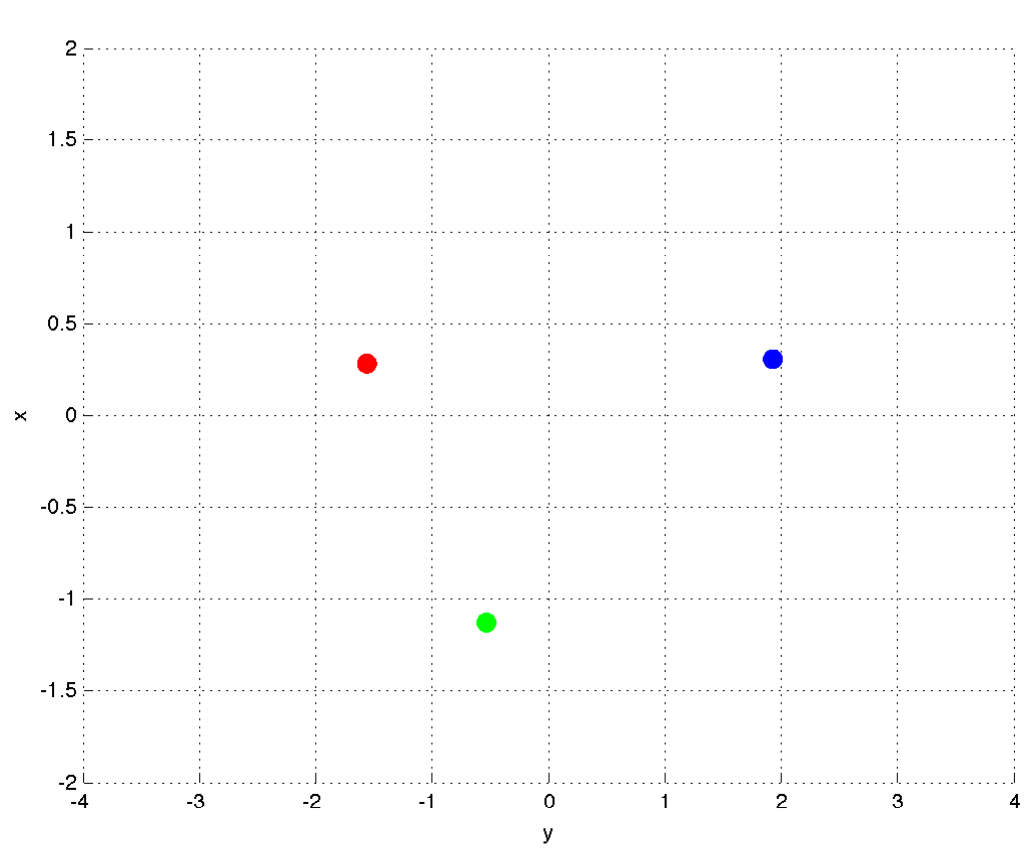}
\label{Fig:Insol_section}}}

\parbox{\textwidth}{
\caption{Illustration of three different paleoclimate dynamical regimes: (\emph{top}) without any forcing, (\emph{middle}) frequency-locking 2:1 onto a 41 kyr periodic forcing (symbol '+' in Figs. \ref{Fig:Fig_bifurc_b} and \ref{Fig:Fig_bifurc_Insol_41k}), (\emph{bottom}) generalized multistable synchronization on the quasiperiodic insolation forcing (symbol '$\times$' in Figs. \ref{Fig:Fig_bifurc_c} and \ref{Fig:Fig_bifurc_Insol}). The ice age model used is Eq. (\ref{Eq:System}) with $\alpha = 11.11$ and $\beta=0.25$.
A set of 70 random initial conditions at $t_0=0$ is used.
\emph{Left}: climatic trajectories in the $(x,y,t)$ space-time space. \emph{Right}: section of these trajectories at  $t=550$ kyr, which reveals clearly the potential formation of the attracting trajectories and allows an easier counting of these.}
\label{Fig:Many_2_3_Attractors}}
\end{figure*}

\paragraph{Illustration of the synchronization phenomenon  \\[0.25cm]}
A typical climatic trajectory for $\gamma=0$ is shown in Fig. \ref{Fig:vdp_03},
from two different points of view: the time series and the phase space portrait.
In the time series, we recognize the slow variable $x$ (the ice volume),
while $y$ exhibits slow-fast dynamics.
This climatic trajectory is also shown in the  two-dimensional $(x,y)$ phase space  of the autonomous system where arrows
indicate direction of the flow. The trajectory converges to the limit cycle with slow-fast dynamics
(the speed along the trajectory can be visually assessed by the circles of the residence plot).
Let us now consider a set of 70 random initial conditions in the $(x,y)$-plane
at time $t_0=0$, and study the resulting climatic trajectories in the three-dimensional phase space $(x,y,t)$
of the nonautonomous system (Fig. \ref{Fig:No_Insol}) for time $t>t_0$.
One clearly sees that all trajectories converge
to a cylinder---the attracting set in the $(x,y,t)$ space.

However, if we consider now an external forcing ($\gamma>0$) then \emph{synchronization} onto this forcing may
occur  under certain conditions \cite{Ashkenazy:2006,Tziperman:2006}. According to our definition
(\ref{Eq:synch}--\ref{Eq:synchstab}),  synchronization is represented by an attracting climatic trajectory
in the  $(x,y,t)$ phase space.

Consider first the case of a purely
periodic forcing with a period  of $T_{F} =$ 41 kyr and strength $\gamma = 3.33$.
The 70 initial conditions give rise to climatic trajectories that, after a sufficiently long integration time, converge to two \emph{attracting trajectories}  (see Fig. \ref{Fig:Insol_41k}). Both attracting trajectories  are periodic with period of $T_R=2 T_F= 82$ kyr, and time-shifted versions of each other.
This phenomenon is described in the literature as $2:1$ \emph{phase-locking} or \emph{frequency-locking}. Generally speaking  a  $n:m$  $frequency$-$locking$ is defined as a fixed-in-time relation between the frequencies of the forcing ($\omega_F$) and the oscillator response ($\omega_R$) of the form
$
n \,\omega_R = m\, \omega_F,
$
where $m$ and $n$ are integers \cite[p.52]{Pikovsky:2001}.

Then consider the case of the quasiperiodic insolation forcing described in Eq.~(\ref{Eq:Insolation_35_terms}) with
$\gamma = 0.75$  and $\tau = 43.86$.
Figure \ref{Fig:Insol} shows that the 70 climatic trajectories now converge onto three \emph{attracting trajectories},
which reveals that synchronization can be \emph{multistable} \cite[p.348]{Pikovsky:2001}, \cite[p.94]{Balanov:2009}.
This phenomenon is described in the literature as \emph{mode-locking}  \cite{Svensson:2009}.
Note that because of the quasiperiodicity of the insolation forcing, these attracting trajectories are  \emph{no longer} periodic
nor time-shifted versions of each other. The number  of attracting trajectories  depends on many factors including the dynamics of the unforced system, the nature of the forcing $F(t) $, and the amplitude $\gamma$ of the forcing. We will study this in more details in Sec.~\ref{Sec:Clustering}.

\paragraph{Detection of  synchronization by the way of the largest Lyapunov exponent (LLE  or $\lambda_{max}$) \\[0.25cm]}

Local or linear stability of an attracting climatic trajectory can be quantified with the largest Lyapunov exponent (LLE) denoted here as $\lambda_{max}$ \cite{Benettin:1980fk}. The quantity $\lambda_{max}$ is a measure of the (average) exponential rate of divergence ($\lambda_{max}>0$) or convergence ($\lambda_{max}<0$) of nearby climatic trajectories. Therefore, a negative value of  $\lambda_{max}$ indicates a locally attracting climatic trajectory or \emph{generalized synchronization} \cite{Pikovsky:2001,Wieczorek:2009}. A transition from
$\lambda_{max}<0$ to $\lambda_{max}=0$ indicates a bifurcation where the attracting climatic trajectory disappears and generalized synchronization is lost. Null and positive values of  $\lambda_{max}$ indicate lack of synchrony (positive $\lambda_{max}$ indicates chaos but this regime is not encountered here). In the case of periodic forcing, computations of $\lambda_{max}$ can be easily validated with more precise and reliable numerical bifurcation continuation techniques (see \S ~'\emph{41 kyr periodic forcing}' below).

\paragraph{Long-term $ \lambda_{max}$ and short-term
$ \lambda^H_{max}$ LLE's \\[0.25cm]}

The largest Lyapunov exponent $\lambda_{max}$ is mathematically defined\footnote{Even if differential versions of the LLE have sometimes been developed mainly for computational efficiency purposes, we however preferred within this article to stick on the original definition of the LLE, because it is more standard and there is no insistent need for lowering computation time in the present framework, as the number of degrees of freedom of the system is reduced.} as  \cite{Ott:2002} :
\begin{equation}
\lambda_{max} =
\lim_{|\delta\bold{Z}(0)|\rightarrow 0}~
\lim_{t \rightarrow \infty}~
\frac{1}{t} \ln
\frac{|\delta\bold{Z}(t)|}
{|\delta\bold{Z}(0)|}
\label{Eq:lyap_global}
\end{equation}
where $\delta\bold{Z}=[\delta x,\delta y]$
are vanishing perturbations about $x$ and $y$, respectively, governed by the linearization of
system Eq. (\ref{Eq:System}).
Whereas this classical $\lambda_{max}$ is defined in long term limit $(t \rightarrow \infty)$, one can also define \cite{Abarbanel:1991} a \emph{short-term} version,  $\lambda^H_{max}$,  by considering a finite time interval $H$ ($H=50$ kyr will be considered in this article):
\begin{equation}
\lambda^H_{max} =
\lim_{|\delta\bold{Z}(0)|\rightarrow 0}~
\frac{1}{H} \ln
\frac{|\delta\bold{Z}(H)|}
{|\delta\bold{Z}(0)|}
\label{Eq:lyap_local}
\end{equation}

While $\lambda_{max}$  gives the average or long-term stability information, $\lambda^H_{max}$ can tell us about the behaviour of nearby trajectories within a short time interval $H$. For example, $\lambda_{max}<0$ does not necessarily imply $\lambda^H_{max}<0$ for some suitably chosen $H$. The definition (\ref{Eq:lyap_local}) will be useful in studying the robustness of generalized synchronization in Sec.~\ref{Sec:Robustness}.

For computing $ \lambda_{max}$,
as the system Eq. (\ref{Eq:System}) 
and the Jacobian have an analytical form,
tangent space methods \cite{Kantz:2004} can be used;
 technical details are given in the Appendix \ref{Annex:Details_Lyapunov}.

\paragraph{Influence of the parameters $\gamma$ and $T_{ULC}$  \\[0.25cm]}

\begin{figure*}
\centering

\mbox{
\subfigure[{Largest Lyapunov Exponent
$\lambda_{max}$ [kyr$^{-1}$]
for the 41 kyr periodic forcing case.
The region with $\lambda_{max} <0$ corresponds to synchronization; we recognize its underlying Arnol'd tongue structure.
The bifurcation boundaries of the Arnol'd tongues obtained 
with the more accurate numerical continuation method AUTO
 are superimposed, for validation purposes (white curve), and match perfectly.}]
{\includegraphics[width=0.485\textwidth]
{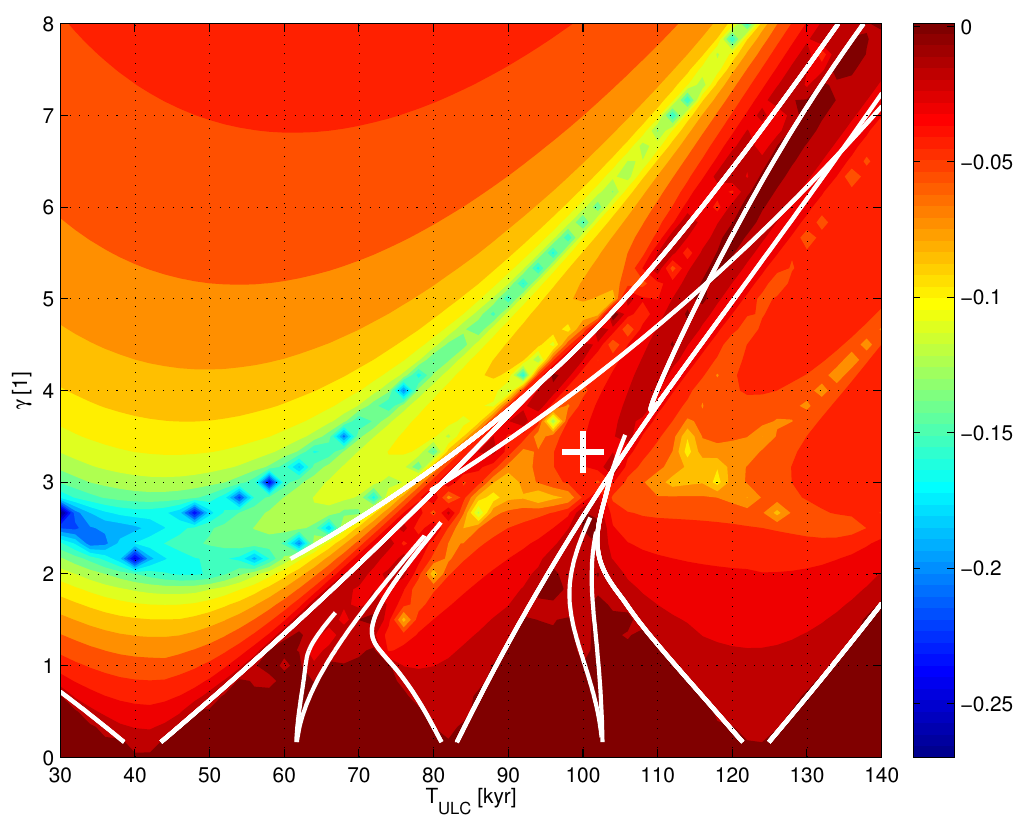}
\label{Fig:Fig_bifurc_b}}\quad
\subfigure[{Largest Lyapunov Exponent
$\lambda_{max}$ [kyr$^{-1}$]
in the case of the quasiperiodic insolation forcing given by Eq.~(\ref{Eq:Insolation_35_terms}): the broad region of synchronization is a typical signature of the quasiperiodic insolation.
The region of synchronization appears to be in one single piece with some indications of well-separated tongues (mode locking) at small $\gamma$.}]
{\includegraphics[width=0.485\textwidth]{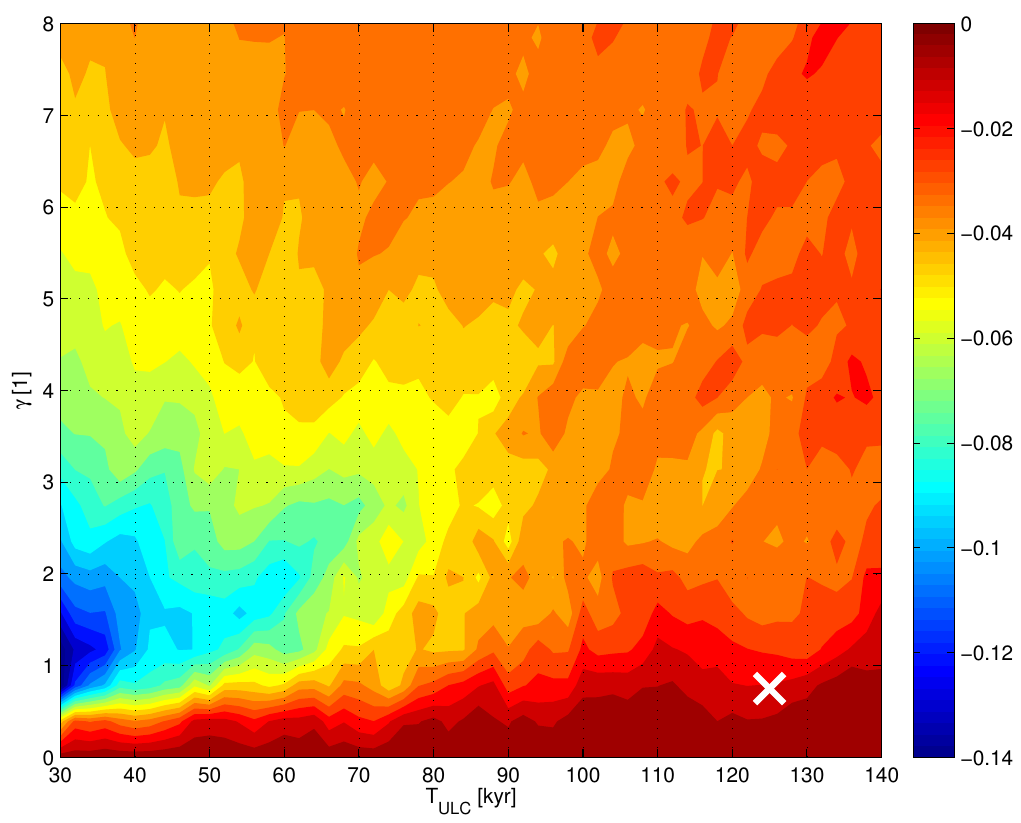}
\label{Fig:Fig_bifurc_c}}}

\mbox{
\subfigure[Numerical estimate of the
number of attracting trajectories $N$  for the 41 kyr periodic forcing case.
In practice, small positive values correspond to synchronization onto a few attracting trajectories, while high values indicate no synchronization.
Now the region inside the synchronization tongues is colored in function of $N$.
For the tongue corresponding to a frequency-locking $n:1$, we have $n$ attracting trajectories.
The bifurcation boundaries of the Arnol'd tongues obtained 
with the more accurate numerical continuation method AUTO
 are superimposed, for validation purposes (black curve), and match perfectly.
]
{\includegraphics[width=0.485\textwidth]{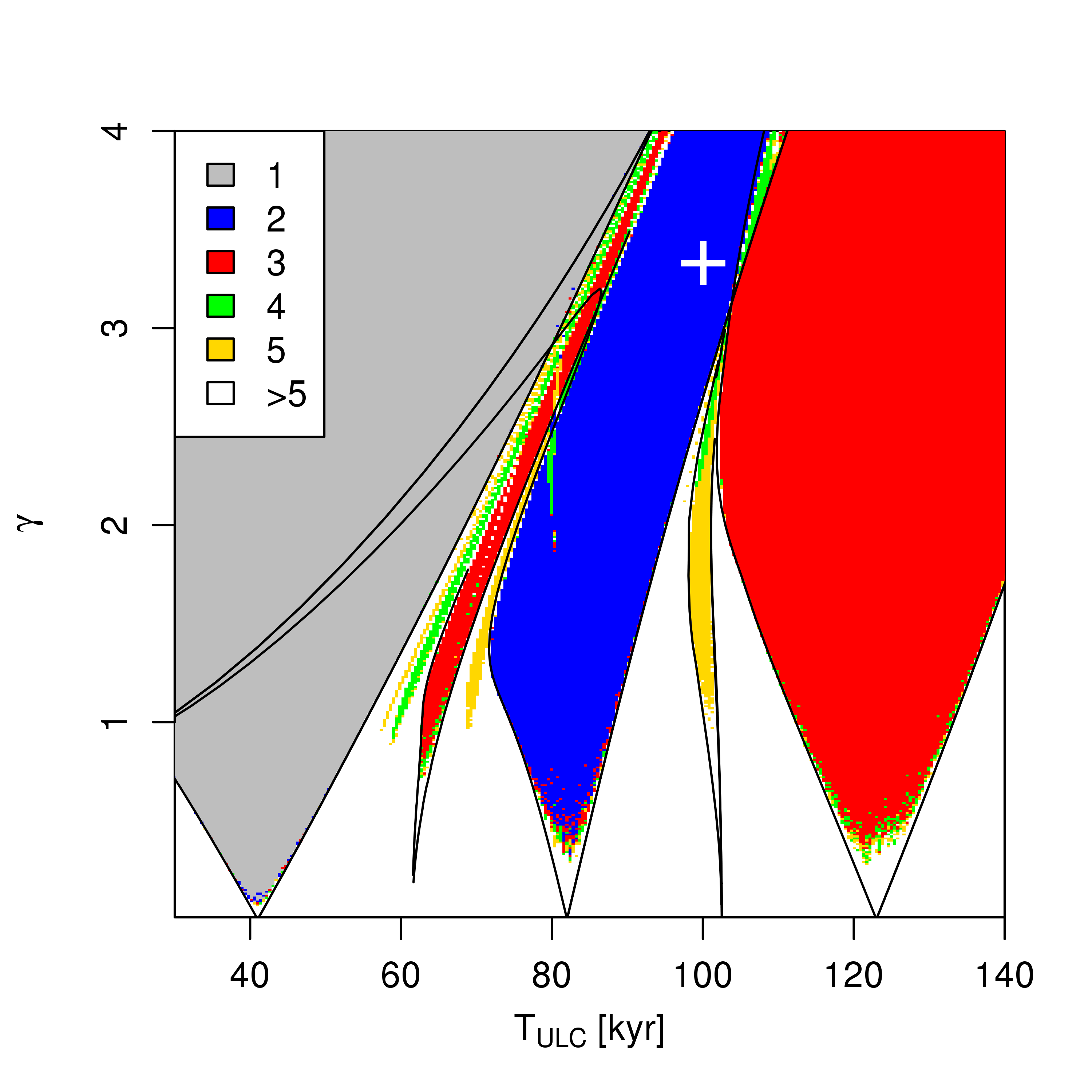}
\label{Fig:Fig_bifurc_Insol_41k}} \quad
\subfigure[{Numerical estimate of the
number of attracting trajectories $N$  in the case of the quasiperiodic insolation forcing given by Eq.~(\ref{Eq:Insolation_35_terms}): the structure is now much more complex, consisting of intermingled series of Arnol'd tongues.
The region with one attracting trajectory, corresponding to {\it unique or monostable generalized synchronization}, is the largest. However, there are also parameter sets with $N = 2, 3$ or even more attracting trajectories.
}]
{\includegraphics[width=0.485\textwidth]{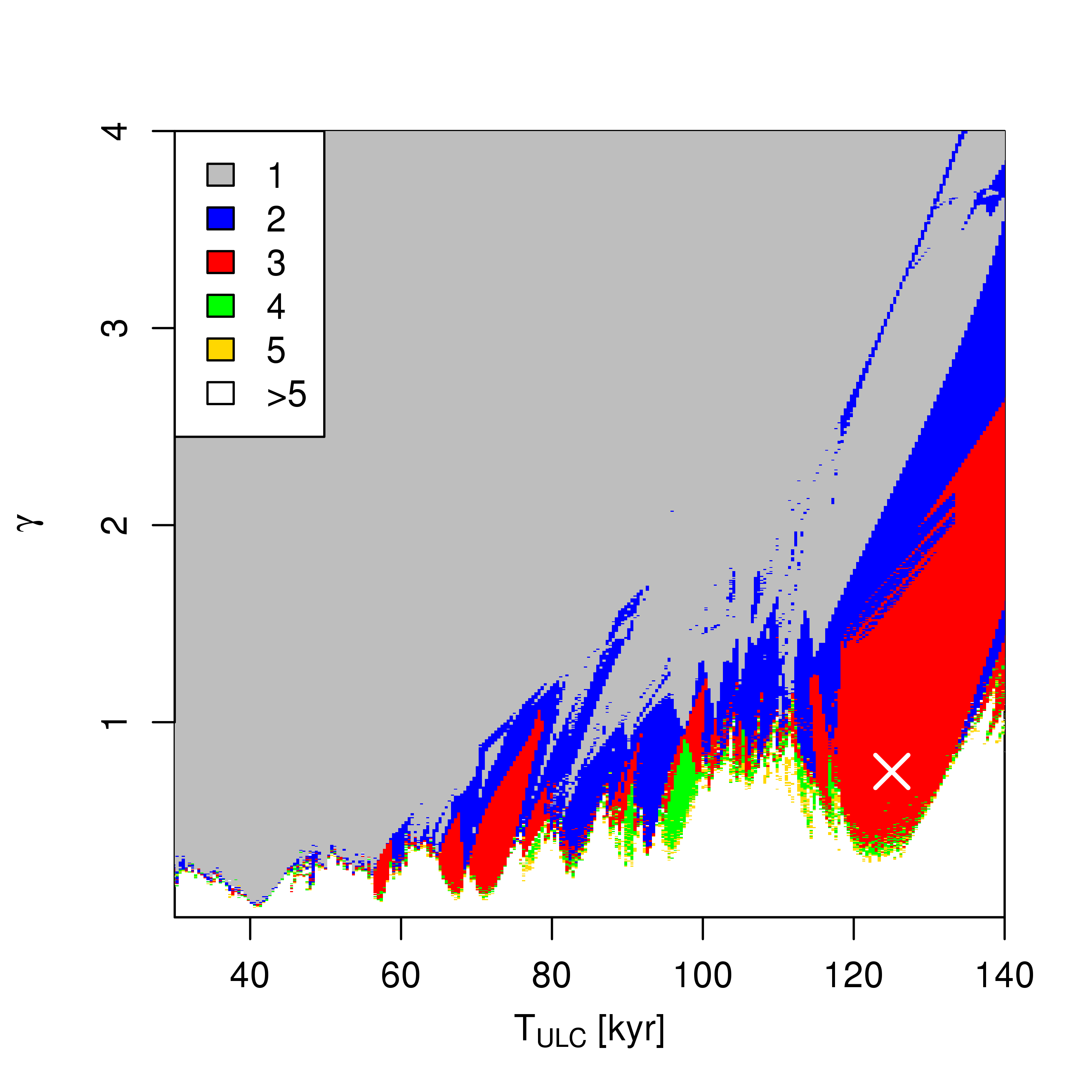}
\label{Fig:Fig_bifurc_Insol}}}

\parbox{\textwidth}{
\caption{Detection of  synchronization
by two different methods: the largest Lyapunov Exponent ($\lambda_{max}<0$) 
(\emph{top}) and the
numerical estimate of the
number of attracting trajectories $N$ by a clustering technique (\emph{bottom}), 
and for two different types of  forcing:
a 41 kyr purely periodic forcing (\emph{left}) and the quasiperiodic insolation forcing given by Eq.~(\ref{Eq:Insolation_35_terms}) (\emph{right}).
These  diagrams
show for which values of the parameters \{$T_{ULC},\gamma$\}  synchronization of the climate system occurs for the  ice age model Eq. (\ref{Eq:System}) with $\alpha = 11.11$, $\beta=0.25$, and $\tau = 35.09$.
The 
symbol (+) refers to the specific climatic attracting trajectories illustrated in Fig. \ref{Fig:Insol_41k} in the 41 kyr periodic forcing case, for which $N=2$. The 
symbol ($\times$) refers to the specific climatic attracting trajectories illustrated in Fig. \ref{Fig:Insol} in the quasiperiodic insolation forcing case,
for which $N=3$.
}
\label{Fig:Fig_bifurc}}
\end{figure*}

The two particular types of synchronization illustrated in Fig. \ref{Fig:Insol_41k} and \ref{Fig:Insol}
have been obtained for a fixed value of the amplitude  $\gamma$ of the external forcing
and of the natural period of the unforced paleoclimatic system $T_{ULC}$.
Now, we are equipped to achieve  a much broader view of the dynamics by
performing a parametric study on these two parameters.
The quantity plotted in Figs.~\ref{Fig:Fig_bifurc_b} and ~\ref{Fig:Fig_bifurc_c} is the largest
Lyapunov Exponent at 3 Myr, far from the transient behaviour so that $\lambda^{H=3~\text{Myr}}_{max}$
may already been considered as a good approximation of  $\lambda_{max} $.

\paragraph{41 kyr periodic forcing \\[0.25cm]}

Fig.~\ref{Fig:Fig_bifurc_b} corresponds to the case of the 41 kyr periodic forcing
($T_F=41$ kyr).
The synchronization region ($\lambda_{max} <0$) is composed of  several V-shape regions, called {\it Arnol'd tongues} (phase- or frequency-locking), originating at 1, 2, 3, etc. times the forcing period $T_F$. These regions correspond to 1:1, 2:1, 3:1 frequency-locking zones (3:2 and 5:2 can also be guessed). Periodic solutions are found within these regions which originate generally speaking at $T_{ULC} = (m/n) \; T_F$.
No synchronization is possible when $\gamma$ is zero  but synchronization may occur  already for infinitesimally small $\gamma$.
Then, for increasing $\gamma$, the synchronization region  widens and synchronization becomes more stable
up to an optimum value $\gamma^*$ of the forcing. When $\gamma > \gamma^* $, the synchronization becomes less and less effective, because at large $\gamma$ the system is too much steered away from its natural dynamics;
it may even be driven into chaos at yet higher forcing amplitude \cite{Mettin:1993} but this case is beyond our focus.


In order to perform an accurate validation of the synchronization region given by the LLE ($\lambda_{max} <0$) method, we computed the main Arnol'd tongues boundaries with the more accurate numerical continuation methods such as AUTO  \cite{Doedel:2009}.
The case of periodic forcing $[F(t) = \sin(\omega t)]$ with $\beta = 0$ has been already extensively studied in the literature, analytically assuming some approximations \cite[p.70--75]{Guckenheimer:1983}, and using numerical algorithms for pseudo-arc length continuation \cite{Mettin:1993}.
The usual approach extends the original nonautonomous system by additional differential equations for the forcing so that the system becomes autonomous, and then explores the $(\omega, \gamma)$  parameter space. Note that the asymmetry introduced here with the parameter $\beta$ adds slightly more complexity and induces additional features to the diagrams documented in these papers.
In this way, we computed  Arnol'd tongue boundaries as saddle-node of limit cycle  bifurcations for the extended system with $\alpha = 11.11$ and $\beta = 0.25$.

Superposition of LLE calculations and bifurcation boundaries in Fig.~\ref{Fig:Fig_bifurc_b} shows that the synchronization
regions obtained with two different techniques match perfectly.
This is a confirmation that the method based on the LLE works fine and we will be able to use it
for the case of the quasiperiodic insolation forcing.
Note that  bifurcation boundaries are also drawn in Fig.~\ref{Fig:Fig_bifurc_Insol_41k}
in order to stress the correspondence  with yet another method of detecting synchronization
that will be discussed in Sec.~\ref{Sec:Clustering}.

\paragraph{Astronomical quasiperiodic forcing \\[0.25cm]}

For the case of the quasiperiodic insolation forcing (Fig.~\ref{Fig:Fig_bifurc_c}),
the region of synchronization appears to be in one single piece with some indications of well-separated tongues (mode locking) at small $\gamma$.
In other words,
whatever the value of the natural period $T_{ULC}$ of the paleoclimatic system,
it has a higher probability of being synchronized onto the insolation forcing.
Of course, for very low values of $\gamma$, there is still no synchronization.
Note that
$\gamma$ has been downscaled by $\sigma_{insol} \approx 5$ (compared to  $\sigma_{sine} = \sqrt 2/2$ for the 41 kyr periodic forcing),
in order to keep a realistic range comparison.

\section{Non uniqueness: multistability and basins of attraction}
\label{Sec:Clustering}

The detection of synchronization using the LLE ($\lambda_{max}<0$)
gives only an Yes/No-type of information, without making any distinction between different tongues as this would require
information about multistability. For example,  Fig. \ref{Fig:Fig_bifurc_c}
indicates synchronization for the parameter settings marked with the symbol '$\times$'' but gives no information
about the number of attracting trajectories (we know that there are three different attracting trajectories from Fig. \ref{Fig:Insol}).
To explore the problem of multistable synchronization, we propose a {\it clustering method} that
not only allows us to detect synchronization, but additionally provides information about the
\emph{number of attracting trajectories} denoted here with $N$.

\paragraph{Multistability Analysis:  numerical estimate of the number of attracting trajectories $N$ by a clustering technique\\[0.25cm]}

Consider the case of the quasiperiodic insolation forcing with the three attracting trajectories,
i.e., $N=3$ (Fig.~\ref{Fig:Insol}).
Although $N$ can often be easily assessed visually, we want to automatically detect and count the number of $AT$s.
As a matter of fact,
$N$ can be easily estimated in the following way. Fix a time $t$ that defines a two-dimensional  $(x,y)$-section in the $(x,y,t)$ phase space. Then start with a grid of initial conditions at some time $t_0<t$ and take $t-t_0$ sufficiently large so that all the initial conditions converge to the attracting trajectories at $t$. Since each $AT$ is represented by a point on the $(x,y)$-section,  the problem of counting attracting trajectories reduces to a simple clustering problem.
We designed a suitable automatic cluster detection algorithm that counts the number of clusters to obtain an estimate of $N$.
For example, Fig.~\ref{Fig:Insol_section} shows\footnote{See also Fig.~\ref{Fig:Fig_basin_IC_t_20b} for a detailed view.} the $(x,y)$-section of the three-dimensional $(x,y,t)$ phase space at $t=550$ kyr, given 70 initial conditions at $t_0=0$.
The 70 trajectories converge onto three (highly concentrated) clusters corresponding to the three attracting trajectories.

The idea
of using clustering analysis for paleoclimatic dynamics comes from the natural fact
that clustering  is another way of looking at generalized synchronization where negative LLE makes the trajectories cluster more efficiently.
This provides another insightful viewpoint on the problem of identification of number of synchronized solutions of the paleoclimatic system: the more stable the synchronization, the more efficient formation of clusters.

\begin{figure}
\begin{center}
\noindent\includegraphics[trim = 4cm 8cm 4cm 10cm, clip, width=0.48\textwidth]{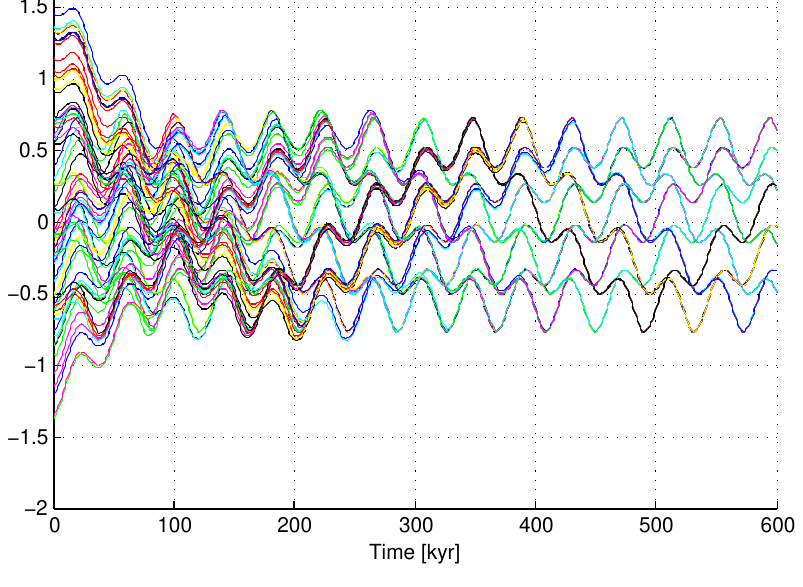}
\caption{Illustration of the importance of choosing the time interval $t-t_0$  large enough so that the  transient behaviour is gone, in order to have sufficiently well formed clusters.
Here, choosing $t=550$ kyr ensures that the eight clusters are already formed, starting from $t_0=0$.
}
\label{Fig:Fig_eps_0_05_8_clusters}
\end{center}
\end{figure}

Two important aspects of cluster analysis\footnote{Cluster analysis or \emph{clustering} is the assignment of a set of observations into subsets (called clusters) so that observations in the same cluster are similar  \emph{in some sense}. This is a common technique for statistical data analysis used in many fields for countless applications. There exists many types of clustering, along with several methods, among which: hierarchical, partitional, spectral, kernel PCA (principal component analysis), k-means, c-means and QT clustering algorithms.} have to be considered  to avoid risks of mis-identification of clusters:
\begin{itemize}
\item
the notion of a cluster is based on the \emph{threshold distance $d_T$}\footnote{This threshold distance $d_T$ appears in any computation related to clustering analysis, or analogically Recurrence Plots (RP) analysis in complex networks, for determining neighbours \cite{Marwan:2009,Donges:2009}.} that has to be carefully chosen.
If $d_T$ is chosen too large, there will be just one cluster including all points;
if it is too small, no clusters will form with more than one point.
\item in order to have sufficiently well formed clusters, the time interval $t-t_0$ must be chosen large enough so that the  transient behaviour is gone; an illustration of the convergence is given in
Fig.~\ref{Fig:Fig_eps_0_05_8_clusters}.
\end{itemize}

Depending on the type and amplitude of the forcing $\gamma$,
we can have potentially a whole range of possible number of attracting trajectories $N$, ranging from
\emph{one} \cite{Tziperman:2006}, to
\emph{a few} (two in the 41 kyr periodic forcing example in Figs.~\ref{Fig:Insol_41k} 
and \ref{Fig:Insol_41k_section},
or three in the quasiperiodic insolation forcing example in Figs.~\ref{Fig:Insol} and \ref{Fig:Insol_section}).
When  no forcing is considered (Figs.~\ref{Fig:No_Insol} and \ref{Fig:No_Insol_section}), or there is forcing but no synchronization occurs, we find no clusters at all. This means that there are as many points in the  $(x,y)$-section at time $t$ as initial conditions at time $t_0$. Clearly, it is difficult to numerically distinguish between  no synchronization and a large number of attracting trajectories ($N \gg 1$). Therefore, we restrict ourselves to just six different regions in Fig.~\ref{Fig:Fig_bifurc}, where we use white to indicate when there are none or more than five attracting trajectories.

Now, we apply the numerical clustering analysis in the case of the periodic forcing (Fig.~\ref{Fig:Fig_bifurc_Insol_41k}) and of the quasiperiodic forcing (Fig.~\ref{Fig:Fig_bifurc_Insol}).
We set $t=0$ and  consider a grid of 49 initial conditions covering $x \in [-2.2;2.2]$ and $y \in [-2.2;2.2]$ at the initial time $t_0 = - 40 T_F$ for the periodic forcing ($T_F$ is the period of the forcing), and  $t_0 =-1600$ kyr  for the astronomical forcing. Two points in the $(x,y)$-section are estimated to belong to a different cluster if their Euclidean distance is greater than 0.1.

\paragraph{41 kyr periodic forcing \\[0.25cm]}

\begin{figure}
\begin{center}
\includegraphics[width=0.48\textwidth]{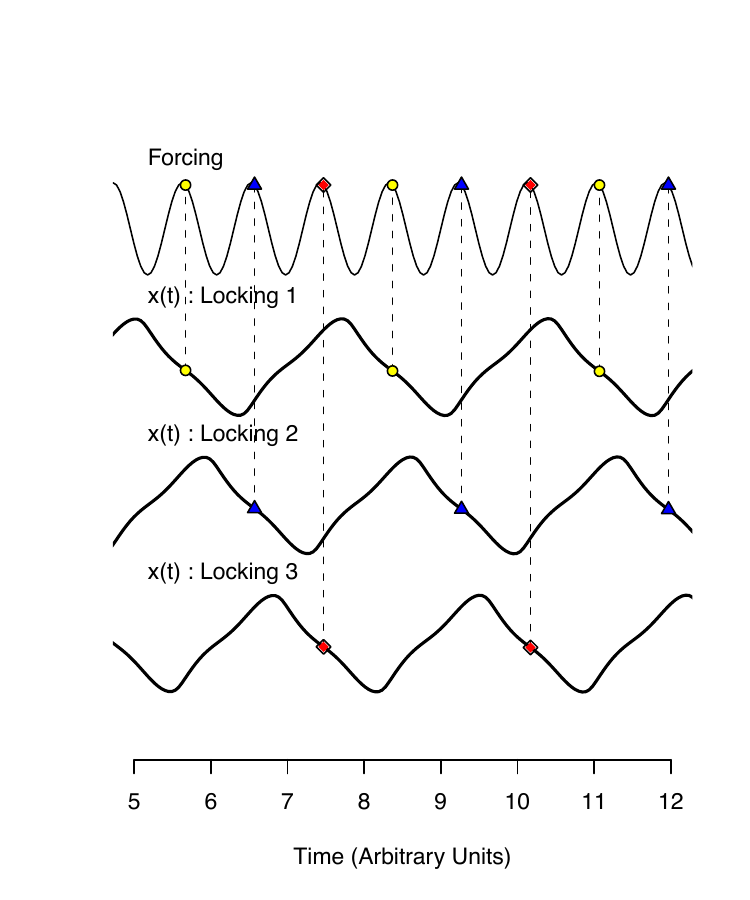}
\caption{
Illustration of the three possible synchronized solutions ($N=3$) existing for the  frequency-locking 3:1 on a purely periodic forcing, in the time series format. The response can be locked on one of the periods of the forcing.
}
\label{Fig:Fig_Locking_3_ways_without_Legend}
\end{center}
\end{figure}

An illustration of the three possible synchronized solutions ($N=3$) existing for the  3:1 frequency-locking on a periodic forcing
is given in Fig. \ref{Fig:Fig_Locking_3_ways_without_Legend}, where the response can be locked on one of the periods of the forcing.
More generally, $N$ corresponds to the number of forcing cycles associated with the synchronization regime ($N=1$ for 1:1; $N=2$ for 2:1; $N=3$ for 3:2, 3:1, etc.)\footnote{This statement relies  on the system invariance with respect to a time-shift of one forcing period (\cite{Tziperman:2006} show a very nice illustration of this point).
}.
The resulting pattern of different $N$  is, as expected, in agreement with the bifurcation diagram (Fig.~\ref{Fig:Fig_bifurc_Insol_41k}). For example, $N = 3$ in the 3:1 tongue. This method allows one to visualize the 4:3 ($N = 4$) and even the 5:4 ($N = 5$) tongues to the left of the 3:2 tongue. It is also seen that $N$ is generally larger where different synchronization regimes co-exist; this is the case between the 2:1 and 3:1 regimes around $\gamma = 4$.

\paragraph{Astronomical quasiperiodic forcing \\[0.25cm]}

Fig.~\ref{Fig:Fig_bifurc_Insol} shows that  synchronization occurs for most parameter configurations. The region with one attracting trajectory ($N = 1$), corresponding to {\it unique or monostable generalized synchronization} \cite{Rulkov:1995}, is the largest. However, there are also parameter sets with $N = 2, 3$ or even more attracting trajectories. They indicate {\it multistable generalized synchronization} where  different possible stable relationships (\ref{Eq:synch}) between the forcing and the oscillator response coexist.

\begin{figure}
\begin{center}
\noindent\includegraphics[width=0.485\textwidth]{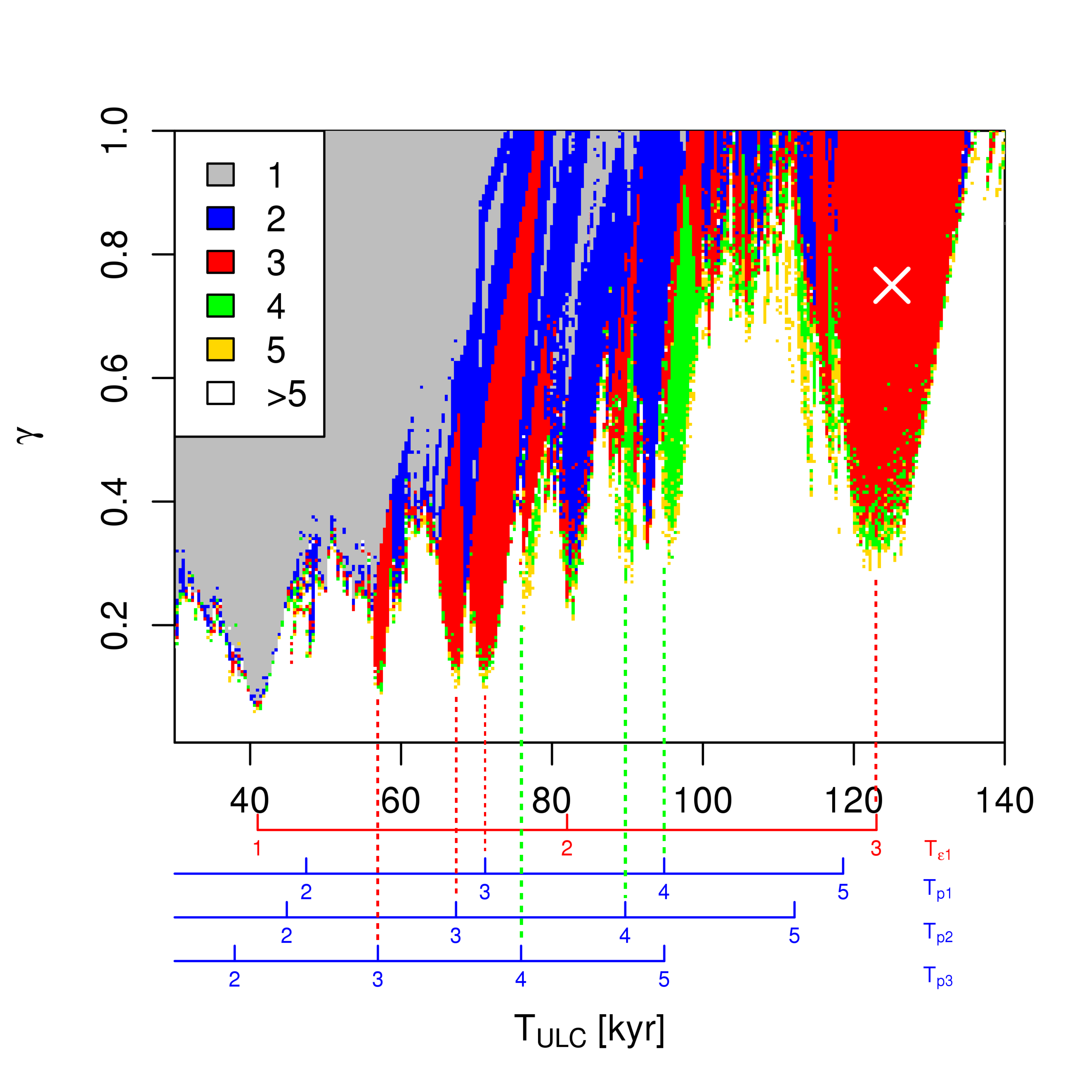}
\caption{A careful examination of the lower part (low $\gamma$) of Fig.~\ref{Fig:Fig_bifurc_Insol}
allows an insightful physical interpretation.The underlying structure of the intermingled series of Arnol'd tongues can be understood:
it is a mixing of several Arnol'd tongues series corresponding to each of the main components  of the
spectrum of the astronomical forcing, yielding a resulting pattern which looks like misaligned combs.
Subscales based on the period of the four main components 
($\epsilon_1$, $p_1$, $p_2$ and $p_3$) of the insolation highlighted in
Fig.~\ref{Fig:Fig_spectre_insol} allows a clear recognition of each individual Arnol'd tongues series.
}
\label{Fig:Cardinality_beta025_alpha_11_astro_zoom}
\end{center}
\end{figure}

A closer view in the lower values of $\gamma$ is given in Fig.~\ref{Fig:Cardinality_beta025_alpha_11_astro_zoom}, which allows an  insightful physical interpretation.

Three tongues with $N = 1, 2, 3$ are rooted at $T_{ULC}/T_{\epsilon_1}$ = 1, 2 and 3, respectively, suggesting a synchronization on the main obliquity component of the astronomical forcing of the same nature as synchronization on a periodic forcing. A series of other synchronization tongues with $N  > 1$ appear; they correspond to 2:1 ($N = 2$), 3:1 ($N = 3$), 4:1 ($N = 4$) and even 5:1 ($N = 5$) synchronization on the three leading components of precession, denoted $p_1$, $p_2$ and $p_3$. Consequently, the richness of the astronomical forcing effectively widens the parameter range for which synchronization occurs, compared to a periodic forcing.
The phenomenon may be understood intuitively: just as you are more likely to tune on some radio station if you are surrounded by a dozen of free FM emitters, the system is more likely to synchronize on the rich astronomical forcing than on a periodic forcing.
Synchronization with $N=1$ or $2$, found for larger $\gamma$, can be interpreted as a form of combined synchronization on both obliquity and precession.

It is crucial to appreciate that synchronized solutions are not periodic and that, unlike the periodic forcing case,  different synchronized solutions for a given set of parameters are not time-shifted versions of each other. The idea that different synchronized solutions co-exist is of practical relevance for paleoclimate theory. Namely, the set of parameters used to obtain the fit to the paleoclimatic records shown in Fig.~\ref{Fig:Fit_VDP} give two distinct solutions at $t = 0$ when started from a grid of initial conditions at $t_0 = -700$ kyr. Sensitivity studies show that  the choice of $t-t_0$ is sometimes important for estimating correctly $N$. However, tests with  $t-t_0$ as large as 200 Myr of astronomical time suggest that several attracting trajectories may co-exist at the asymptotic limit of  $t_0\rightarrow -\infty$.

\begin{figure}
\begin{center}
\noindent\includegraphics[width=0.485\textwidth]{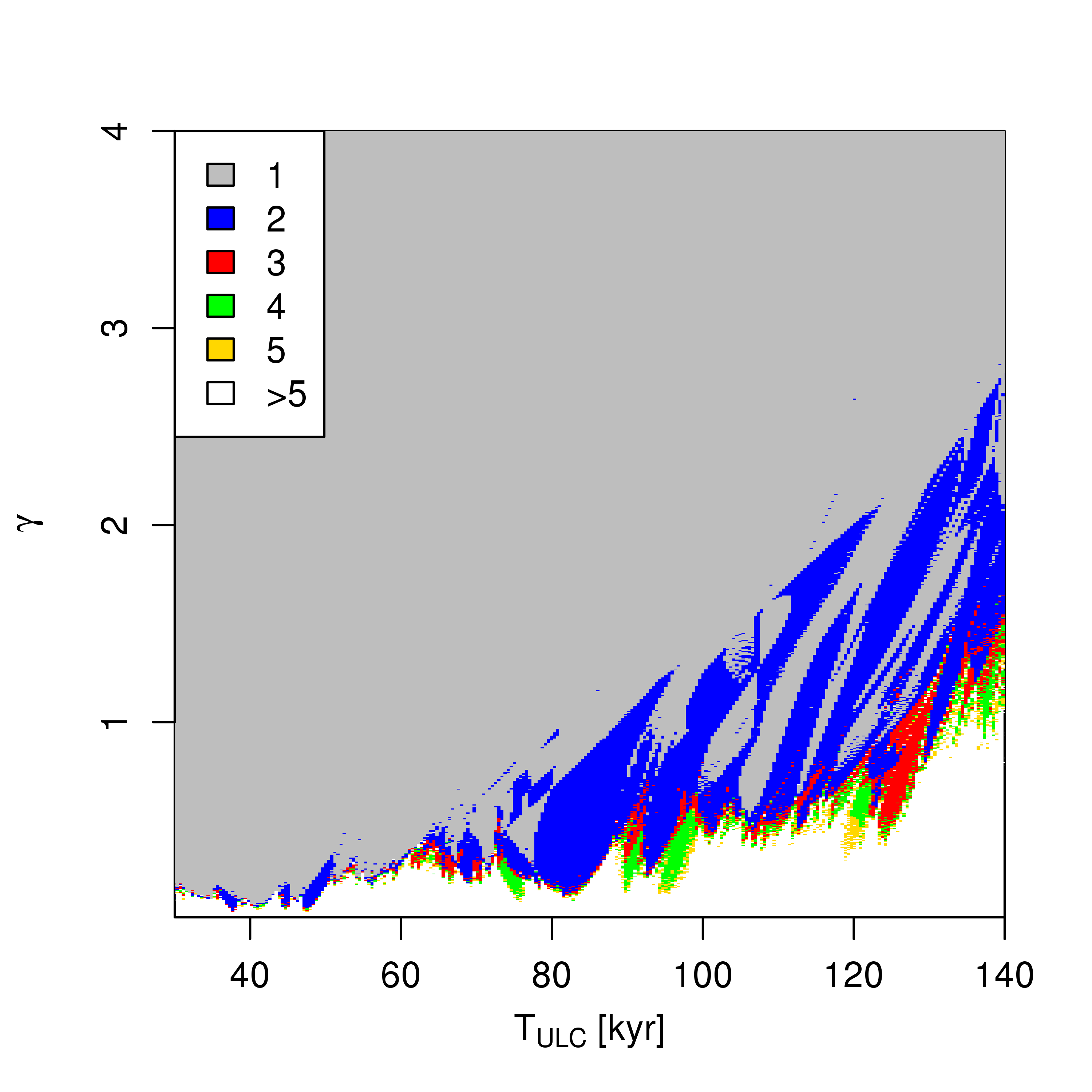}
\caption{Same as Fig. \ref{Fig:Fig_bifurc_Insol}, but now with $\beta = 0.6$ instead of $\beta=0.25$.
The particular values of $N$ change, as the pattern of intermingled tongue series is different,
but the main conclusion about the multistability remains.}
\label{Fig:Cardinality_beta06_alpha_11_astro}
\end{center}
\end{figure}

A similar numerical clustering analysis plot for $\beta = 0.6$ is shown in Fig.~\ref{Fig:Cardinality_beta06_alpha_11_astro} in order to give an idea of the effect of this parameter
(a more detailed analysis is performed in Sect.~\ref{Sec:beta}).
The main conclusion about the multistability remains, but the particular values of $N$ change, as the intermingled tongue series are different.

\paragraph{Evolving shape of the basins of attraction \\[0.25cm]}

Each $AT_i ~(i=1 \dots N)$ has its own \emph{basin of attraction}\footnote{
A more formal definition of the \emph{basin of attraction}
for  nonautonomous dynamical systems is given  in \cite{Kloeden:2000,Langa:2002}.
}
\cite{Barnes:1997}, that is defined as
the set of all initial conditions in the $(x,y,t)$ phase space that converge to that $AT_i$ as time tends to infinity.
For our nonautonomous system Eq.~(\ref{Eq:System}), we can study basins of attraction in the  $(x,y)$-section for different but fixed values of initial time $t_0$, and observe how they vary with $t_0$.
A given initial condition
at time $t_0$ lies in the basin of attraction of $AT_i$ if it approaches  $AT_i$ as time tends to infinity.
Technical details about the computation of the basins of attraction
by use of the specific classification algorithm developed (see Fig. \ref{Fig:Fig_basin_IC_t_20b})
are given in the Appendix \ref{Annex:Details_Basins}.
Basins of attraction are of major importance because they provide the information about global or nonlinear stability of synchronization. If we care about predictability, basin boundaries indicate when a change in the attracting climatic history is likely.

\begin{figure*}
\begin{center}
\includegraphics[width=\textwidth]{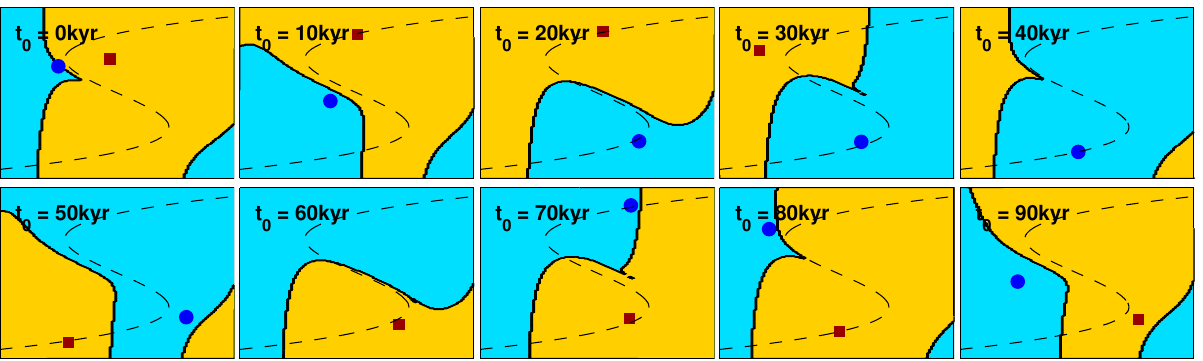}
\caption{Evolving shape of the two basins of attraction for system Eq.~(\ref{Eq:System}), with $\alpha = 11.11$, $\beta=0.25$, $\gamma = 3.33$ and $\tau = 35.09$. Case of a 41 kyr periodic forcing : the pattern is 82 kyr periodic (frequency-locking 2:1). The function $\Phi'(y) = y^3 / 3 -y = x$, corresponding to the slow manifold, is also shown (dashed curve). The attracting trajectories (symbol) sometimes lie  close to the boundary of their own basin of attraction.
}
\label{Fig:Strip_41k}
\end{center}
\end{figure*}
\begin{figure*}
\begin{center}
\includegraphics[width=\textwidth]{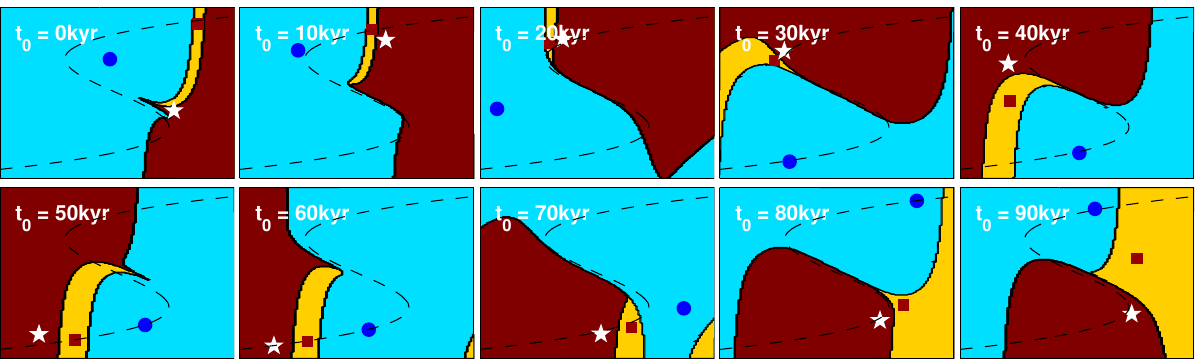}
\caption{Evolving shape of the three basins of attraction for system Eq.~(\ref{Eq:System}), with $\alpha = 11.11$, $\beta=0.25$, $\gamma = 0.12$ and $\tau = 35.09$. Case of the quasiperiodic insolation given by Eq.~(\ref{Eq:Insolation_35_terms}): the pattern is quasiperiodic, constituting the specific signature of the insolation.
The function $\Phi'(y) = y^3 / 3 -y = x$, corresponding to the slow manifold, is also shown (dashed curve). The attracting trajectories (symbol) sometimes lie  close to the boundary of their own basin of attraction.
}
\label{Fig:Strip_insol}
\end{center}
\end{figure*}

The evolving shape of the basins of attraction is shown in Figs. \ref{Fig:Strip_41k} and \ref{Fig:Strip_insol}, for the case of 41 kyr periodic forcing and the quasiperiodic astronomical forcing, respectively. The evolution is shown as a comic strip,
where each subfigure has an $x-$axis ranging  from -1.5 to 1.5, and a $y-$axis ranging from -2.5 to 2.5 (the axis labels have been removed for a better readability).

In the case of a periodic forcing (two basins), the pattern repeats itself periodically (compare
the $t_0=0$ kyr to $t_0=40$ kyr, and to $t_0=80$ kyr subfigures) in Fig.~\ref{Fig:Strip_41k}.
However, in the case of the quasiperiodic forcing (three basins),
the pattern is much more intricate and seems not to repeat itself  for the time horizon considered here.

The ratio between the area of a basin of attraction and the considered area of the phase space can be
interpreted as a probability to converge to the corresponding attracting trajectory when starting
from a randomly chosen  initial condition.
In the case of the periodic forcing, the two  $AT$s are roughly equally likely for all $t_0$
as could be guessed from Fig.~\ref{Fig:Insol_41k}.
However, this is not  the case for the quasiperiodic forcing where the probability to reach the same attracting trajectory
may vary significantly in time. For example,
the yellow basin is rather small at $t_0=0$ kyr but becomes much larger at a later time $t_0=90$ kyr.

In the multistable regime, if an $AT_i$ happens to lie sufficiently close to its basin boundary,
then small perturbations could make the  climate  jump to another (coexisting) $AT_{j\neq i}$, reducing predictability.
This phenomenon is illustrated  in Sec.~\ref{Sec:Robustness}.

\section{Influence of the symmetry-breaking parameter $\beta$}
\label{Sec:beta}

\begin{figure*}
\centering
\mbox{
\subfigure[]
{\includegraphics[width=0.485\textwidth]{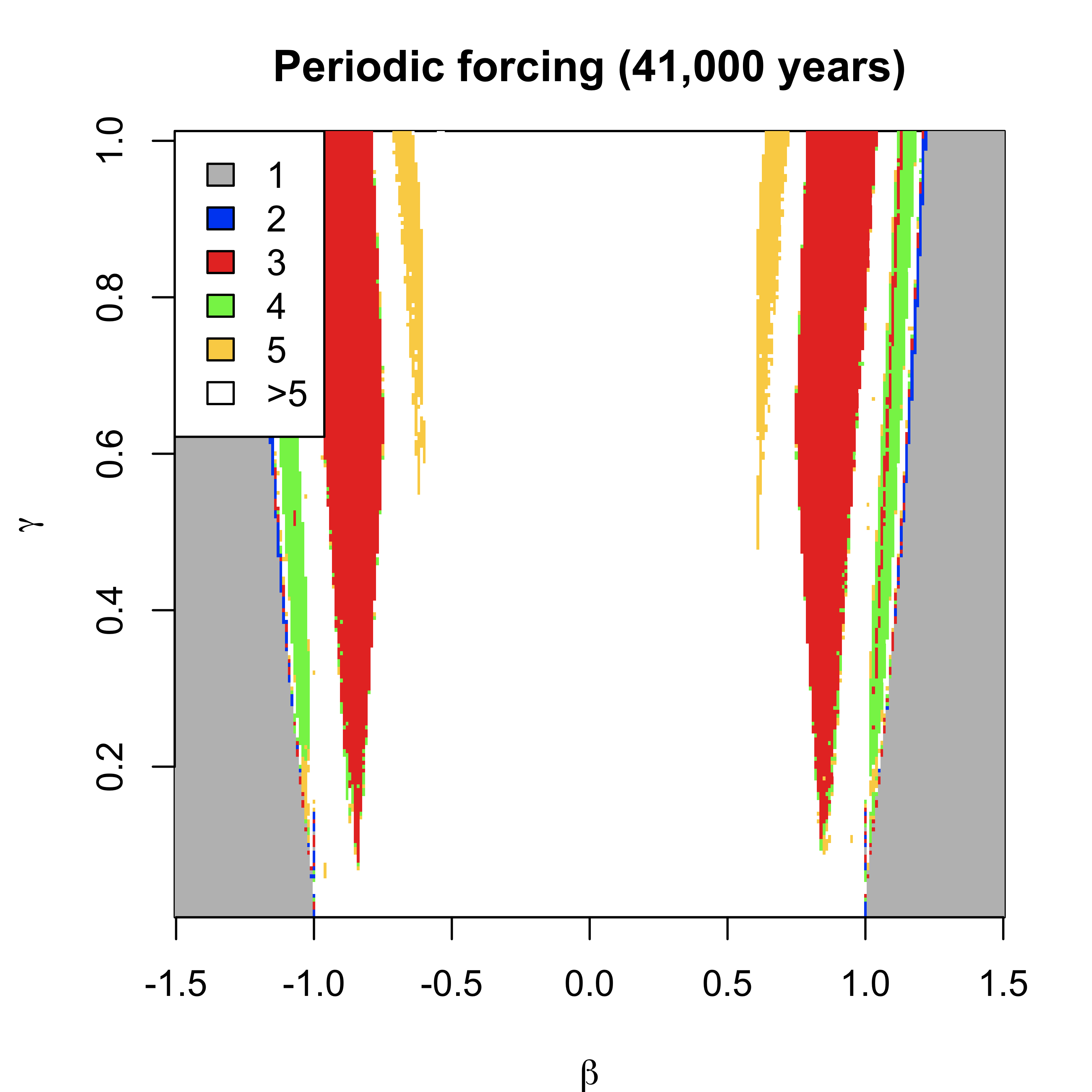}
\label{Fig:Cardinality_beta_gamma}} \quad
\subfigure[]
{\includegraphics[width=0.485\textwidth]{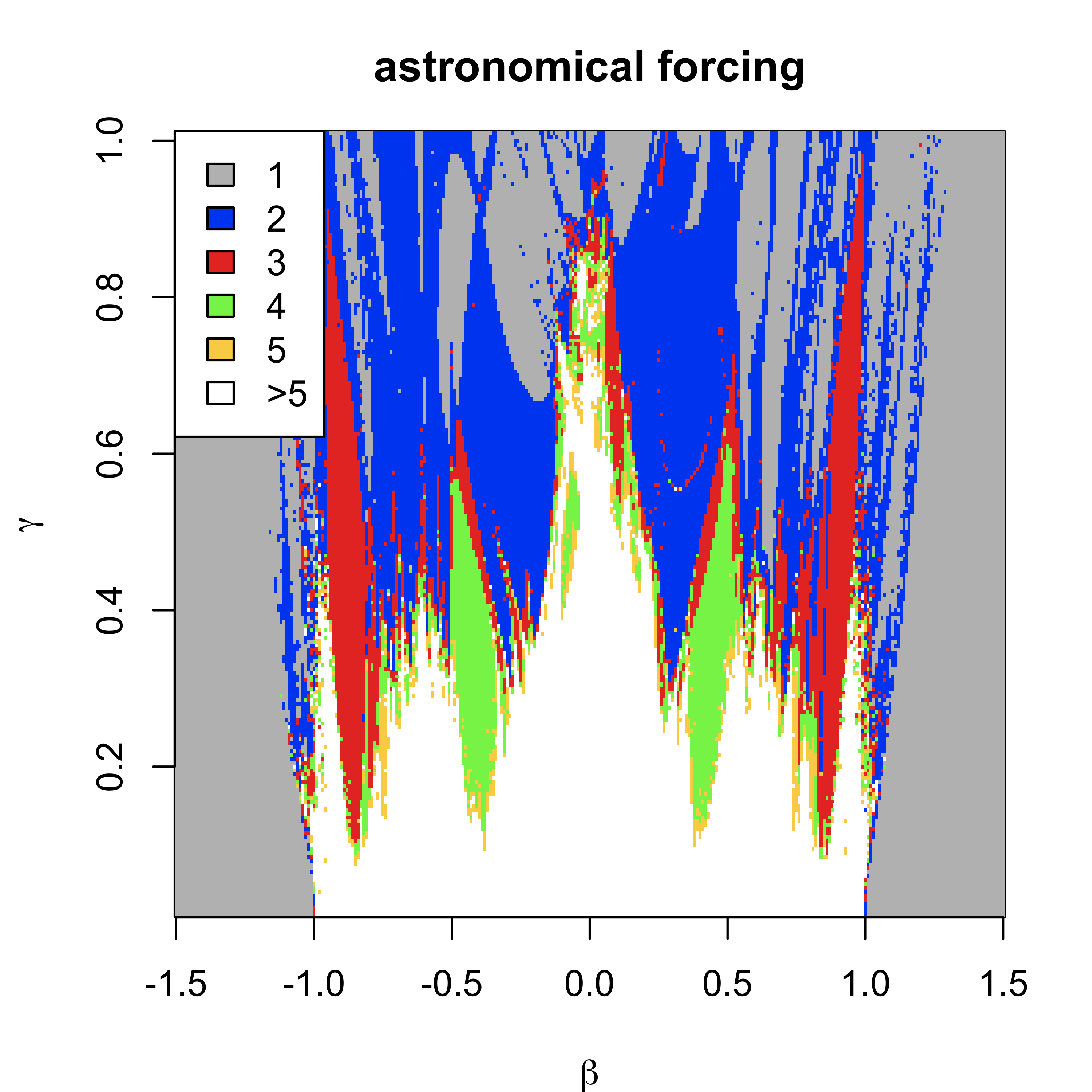}
\label{Fig:Cardinality_beta_gamma_astro}}}
\parbox{\textwidth}{
\caption{Numerical estimate of the number of attracting trajectories $N$ plotted as a function of the asymmetry parameter $\beta$ and the amplitude of the forcing $\gamma$ for the ice age model Eq. (\ref{Eq:System}) with $\alpha = 11.11$, $\beta=0.25$, and $\tau = 35.09$, assuming (\emph{left}) a 41 kyr periodic forcing and (\emph{right}) the quasiperiodic astronomical forcing given by Eq. (\ref{Eq:Insolation_35_terms}). Synchronization occurs in a larger area of the parameter space
in response to the astronomical forcing than in response to the periodic forcing,
and the pattern of different $N$ for larger $\gamma$    is much more richer and complex; this 
structure emerges from the interaction with different harmonics of the astronomical spectrum and their beatnotes.}
\label{Fig:Cardinality_beta_gamma_both}
}
\end{figure*}

As the parameter $\beta$ controls the asymmetry of the glaciation/deglaciation saw-tooth structure (a higher value of $\beta$ leads to an enhanced asymmetry), it is useful to
investigate its effect.
We have already indicated in Fig. \ref{Fig:Cardinality_beta06_alpha_11_astro} that multistability depends on  $\beta$ in the case of the quasiperiodic insolation forcing. A more systematic approach encompassing the whole range of
$\beta$ is shown in Figs. \ref{Fig:Cardinality_beta_gamma} and \ref{Fig:Cardinality_beta_gamma_astro}, for the case of 41 kyr periodic forcing and the quasiperiodic astronomical forcing case, respectively.

First consider the 41 kyr periodic forcing (Fig.~\ref{Fig:Cardinality_beta_gamma}).
To understand this Figure, recall that the unforced oscillator (i.e., $\gamma = 0$) has a stable fixed point for $|\beta| > 1$ and  a stable limit cycle for $|\beta| < 1$.

The system responds almost \emph{linearly}
 to the forcing when $|\beta|$ is sufficiently large. This explains regions of unique synchronization ($N = 1$) where only one climate response is possible. The system becomes \emph{excitable}  when $|\beta|$ is just slightly greater than one. If the forcing is large enough it will excite oscillations. In this case, $N$ is equal to the number of initial conditions if synchronization is lost, or to a smaller number if synchronization occurs.

Consider now the interval $-1 < \beta < 1$. For this, keep in mind (i) that the period of the unforced oscillation varies by almost a factor of two within the range $0 < |\beta| < 1$, and (ii) that synchronization requires some relation between the period of the unforced oscillations and  the forcing period. Consequently, synchronization on the periodic forcing occurs  only for fairly narrow ranges of $\beta$ that are symmetric around zero. The figure reminds us of Arnol'd tongues. The main synchronization regimes detected here correspond to 4:1, 3:1 and 5:2 frequency-locking. Outside these synchronization regimes, the system fails to converge to a sufficiently small set of attracting trajectories, meaning that the forcing is not as an efficient pacemaker.

Finally, compare this situation with that obtained with the astronomical forcing (Fig.~\ref{Fig:Cardinality_beta_gamma_astro}).
Synchronization now occurs in a larger area of the parameter space. Whereas the structure  of the periodic forcing is preserved as long as the forcing amplitude is low enough, there is a much more richer and complex  pattern of different $N$ for larger $\gamma$.
This pattern emerges from the interaction with different harmonics and their beatnotes.

\section{Robustness of  synchronization}
\label{Sec:Robustness}

{\it Robustness} or {\it reliability} of synchronization can be studied in terms of two properties of an attracting climatic trajectory. Local stability analysis based on the short-term LLE ($\lambda^H_{max}$) provides information about the short-term local convergence towards the $AT$. For example, a temporary loss of local stability indicated by $\lambda^H_{max}>0$ will cause a temporary loss of synchrony and divergence from the $AT$ even though the trajectory is stable on average ($\lambda_{max}<0$). Global stability analysis based on the geometry of basins of attraction for different $AT$s provides information about the system response to external perturbations such as random  fluctuations. 
For example, an external perturbation may push a climatic trajectory outside of the basin of attraction of the $AT$. 
Robustness and uniqueness of synchronization are closely linked in the sense that the global stability is a factor only when there are coexisting attracting trajectories. Robustness is compromised most when a temporary loss of local stability coalesces with a weakening of the global stability. We will now briefly discuss these two effects that could restrict the prediction horizon for the evolution of climatic trajectories.



\paragraph{Temporary desynchronization via loss of local stability \\[0.25cm]}

\begin{figure}
\begin{center}
\includegraphics[width=0.48\textwidth]{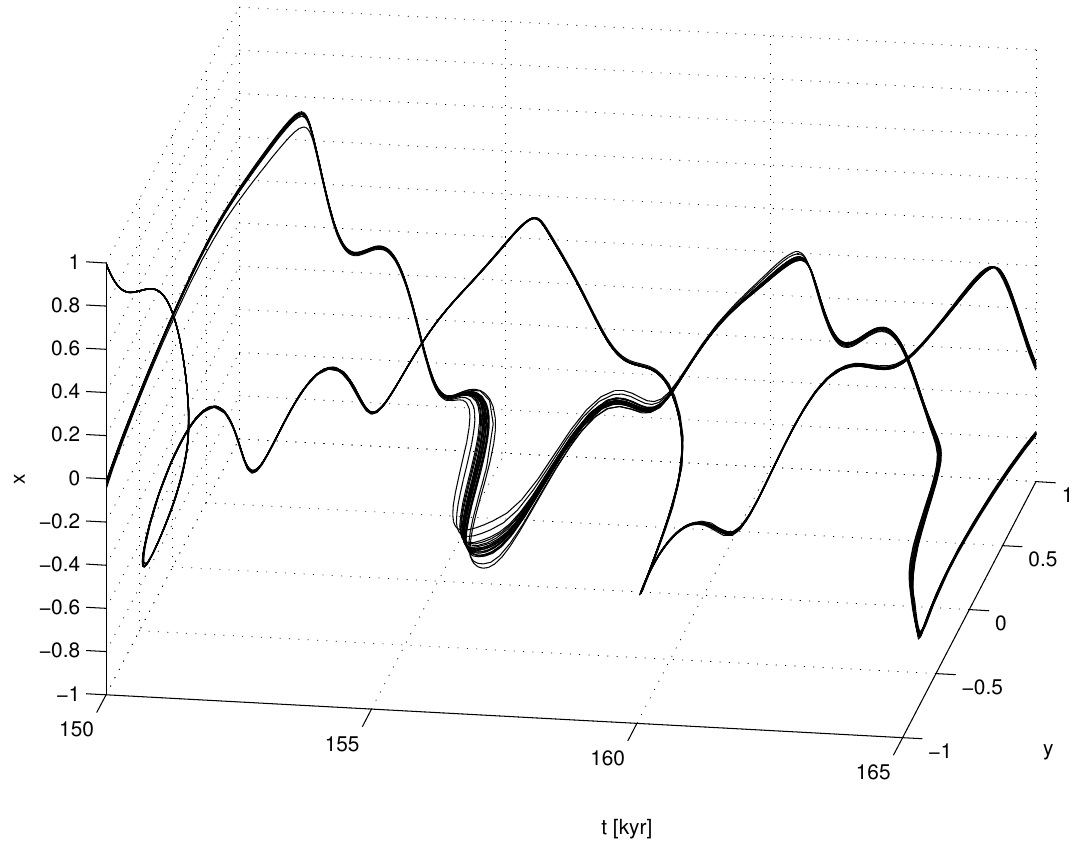}
\caption{The temporary divergence of nearby climatic trajectories reveals short-term local instabilities (e.g. around $t \approx157$ kyr).
This illustration was obtained by considering a set of 50 random initial conditions at $t_0=0$ within $x \times y \in  [-1,1] \times [-1,1] $, and by adding some  noise of small amplitude in the system at every time step in order to trigger the instabilities. The model used is Eq. (\ref{Eq:System}) with $\alpha = 11.11$, $\beta=0.25$, $\gamma = 0.033$, and $\tau = 33.33$.
}
\label{Fig:Spread}
\end{center}
\end{figure}

Some additional experiments made in our paleoclimatic framework
reveal another strange behaviour in the system
 Eq. (\ref{Eq:System}), that can be deduced from a careful inspection of Fig. \ref{Fig:Spread}.
In the presence of small additive noise, we notice that nearby trajectories could diverge for some time, like those around $t \approx 157$ kyr.

\begin{figure}
\begin{center}
\includegraphics[width=0.50\textwidth]{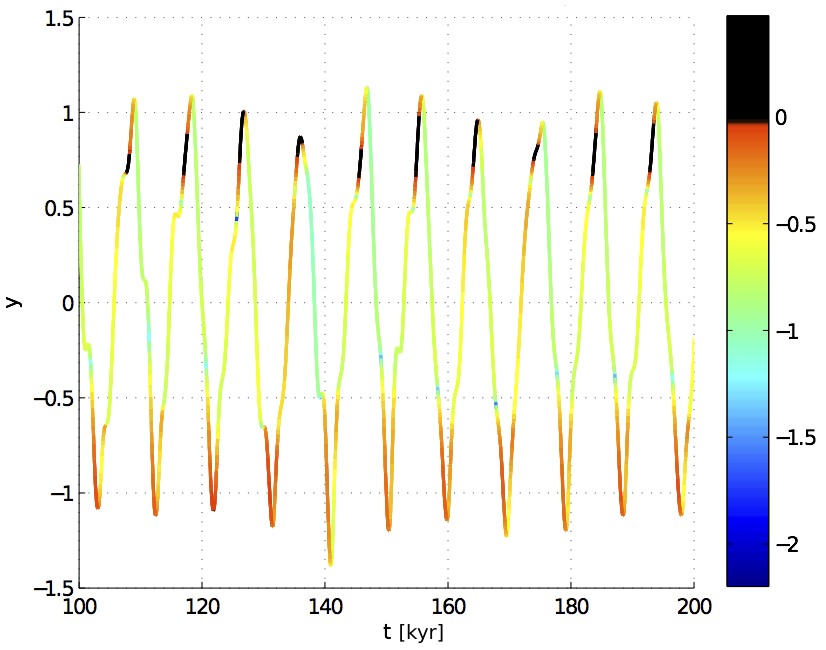}
\caption{Short-term largest Lyapunov Exponent $\lambda^{H=50~\text{kyr}}_{max}$ of one of the two attracting trajectories of system Eq. (\ref{Eq:System}), subject to insolation forcing given by Eq.~(\ref{Eq:Insolation_35_terms}), and for the same parameters as in Fig.~\ref{Fig:Spread}.
The attracting climatic trajectory has been coloured according to the values of $\lambda^{H=50~\text{kyr}}_{max}$, revealing temporary desynchronizations; this explains the divergence of nearby trajectories found in Fig.~\ref{Fig:Spread}, e.g. around $t \approx157$ kyr.}
\label{Fig:MT_Lyap}
\end{center}
\end{figure}

Such temporary divergence is  similar to desynchronization bursts \cite{Rulkov:1995}
and strongly suggests to investigate the evolving sign of the short-term LLE $\lambda^H_{max}$ along the attracting climatic trajectory.
We computed  $\lambda^{H=50~\text{kyr}}_{max}$ along one of the two attracting trajectories
of system Eq. (\ref{Eq:System}), subject to insolation forcing given by Eq.~(\ref{Eq:Insolation_35_terms}).
The result is shown in Fig. \ref{Fig:MT_Lyap}, where
the attracting climatic trajectory has been coloured according to the values of $\lambda^{H=50~\text{kyr}}_{max}$.
Although
the system is synchronized on a long term ($ \lambda_{max} = -0.2$ kyr$^{-1}$, see Fig.~\ref{Fig:LT_Lyap}), we see here that there exist episodes with positive values of the short-term\footnote{
For completeness, the very short-term (\emph{instantaneous}) stability has also been
investigated (see Fig. \ref{Fig:Link_Lyap_Jac} in the Appendix \ref{Annex:Instantaneous_instability}), as a limit case $H \rightarrow 0$,
but it is less physically relevant within the paleoclimatic context.
}  LLE $\lambda^H_{max}$, revealing temporary desynchronizations \cite{Wieczorek:2009}.  This explains the divergence of nearby trajectories found in Fig.~\ref{Fig:Spread}.

These results remained unchanged with respect to the most important parameters of the model.
For example,
our main conclusions about the stability remain qualitatively valid, even for  different values of $\alpha$ (like $\alpha = 100$), or with a different type of potential
($\Phi'_5(y) = (y+1.7)(y+1.58)(y+0.8)(y)(y-0.5)$),
even if the shape and size of the limit cycle and the boundaries  of the basins of attraction are of course  different.
The effect of the insolation function $F(t)$ has also been checked:
we compared the attracting trajectories for the insolation given by Eq.~(\ref{Eq:Insolation_35_terms}), compared to those for the insolation given by  \cite{Laskar:2004}.
As these insolation functions are very similar, the results are also very similar, and no difference was noticed.

At first glance, it may appear that these episodes of temporary divergence are not relevant to the robustness of synchronization because  climatic trajectories converge back the the attracting trajectory on a long term. However, other effects may be present that could strongly amplify such temporary divergence. They are identified below.


\paragraph{Sensitivity to perturbations: preliminary results  \\[0.25cm]}

Consider again  Fig.~\ref{Fig:Strip_insol} showing  $(x,y)$-sections with coexisting attracting trajectories 
in the case of the quasiperiodic insolation, 
and their basins of attraction for different values of $t_0$. Suppose now that the system is subject to additive fluctuations (for example, these may represent volcanic eruptions).
Under certain conditions, such external perturbations may cause a displacement of the trajectory to a different basin of attraction, causing a jump\footnote{
In the periodic forcing case the phenomenon of jumping from one attracting trajectory to another  in response to a perturbation is called  a \emph{phase slip} \cite[p.238]{Pikovsky:2001}.
} to another attracting trajectory. 

\begin{figure}
\begin{center}
\includegraphics[width=0.48\textwidth]{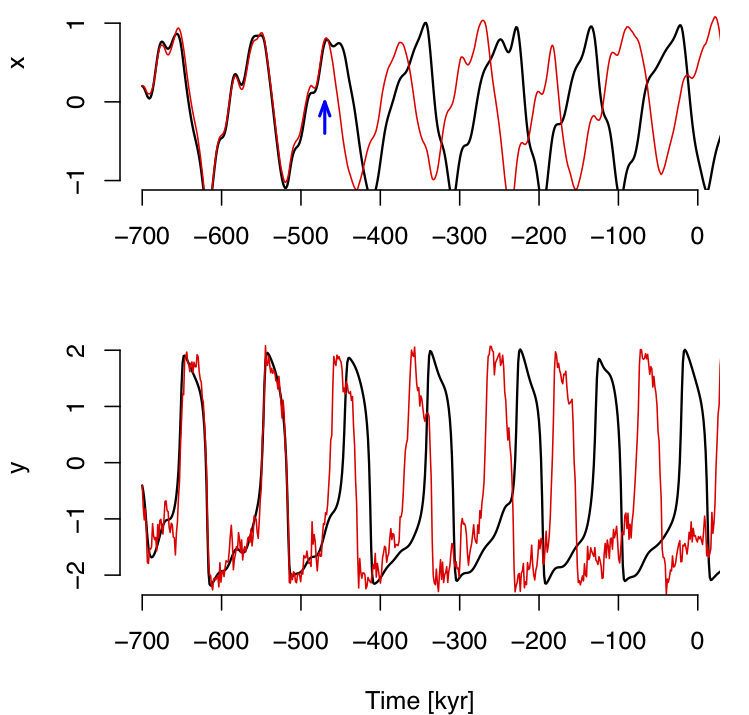}
\caption{Sensitivity of the climatic system to perturbations. The same solution of the ice ages model Eq. (\ref{Eq:System})  as in Fig.~\ref{Fig:Fit_VDP} is plotted (black) along with one sample trajectory of the same system (red), but with additive fluctuations added to the fast variable: $dy = - \tau^{-1}(\alpha(\Phi'(y) - x))dt + b~dW$, with $b = 0.5 \sqrt{\omega_{\epsilon_1}}$ and $W$ a Wiener-process. The arrow shows the time of the  jumping from one to the other climatic attracting trajectory, which reduces the predictability of the timing of the glaciations.
}
\label{Fig:stochastic}
\end{center}
\end{figure}

As a further illustration of this idea we show in Fig.~\ref{Fig:stochastic} two attracting trajectories (in the time series format) that coexist for the same system parameters as those used for the fit of Fig.~\ref{Fig:Fit_VDP}, but with additive fluctuations added to the fast variable (see legend for details).   
A jump from one trajectory to another at around -475 kyr (arrow) may clearly been identified.
This shows that a climatic trajectory is robust against fluctuations if it stays away from the basin boundary but its robustness can weaken significantly due to the weakening of the global stability near the basin boundary.

We conclude that externally triggered jumps between coexisting attracting climatic trajectories seem to be most likely when the temporary desynchronization due to the loss of local stability coalesces with the weakening of the global stability due to the proximity to the basin boundary.

\section{Conclusions}
\label{Sec:Conclusions}

We have identified, illustrated, and provided a systematic study of {\it generalized and multistable synchronization} between the climatic glacial/interglacial oscillations and the  astronomical forcing. 
For doing so, a series of appropriate concepts and tools have been developed.
A van der Pol-type relaxation oscillator, designed to reproduce the slow-fast dynamics of the paleoclimatic records, has been used for illustration purposes, but the methodology proposed here may of course be applied to other paleoclimatic models.

To study the uniqueness of synchronization, we  proposed a convenient concept of the number of attracting trajectories in the phase space of the nonautonomous forced system,  each of which corresponds to a synchronized solution. We computed the number of synchronized solutions using a numerical clustering technique,  and  uncovered that in addition to a {\it unique} or {\it monostable} synchronization, there are parameter settings where one finds a {\it nonunique} or {\it multistable} synchronization. At low forcing amplitude we found regions of {\it mode locking} where the system synchronizes on the individual components of the astronomical forcing in a way that is similar to frequency-locking on periodic forcing (Arnol'd tongues), giving rise to coexisting  synchronized solutions. As the forcing amplitude is increased, the combined effects of precession and obliquity restrict the number of possible synchronized solutions. The emerging stability diagram consists of a large region of monostable synchronization mixed with smaller regions of multistable synchronization. A comparison with periodic forcing shows that the system finds it easier to synchronize to quasiperiodic insolation forcing. It is therefore conceivable that the climate system wandered throughout preferential synchronization regimes on obliquity, precession, or combinations of both, as environmental parameters varied throughout the history of the Pleistocene.

The robustness of generalized synchronization was investigated in terms of the key indicators of stability of synchronized
solutions: the  long- and short-term largest Lyapunov exponent (local stability), and the geometry of the basins of attraction
(global stability). We found that even though the synchronized solutions are locally stable on a long term, there exists
episodes where the  short-term largest Lyapunov exponent becomes positive, leading to temporary desynchronizations.
As a result, climatic trajectories could diverge from the synchronized solution for some period of time (50 kyr typically).
Moreover, we computed the evolving shape of the basins of attraction for the coexisting synchronized solutions, and uncovered that
these solutions sometimes approach the basin boundary where they become very susceptible to external perturbations.
As a result, a small perturbation could make the climate jump from one synchronized solution to another, reducing predictability.
Such jumps seem to be most likely when the temporary loss of the local stability coalesces with the proximity to the basin boundary.
In this context, we briefly discussed the effect of stochastic perturbations on the timing of the deglaciations.
We also illustrated the difference between the evolving shape of the basins of attraction for periodic and quasiperiodic
insolation forcing. In the case of the insolation forcing, we obtained an intricate pattern of basins of attraction that
does not appear to repeat itself in time.


%
%

\vspace*{1cm}
{\bf \Large \begin{center} Acknowledgements \end{center} }
\vspace*{1cm}

We are grateful to Guillaume Lenoir for his thorough review of the several versions of the paper.
The original idea of using clustering analysis for automatically identifying the number of stable locking states came to the main author, after presentations and discussions, especially at the 1$^{st}$ ITOP Workshop, held on 23--26 February, 2010 in Marche-en-Famenne, Belgium, and also from the 458. WE-Heraeus-Seminar on 'SYNCLINE 2010: Synchronization in Complex Networks', held on 26--29 May 2010 at the Physikzentrum Bad Honnef (Germany), where
some preliminary results of this research have been presented in a poster \cite{De-Saedeleer:2010}.
The
project is funded by the ERC (European Research Council) starting grant ITOP ('Integrated Theory and Observations of the Pleistocene') under the convention ERC-2009-5tg 239604-ITOP. M. Crucifix is Research Associate with the Belgian National Fund of Scientific Research, and B. De Saedeleer is Post Doctoral Research Assistant with the ITOP Project.
Some
Figures and calculations where made with the R language and the Intel Fortran Compiler.

\vspace{1 cm}
\appendix

\numberwithin{equation}{section} 

{\LARGE\bf Appendix}

\section{Insolation model using 35 terms}
\label{Annex:Insolation}

We give at the Table \ref{Tab:Insolation35}
the numerical values of the $3 \times 35$ terms for computing the insolation $F(t)$ following
Eq.~(\ref{Eq:Insolation_35_terms}) in Sect.~\ref{Sec:Introduction}.

\begin{table}[!htp]
\caption{Coefficients for the insolation model described in Eq.~(\ref{Eq:Insolation_35_terms}) used for the incoming solar radiation anomaly at Summer Solstice and at 65$^\circ$N. The first 15 terms correspond to the obliquity component, while the other 20 terms correspond to the precession.}
\label{Tab:Insolation35}       
\begin{tabular}{rrr}
\hline\noalign{\smallskip}
$\omega_i $ [rad/kyr] & $s_i$ [W/m$^2$] & $c_i$  [W/m$^2$] \\
\noalign{\smallskip}\hline\noalign{\smallskip}
0.153249478547167 &-11.2287376815124 &3.51682075211241 \\
0.158148666238883 &-3.82499371467540 &-0.761851750263805 \\
0.117190147169570 &2.28814805956066 &1.80233702684623 \\
0.155061775112933 &-1.29770081956440 &-0.635152963728496 \\
0.217333905941751 &0.380973541305497 &-1.46301711999210 \\
0.150162587421217 &1.54904176353302 &-0.0883941912769817 \\
0.211709630908568 &-0.810768209286259 &-0.577980646565494 \\
0.156336369673117 &-0.918358442095885 &0.196083726889428 \\
0.148350290855451 &0.256895610735773 &-0.524697312305024 \\
0.206924898030688 &-0.335783913402678 &-0.0194792150128644 \\
0.212525165090383 &0.267659228540196 &0.128915417116900 \\
0.229992875969202 &0.0696189733188958 &0.0746231714061285 \\
0.306498957094334 &0.0247349748169616 &0.0140464395340974 \\
0.311398144786051 &0.0138353727621181 &0.0304736668840422 \\
0.004899187691716 &-0.160479848721994 &0.0594077968934257 \\ \hline
0.264933601588513 &-15.5490493322904 &-9.70406287110532 \\
0.280151350350945 &15.4319556361701 &4.75247271131525 \\
0.331110950251899 &9.0992249352734 &-10.6115244887390 \\
0.328024059125949 &-7.87065384013669 &6.61544246063503 \\
0.326211762560183 &0.813786144754451 &-4.52641408099246 \\
0.269742342439881 &0.0690448504314857 &-3.31639260969558 \\
0.332923246817665 &1.44050770785967 &1.06339286050120 \\
0.371638925683567 &0.925324276580528 &-1.02066758672154 \\
0.275366617473065 &0.997628846513796 &-0.362906496840039 \\
0.323124871434233 &-0.378637986107629 &0.527217891742183 \\
0.259396912994958 &0.339477750517033 &-0.560509461538342 \\
0.324937167999999 &-0.576082669762308 &1.18669572739338 \\
0.334197841377850 &0.346906064369828 &-0.648189701487285 \\
0.274551083291250 &-0.441772417569753 &0.289576210423804 \\
0.418183080135680 &-0.0184884064645011 &0.109632390175297 \\
0.111684123041346 &-0.428006728186239 &0.357006342316690 \\
0.433400828898112 &-0.0049199219454561 &-0.106148873639336 \\
0.126901871803777 &0.257509918217341 &-0.377639794223366 \\
0.336010137943616 &-0.421809264016129 &0.324327509437558 \\
0.177861471704732 &-0.161827722328271 &-0.362683869407858 \\
\noalign{\smallskip}\hline
\end{tabular}
\end{table}


\section{The classical van der Pol oscillator : description and dynamical behaviour}
\label{Annex:vdp}

\paragraph{Classical van der Pol oscillator\\[0.25cm]}

The classical van der Pol oscillator \cite{vanderPol:1926} is very well known, widely used and deeply studied; here the necessary description which is relevant to the scope of this article is given.
For more details the reader is referred to the literature, see e.g. \cite{Strogatz:1994}.

This van der Pol-type oscillator model can be regarded as a special case of the FitzHugh-Nagumo (FHN) model  \cite{Kosmidis:2003}, also known as Bonhoeffer-van der Pol (BVP) model
\cite{Barnes:1997}.

Historically, this model was derived by van der Pol, when he discovered
the existence of oscillations in electrical circuits.
He found that the oscillation period is determined by $\tau^* = RC$ (time constant of relaxation)
in $RC$ circuits, or by $\tau^* = L/R$ in $RL$ circuits ;
hence he named this oscillation as \emph{relaxation oscillations}.
The interesting characteristics of the relaxation oscillation are the slow asymptotic behavior and the sudden discontinuous jump to another value.
This oscillator has therefore been widely used in many fields like
in physical and biological sciences, neurology, seismology, lasers, optoelectronics, etc.

This two-dimensional system exhibits the same basic dynamical features (limit cycle and slow-fast dynamics) necessary to fit the paleoclimatic data, except the asymmetry, for which the system will be corrected for, by the introduction of an additional parameter $\beta$ (see Eq. \ref{Eq:x} in Sect.~\ref{Sec:Introduction}).

This van der Pol oscillator  is non conservative with a nonlinear damping, governed by the following second-order differential equation \cite{Balanov:2009} --- in the case without forcing :
\begin{equation}
\ddot x -  \mu(1-x^2)  \dot x +x = 0
\label{Eq:vdp_ second-order}
\end{equation}
where $ \mu > 0$ is a positive constant proportional to the damping.

The van der Pol oscillator  is a Li\'enard system  \cite{Strogatz:1994},
because $f(x) = -  \mu(1-x^2)$ is an even function and because $g(x) =x$ is an odd function
(see Eq. \ref{Eq:vdp_ second-order}), in the canonical form :
\begin{equation}
\ddot x + f(x)~ \dot x +g(x) = 0
\label{Eq:Lienard}
\end{equation}

Moreover, since this Li\'enard system satisfies additionally the Li\'enard theorem,
it has a unique and stable limit cycle\footnote{
The nowadays still unsolved Hilbert's 16th Problem is related to determining the number and location of limit cycles for an autonomous planar vector field for which both functions are real polynomials of fixed degree.}
in the phase space surrounding the origin.
This is a feature which is absolutely required in order to model the glacial-interglacial oscillations.

By way of the Li\'enard transformation $ y = x -x^3/3 - \dot x /  \mu$,
the second-order Eq. \ref{Eq:vdp_ second-order} can be transformed into an equivalent two-dimensional system of ODE's :
\begin{eqnarray}
\dot x &=&\mu ( x -x^3/3 - y) \label{Eq:vdp_first-order-x}   \\
\dot y &=&  x / \mu \label{Eq:vdp_first-order-y}
\end{eqnarray}
which
can also be expressed under the equivalent  following
form (using the convention $  \varepsilon  = -1 / \mu $, so that $ \varepsilon \rightarrow 0$ means more and more damping):
\begin{eqnarray}
\dot x &=&(y - \Phi'(x))/ \varepsilon \label{Eq:vdp-x} \\
\dot y &=& -\varepsilon x \label{Eq:vdp-y}
\end{eqnarray}
introducing the function
$\Phi'(x)$ :
\begin{equation}
\Phi'(x)  = x^3 / 3 -x
\label{Eq:phi3}
\end{equation}
the associated potential $\Phi(x)$ having the shape of a 2-well potential.
By varying $ \varepsilon $, we can adjust
the respective time scales of the $x$ and $y$ variables;
for $ \varepsilon < 1$, the slow-fast variable is $x$ and the slow variable is $y$.
Note that in the main body of this paper,   $x$ and $y$ are inverted, so that $x$ will represent the slowly varying ice volume  (see Fig. \ref{Fig:vdp_03}), as it is the classical convention for the  dynamical theory of paleoclimates.

When $\mu > 0$, the system exhibits a stable limit cycle, where energy is conserved. Near the origin $x= \dot x = 0$, the system is unstable and energy is gained, and far from the origin the system is damped and energy is lost (while when  $\mu = 0$, there is no damping,
the solution is a pure harmonic signal, and there is conservation of energy everywhere).

There is only one fixed point which is the origin $(x, \dot x) = (0,0)$.

When $x$ is small,
$f(x) \approx - \mu$ (negative damping).
Thus, the fixed point $(x, \dot x) = (0,0)$ is unstable
(an unstable focus when $0 < \mu < 2$, and an unstable node otherwise);
see  \cite{Hilborn:2000}.
On the other hand,
when $x$ is large,
$f(x) \approx + \mu x^2$ (positive damping).
Therefore, the dynamics of the system is expected to be restricted in some area around the fixed point.

When $0 < \mu \ll 1$ (small damping),
the system can be rewritten in order to avoid division by $\mu$.
When $ \mu \gg 1$ (large damping),
the oscillations become less and less symmetric.

When driven, the van der Pol oscillator can lead to synchronization
\cite{Balanov:2009},
but also to deterministic chaos \cite{Ruihong:2008},
depending of the level of the driving force.
Since the paleoclimatic system is driven by the insolation, one could also have synchronization and perhaps chaos -- or at least some path on the routes leading to chaos; i.e. Lyapunov exponents
becoming positive.

This oscillator, or slightly different versions of it (a similar one is the Poincar\'e oscillator \cite{Glass:1994}), has been mathematically largely studied under many aspects: bifurcation structure \cite{Mettin:1993}, fixed points and Arnol'd tongues, chaotic dynamics \cite{Chen:2008, Parlitz:1987},
additive noise \cite{Degli-Esposti-Boschi:2002},
basins of attraction  \cite{Barnes:1997}, etc.
In this article, we  focus on the synchronization, multistability and
predictability properties of a slightly modified version of the classical van der Pol oscillator.

But there is also a big difference of frameworks: usually, a simple periodic forcing is considered,
while in the case of the astronomical forcing,
the forcing has a much more complex form (quasiperiodic), like the one given in Eq.~(\ref{Eq:Insolation_35_terms}).


\paragraph{Time spent on the slow manifold\\[0.25cm]}

In the case of large damping $( \mu \gg 1$, or $\varepsilon \rightarrow 0$),
the oscillations become less and less symmetric,
and significant differences becomes to appear
between $\tau_{slow}$ and $\tau_{fast}$, as defined in Section \ref{Sec:Introduction}.

One can more precisely characterize the slow-fast dynamics of the system
by considering Eqs.  \ref{Eq:vdp-x}--\ref{Eq:vdp-y}.
We see that we have in any case $| \dot y |  \sim \varepsilon $,
which means that $y$ is always slow.

When the system is away from the curve $y = \Phi'(x)$,
we have $(| \dot x|  \sim 1/\varepsilon)  \gg (| \dot y |  \sim \varepsilon)$ :
the vector field is mostly horizontal.
And the travel speed is $v=\sqrt{| \dot x |^2+| \dot y |^2}   \sim 1/\varepsilon \gg 1$,
so that the system moves \emph{quickly} in the horizontal direction
($y$ is the slow variable, and $x$ is the slow-fast variable).
If  $ y > \Phi'(x) $, then $\dot x > 0$, and the trajectory moves clockwise.
See Fig. \ref{Fig:vdp_03}  for an illustration and \cite{Strogatz:1994} for a discussion about the separation between the two time scales.

When the system enters the region\footnote{This region is sometimes called "boundary layer" \cite{Guckenheimer:1983}, by analogy with Fluid Dynamics.} where $ | y - \Phi'(x) |  \sim \varepsilon^2 $,
we then have $(| \dot x|  \sim \varepsilon)  \sim (| \dot y |  \sim \varepsilon)$ :
$| \dot x| $ an $| \dot y| $ are about the same order of magnitude,
and the travel speed is $v  \sim \varepsilon  \rightarrow 0$ :
the system moves slowly along the curve $ y = \Phi'(x)$,
and eventually exits from this region, and so on :
the system has a stable limit cycle.

Hence we have the following relationships:
\begin{eqnarray}
\tau_{slow} \sim \mu \quad &\text{or}&  \quad \tau_{slow} \sim 1 / \varepsilon   \label{Eq:Period_slow}   \\
\tau_{fast} \sim 1/\mu \quad &\text{or}&  \quad \tau_{fast} \sim \varepsilon \quad . \label{Eq:Period_fast}
\end{eqnarray}

The period of the oscillations (or of the unperturbed limit cycle), given by
\begin{equation}
T_{Cycle} =  \tau_{slow} + \tau_{fast} \quad ,
\label{Eq:Period_Limit_Cycle}
\end{equation}
is mainly determined by the time during which the system stays around
the curve $y = \Phi'(x)$,
which can be roughly estimated to be $T \propto  1/\varepsilon$,
as $v  \sim \varepsilon  \rightarrow 0$.
So, the larger the damping (lower $\varepsilon$), the longer the period of the oscillations;
that's why an additional time scaling (factor $\tau$) must sometimes be done in order to keep the 100 kyr period for the limit cycle (see Fig \ref{Fig:vdp_03}).

Let us mention also that analytical expressions have been derived for the amplitude and period of the limit cycle \cite{DAcunto:2006} in both cases ($\mu $ small or big),
and also for the slow manifold equation  \cite{Ginoux:2006}.
There exists also a very fast transition upon variation of a parameter $\beta$, from a small amplitude limit cycle to a large amplitude relaxation cycle, explained by the
so-called canard phenomenon, cycles, and explosion
\cite{Benoit:1981,Guckenheimer:2000,Guckenheimer:2005}.


\section{Technical details about the calculation of the Lyapunov exponents}
\label{Annex:Details_Lyapunov}

We first refer the reader to the seminal papers \cite{Shimada:1979,Benettin:1980fk}.
The methods for computing the Lyapunov (Characteristic) Exponents (LCE) vary depending on the fact that one wishes to achieve the full spectrum of the LCE  \cite{Wolf:1985}, or only the largest one \cite{Rosenstein:1993}.
Analytical derivations of the LCE's of the van der Pol oscillator do also exist \cite{Grasman:2005}.
In this research,
we computed $\lambda_{max}$  using the  standard method
involving a Gram-Schmidt Reorthonormalizaton (GSR) of the 'tangent vectors'
\cite{Shimada:1979,Benettin:1980fk,Wolf:1985};
which is described in the review paper \cite{Ramasubramanian:2000}.

We remind here the fundamental principles. Consider an $n$-dimensional continuous-time dynamical system :
\begin{equation}
\frac{d\vec{Z}}{dt} = \vec{f} (\vec{Z},t)
\label{Eq:dynamical_system}
\end{equation}
where $\vec{Z}$ and $\vec{f}$ are $n$-dimensional vector fields. To determine the LCE's corresponding to some initial condition $\vec{Z}(0)$, we have to find the long term evolution of the axes of an infinitesimal sphere of states around $\vec{Z}$(0).
That is to say that we assume \cite{Ott:2002}
$\delta\vec{Z}(t) = \delta\vec{Z}(0) ~ \text{e}^{\Lambda t} $, with
$\Lambda = diag(\lambda_1, \dots, \lambda_n)$, where $\lambda_i$ are the eigenvalues of the system.
For this, we consider the linearization of Eq. \ref{Eq:dynamical_system}, given by:
\begin{equation}
\frac{d~\delta\vec{Z}}{dt} = \vec{J}~ \delta \vec{Z}
\label{Eq:tangent_map}
\end{equation}
where $ \vec{J}$ is the $n \times n$ Jacobian matrix defined by
$J_{ij} = \partial f_i / \partial Z_j$.
Then, starting from a unit vector $\delta  \vec{Z} (0)$, the original system given by Eq. \ref{Eq:dynamical_system}
is integrated for  $\vec{Z}$  together with the $\delta  \vec{Z} $ tangent system
given by Eq. \ref{Eq:tangent_map}.
The evolution of $\delta  \vec{Z} $ is such that it tends to align with the most unstable direction (the most rapidly growing one).
The choice of the initial vector of the tangent manifold may influence the convergence,
but in practice a spin-up phase can be performed in order to find the good direction.

The largest Lyapunov exponent $\lambda_{max}$ is then defined as :
\begin{equation}
\lambda_{max} =
\lim_{|\delta\bold{Z}(0)|\rightarrow 0}~
\lim_{t \rightarrow \infty}~
\frac{1}{t} \ln
\frac{|\delta\bold{Z}(t)|}
{|\delta\bold{Z}(0)|}
\tag{\ref{Eq:lyap_global}}
\end{equation}

\begin{figure}
\begin{center}
\includegraphics[width=0.40\textwidth]{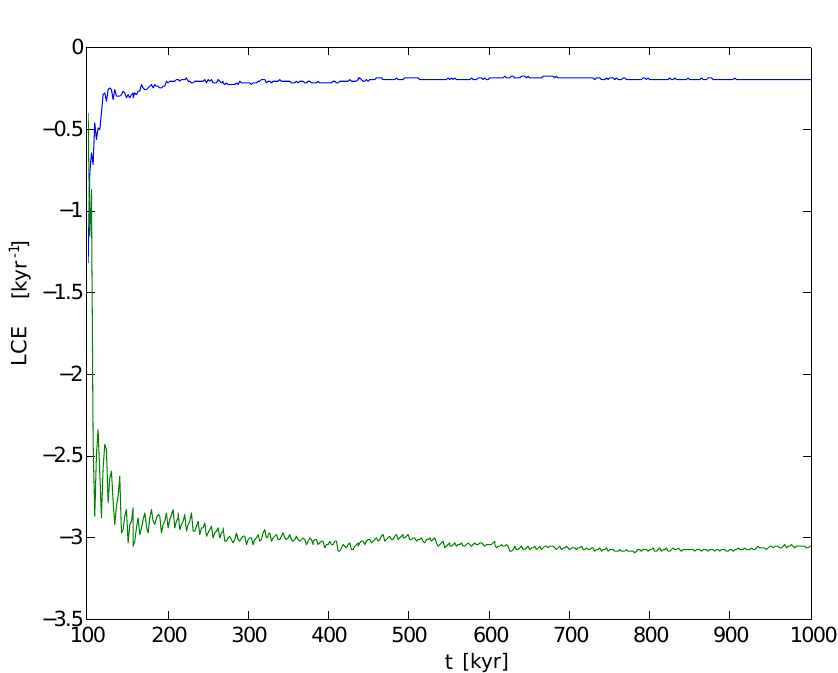}

\caption{Convergence of the two LCE's of
one of the two attracting trajectories of the system Eq. (\ref{Eq:System}), subject to the insolation (Eq.~(\ref{Eq:Insolation_35_terms})); some transient is skipped.
For $\tau = 3.33$, we have  $ \lambda_{max} = -0.2 $ kyr$^{-1} < 0$
(for $\tau = 35.09$, the value must be scaled accordingly, which would give -0.019 kyr$^{-1}$), so
the system is synchronized on a long term.}
\label{Fig:LT_Lyap}
\end{center}
\end{figure}

It is of course impossible in practice
to go to infinity\footnote{The proof of the existence of such a limit has been made by
\cite{Oseledec:1968}.}; so the computation is always
truncated to some finite final time, usually of the order of $10^{4-5}$ times the period of the forcing.
The convergence
can be more rapidly achieved if some  transient behaviour is skipped (like illustrated in Fig.~\ref{Fig:LT_Lyap}),
i.e. we compute the LCE's only when we are quite sure to be on the attracting trajectory.

Note that there exists a whole spectrum of Lyapunov exponents $\lambda_{i}$ ($i=1,2,...,n$), that can be computed with a unit vector basis, and with a renormalization procedure.
Although
we are mostly interested in $\lambda_{max}$ in this article
--- a positive $\lambda_{max}$ is associated to a desynchronization ---,
our subroutine allows to compute all the spectrum of $\lambda_{i}$'s of any system.
For the sake of flexibility, we used a symbolic software, so that the model could be very easily changed, and all functions (like the Jacobian matrix) are automatically derived once the initial system is given.

One
of the standard and popular methods
to compute the Lyapunov spectrum of a dynamical system involves a Gram-Schmidt Reorthonormalizaton (GSR) of the 'tangent vectors'
\cite{Shimada:1979,Benettin:1980fk,Wolf:1985};
differential versions of this have also been formulated.
Our
subroutine includes the GSR procedure, which is required in order to avoid computational overflows, and degeneracy into a single vector.
The frequency of reorthonormalization is not critical; as a rule of thumb, GSR is usually performed on the order of once per characteristic period.
Here, the normalization time step has been chosen optimally, that is to say the largest possible which preserves the accuracy.
Since the GSR never affects the direction of the first vector in a system, this vector tends to seek out the direction in tangent space which is most rapidly growing.

Our subroutine has been  validated by comparing our LCE's to those of \cite{Ramasubramanian:2000} for several systems (driven van der Pol, Lorenz '63, etc.); the order of the accuracy achieved  is 1\permil ~for the Lorenz system\footnote{For the Lorenz system, one must pay attention to the spurious trivial set of LCE's corresponding to the origin \cite{Shimada:1979,Bryant:1990}.}.
One also checked that $\text{div} ~\vec{f} = Tr(\mathcal{J}) = - \sigma -1 - \beta = -13.66 =  \sum_{i=1}^n \lambda_i $
holds.

Coming back to our system of Eqs \ref{Eq:x}--\ref{Eq:y}, for which the Jacobian matrix is:
\begin{equation}
\mathcal{J} = - 
\left[
\begin{array}{ccc}
0 && 1  \\
-\alpha && \alpha \Phi ''(y)
\end{array}
\right] / \tau
=
\mathcal{J} (y,  \alpha  )
\label{Eq:Jacobian}
\end{equation}
we end up to the following results (Fig.~\ref{Fig:LT_Lyap}):
the system Eq. (\ref{Eq:System}) subject to the insolation (Eq.~(\ref{Eq:Insolation_35_terms}))
is synchronized on a long term, since
$ \lambda_{max} = -0.2 $ kyr$^{-1} < 0$.

\paragraph{Some properties of the LCE's ($\lambda_i$) \\[0.25cm]}

The LCE's are very useful in order to characterize the dynamical behaviour of a system; for example,
here are a few interesting properties, valid in the case of the dissipative system Eq. (\ref{Eq:System})
in the autonomous case :
\renewcommand{\labelitemi}{$\bullet$}
\begin{itemize}
\item a dynamical system of dimension $n$ has $n$ LCE's and $n$ eigenvectors \cite{Lichtenberg:1983},
\item $ \sum_{i=1}^n \lambda_i = Tr(\vec{J}) = \text{div} ~\vec{f} $, both being related to the growth of a volume of dimension $n$ of the phase space  \cite{Shimada:1979},
\item $\text{div}(\vec{f}) < 0$ for a dissipative system (hence the system has at least one negative exponent),
\item $\lambda_i = 0$ along a limit cycle (tangent direction of the attractor) \cite{Kantz:2004},
\item at least one LCE vanishes if the trajectory of an attractor does not contain a fixed point \cite{Haken:1983},
\item the nature of the attractor can be described by analyzing the sign of the LCE's,
which gives a qualitative picture of the dynamics \cite{Wolf:1985,Muller:1995},
also function of the dimension (1D to 4D), e.g.
an attractor for a dissipative system with one or more $\lambda_i > 0$ is said to be \emph{strange}
or \emph{chaotic}; if more than one $\lambda_i > 0$, there is \emph{hyperchaoticity} \cite{Rossler:1979}.
\item the existence of $\lambda_i > 0$ is mathematically related to the theory of ergodicity\footnote{An attractor is said \emph{ergodic} (or transitive) if the points fill an entire physically realizable domain; \emph{intransitive} if there are distinct attractors in several closed subdomains \cite{Saltzman:2002},
and a system is \emph{almost intransitive} if it resides for long periods of time in one or another domain that is not fully closed off, so that occasional exits from one domain to another can occur.
} of dynamical systems \cite{Eckmann:1985},
\item local bifurcations can be detected by detecting changes of signs of $\lambda_i$.

\end{itemize}

The LCE spectrum is closely related to the fractional dimension of the associated strange
attractor by the Kaplan-Yorke formula \cite{Kaplan:1979}. There are a number of similar measures of the "strangeness" of strange attractors, like the fractal dimension   \cite{Theiler:1990}, information dimension, box-counting dimension, and the correlation dimension \cite{Tsiganis:1999} and exponent \cite{Grassberger:1983},
which allows to distinguish between deterministic chaos and random noise,
and is computationally easier.

It is also possible to estimate the LCE of a system by analyzing its time series, but
with limited data, or a system subject to non negligible stochastic perturbation, the classical methods may provide incorrect or ambiguous results \cite{Ruelle:1990},
hence require specific methods, like  \cite{McCaffrey:1992, Liu:2005}.

\section{Technical details about the computation of the basins of attraction}
\label{Annex:Details_Basins}

The practical computation of a basin of attraction is done as follows.

Let us for example come back to the Fig. \ref{Fig:Insol}, with the three attracting trajectories $AT_i(x,y,t)$ due to the insolation forcing.

Now, we wonder which initial condition leads ultimately to which of the three $AT_i$. This is the concept of the \emph{basin of attraction}
of a given $AT_i$, classically defined by the set of states that leads to a given $AT_i$?.
Let us more precisely define the basin of attraction of a given $AT_i$ as the locus of all points in the ($x,y$) plane which lead to motion which ultimately converges on that $AT_i$.

We initialize many initial conditions on a fine rectangular grid covering the phase space.
Each initial condition is then integrated forward to see which $AT_i$ its trajectory approached. If the trajectory approached a particular one of the three $AT_i$'s, a dot coloured by the color identifying the $AT_i$ is plotted on the grid.
For doing this, we need to define a target time at which we will do the classification,
and a criteria for the classification.

\begin{figure}
\begin{center}
\includegraphics[width=0.40\textwidth]{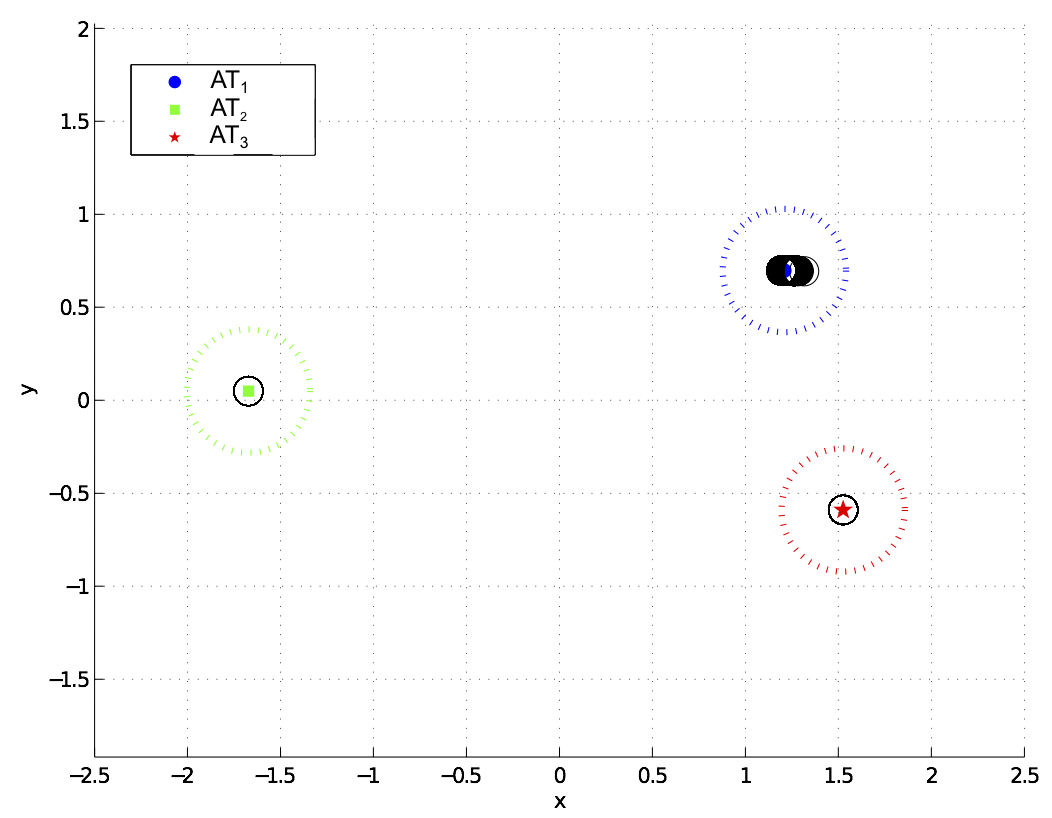}
\caption{Illustration of the classification algorithm, which allows to compute the basins of attractions. The position of the 70 trajectories (black circles) with respect to the three attracting trajectories (symbols) at time $t=550$ kyr,
starting from $t_0 = 20$ kyr, are shown. If the trajectory falls into a dotted circle, classification is possible.
The 70 trajectories form three highly concentrated clusters, corresponding to the  three attracting trajectories.}
\label{Fig:Fig_basin_IC_t_20b}
\end{center}
\end{figure}

The classification algorithm is illustrated in Fig. \ref{Fig:Fig_basin_IC_t_20b},
where a cut has been made at a target time $t=550$ kyr.
The three  $AT_i$'s are displayed, together with the location of the trajectories
(black circles).
To decide if a given trajectory ends up onto a given  $AT_i$,
we choose a maximum distance from the $AT_i$ (dotted circle);
taken here to be $1/4$ of the minimum distance between two $AT_i$'s.
Trajectories falling into that dotted circle are classified as ending into $AT_i$.

We consider a given trajectory starting from $t_0=0$ kyr, integrate it up to $t=550$ kyr, and then we examine its position with respect to the $AT_i$ at the same time T=550 kyr. If the distance to a given $AT_k$ is 'sufficiently small' (see the circles around the $AT_i$), then this initial condition is coloured in the k$^{th}$ color, associated with the k$^{th}$ attracting trajectories $AT_i$.

\begin{figure}
\begin{center}
\includegraphics[width=0.40\textwidth]{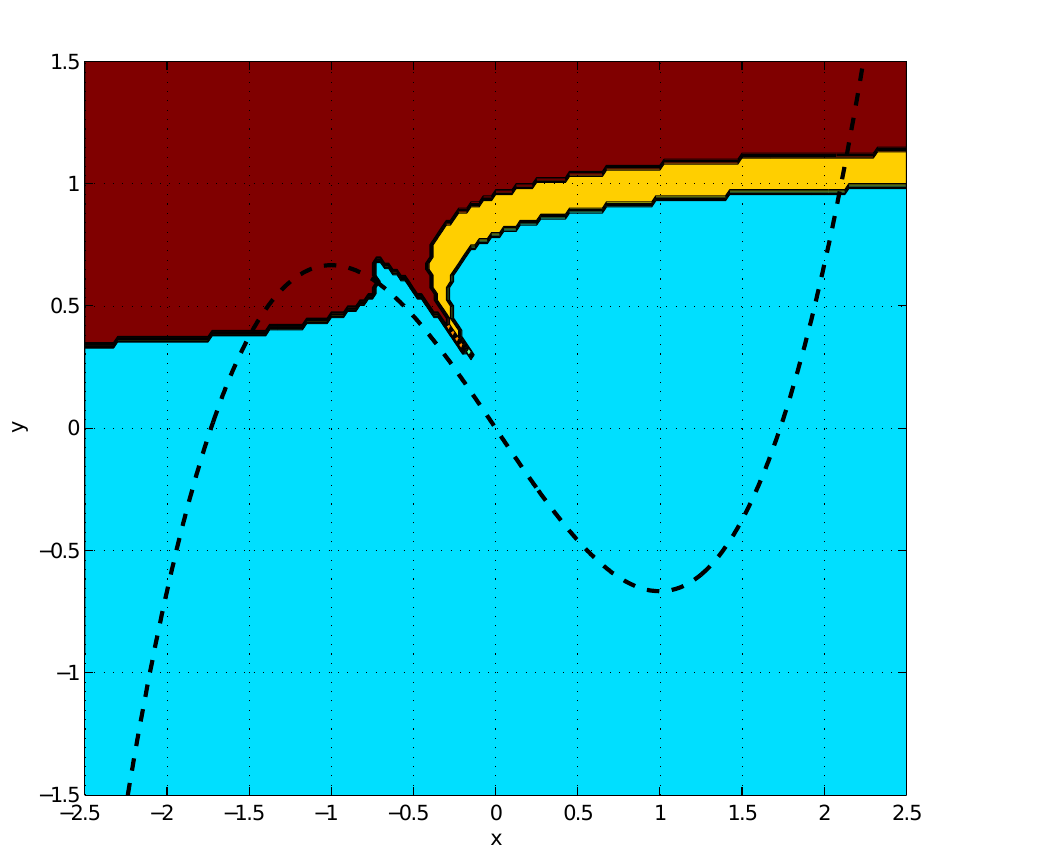}
\caption{When integrating the system system Eq. (\ref{Eq:System}), with many initial conditions (grid $201 \times 121$ = 24 321) starting from time $t_0$ = 0 kyr, with the insolation forcing given by Eq.~(\ref{Eq:Insolation_35_terms}), we obtain the three basins of attraction above; each basin is coloured with the colour of the corresponding $AT_i$. The function $\Phi'(y) = y^3 / 3 -y = x$, corresponding to the slow manifold, is also shown (dashed curve). }
\label{Fig:Basin_Attraction}
\end{center}
\end{figure}

Repeating this process for each initial condition on a fine grid of the whole phase space gives the shape of the basins of attraction. As we have three attracting trajectories, we will have three basins of attraction.
Such basins of attraction obtained are shown in Fig. \ref{Fig:Basin_Attraction}.

As we are in the case of a non-autonomous system, we then have to repeat this procedure for several starting times $t_0$ (the position of the $AT_i$'s are constantly evolving, hence the shape of the basins are also varying with time).
This
has been done to produce the evolving shape of the basins of attraction
in Figs. \ref{Fig:Strip_41k} and \ref{Fig:Strip_insol}.

Note that if $t$ is too close from $t_0$,
transient behaviours predominate, and not enough time has elapsed in order for the trajectory to be attracted by a given $AT_i$, hence the basins of attraction cannot be defined in that case.

The glacial/interglacial cycles do exist since about 3 million years, but  only 8 limit cycles of
100 kyr period have been performed, since the Mid-Pleistocene Transition (MPT). This should however be sufficient to converge onto the attracting trajectories, because the climatic trajectories are rapidly attracted on the limit cycle.
So, if the ice age model Eq. (\ref{Eq:System}) is a realistic one, it would be reasonable to state that the transient of climate dynamics has gone, and that we are currently probably somewhere on an attracting trajectory, if any.

\section{Instantaneous stability}
\label{Annex:Instantaneous_instability}

As we considered the long-term (3 Myr) and short-term ($H=50$ kyr) local stabilities,
we could wonder why not to consider also, at the other extreme, the 'very short-term' local stability
($H \rightarrow 0$).
Even if it has probably less relevance in the context of predictability of glaciations,
it may however still be computed and understood as an
\emph{instantaneous} stability.
As the horizon of time $H$ tends to zero, it means that we no longer consider any motion in the phase space (no average in time), hence the instantaneous stability becomes also a local property in the phase space.

\begin{figure}
\begin{center}
\includegraphics[width=0.48\textwidth]{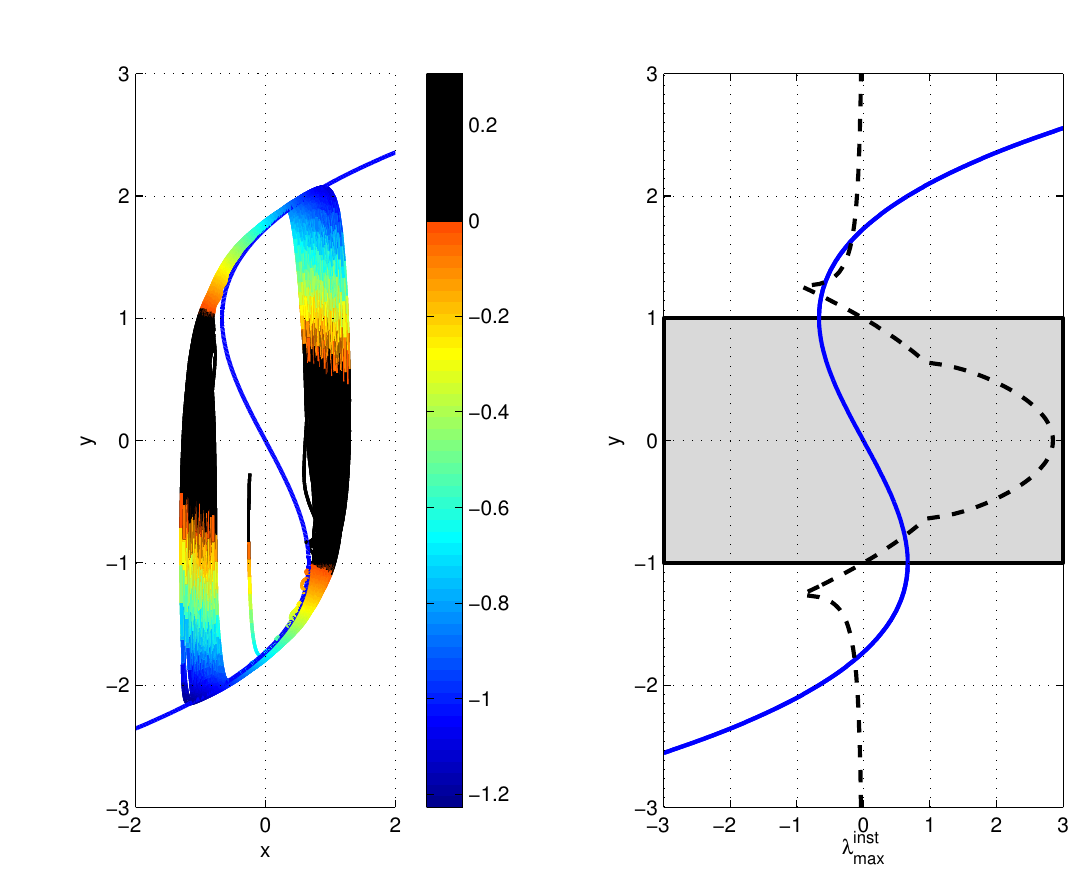}
\caption{Very short-term stability. \emph{Left}
: the attracting climatic trajectory (initial conditions ($x_{0},y_{0}) = (-0.24,-0.27)$) is coloured by the values of $\lambda_{max}^{H = 1}$ for the system Eq. (\ref{Eq:System}) with $\alpha = 11.11$, $\beta=0.25$, $\gamma = 0.033$ and $\tau = 35.09$.
Instantaneous instability ($\lambda_{max}^{H = 1} > 0$)  occurs for $y \in I^* = [-1,1]$.
\emph{Right}: the function $10\times \lambda_{max}^{inst} $ (black dashed curve) has positive values for the same range $y \in I^* = [-1,1]$. For both graphs, the function $\Phi'(y)=x$ has been superimposed for an improved reading and easier interpretation.}
\label{Fig:Link_Lyap_Jac}
\end{center}
\end{figure}

The same formula (Eq.~\ref{Eq:lyap_local}) as for the short-term stability is used,
but now with an Horizon of $H = 1$ kyr, representing the very short-term.
The result is plotted in Fig. \ref{Fig:Link_Lyap_Jac}, where the trajectory has been coloured with the
$\lambda_{max}^{H = 1}$ values.
There are black areas, which means local instantaneous instability.
This black zone is always in the same region of the phase space: when $y$ lies within the interval $I^* = [-1,1]$, whatever the initial condition.

We now demonstrate mathematically why the black zone in Fig.~\ref{Fig:Link_Lyap_Jac}
\emph{left}
is located within the interval $I^* = [-1,1]$.

Let us derive the theoretical expression of the exactly instantaneous LCE $\lambda_i^{inst} $.
If we call $\Lambda_i$ the eigenvalues of the jacobian matrix $J$,
then the  $\lambda_i^{inst} $ are  defined as
$\lambda_i^{inst} =
\lim_{t \rightarrow 0}~
\frac{1}{t} \ln \Lambda_i
$ \cite{Drazin:1992}, so that we can obtain the largest $\lambda_i^{inst}$ by computing
$\lambda_{max}^{inst} = max(\Re(\Lambda_i)) $.

As the Jacobian matrix $J$ is a function of $y$ only (and not $x$), see Eq.~\ref{Eq:Jacobian},
so do  $\lambda_{max}^{inst}$.
The function  $\lambda_{max}^{inst}(y)$ is plotted (Fig. \ref{Fig:Link_Lyap_Jac}, \emph{right}, black dotted curve),
where we clearly see that it is positive on the interval $[-1,1]$,
which is precisely the same interval as $I^*$.

{\begin{flushright} \tiny $\square$ \end{flushright}}

Moreover, the numerical values also match (e.g. the maximum of the function $\lambda_{max}^{inst}$ has a value of about 0.28, which correspond to the maximum of the color scale).

An interpretation of $\lambda_{max}^{inst} $ is that a particular trajectory makes a travel through the phase space  'painted' by $\lambda_{max}^{inst} $, integrating so the instantaneous stability to achieve the short-term and long-term stabilities.


\bibliographystyle{spmpsci}      

\bibliography{BibDesk_BDS_19_sept_2011} 

\begin{thebibliography}{}

\bibitem[Abarbanel et~al., 1991]{Abarbanel:1991}
Abarbanel, H. D.~I., Brown, R., and Kennel, M.~B. (1991).
\newblock Variation of lyapunov exponents on a strange attractor.
\newblock {\em Journal of Nonlinear Science}, 1(2):175--199.

\bibitem[Abarbanel et~al., 1996]{Abarbanel:1996}
Abarbanel, H. D.~I., Rulkov, N.~F., and Sushchik, M.~M. (1996).
\newblock Generalized synchronization of chaos: The auxiliary system approach.
\newblock {\em Phys. Rev. E}, 53(5):4528--4535.

\bibitem[Arnold, 1983]{Arnold:1983}
Arnold, V. (1983).
\newblock {\em Geometrical Methods in the Theory of Ordinary Differential
  Equations}.
\newblock Springer-Verlag, New York, 1988 second edition.
\newblock English translation of the original russian publication:
  "Dopolnitel'nye Glavy Teorii Obyknovennykh Differentsial'nykh Uravneni{\^\i}"
  (Additional Chapters to the Theory of Ordinary Differential Equations,
  Moscow: Nauka, 1978).

\bibitem[Ashkenazy, 2006]{Ashkenazy:2006}
Ashkenazy, Y. (2006).
\newblock The role of phase locking in a simple model for glacial dynamics.
\newblock {\em Climate Dynamics}, 27(4):421--431.

\bibitem[Balanov et~al., 2009]{Balanov:2009}
Balanov, A., Janson, N., Postnov, D., and Sosnovtseva, O. (2009).
\newblock {\em Synchronization: From Simple to Complex}.
\newblock Springer-Verlag, Berlin, Germany.

\bibitem[Barnes and Grimshaw, 1997]{Barnes:1997}
Barnes, B. and Grimshaw, R. (1997).
\newblock Analytical and numerical studies of the bonhoeffer van der pol
  system.
\newblock {\em The ANZIAM Journal}, 38(04):427--453.

\bibitem[Belogortsev, 1992]{Belogortsev:1992}
Belogortsev, A.~B. (1992).
\newblock Quasiperiodic resonance and bifurcations of tori in the weakly
  nonlinear duffing oscillator.
\newblock {\em Physica D: Nonlinear Phenomena}, 59(4):417--429.

\bibitem[Benettin et~al., 1980]{Benettin:1980fk}
Benettin, G., Galgani, L., Giorgilli, A., and Strelcyn, J.-M. (1980).
\newblock Lyapunov characteristic exponents for smooth dynamical systems and
  for hamiltonian systems; a method for computing all of them. part 2:
  Numerical application.
\newblock {\em Meccanica}, 15(1):21--30.

\bibitem[Beno{\^\i}t et~al., 1981]{Benoit:1981}
Beno{\^\i}t, E., Callot, J., Diener, F., and Diener, M. (1981).
\newblock Chasse au canard.
\newblock {\em Collectanea Mathematica}, 31-32(1-3):37--119.

\bibitem[Berger, 1978]{Berger:1978}
Berger, A.~L. (1978).
\newblock Long-term variations of daily insolation and {Q}uaternary climatic
  changes.
\newblock {\em J. Atmos. Sci.}, 35:2362--2367.

\bibitem[Broecker and van Donk, 1970]{Broecker:1970}
Broecker, W.~S. and van Donk, J. (1970).
\newblock Insolation changes, ice volumes, and the {O18} record in deep-sea
  cores.
\newblock {\em Rev. Geophys.}, 8(1):169--198.

\bibitem[Broer and Sim{\'o}, 1998]{Broer:1998}
Broer, H.~W. and Sim{\'o}, C. (1998).
\newblock Hill's equation with quasi-periodic forcing: resonance tongues,
  instability pockets and global phenomena.
\newblock {\em Soc. Brasil Mat}, (29):253293.

\bibitem[Brown and Kocarev, 2000]{Brown:2000}
Brown, R. and Kocarev, L. (2000).
\newblock A unifying definition of synchronization for dynamical systems.
\newblock {\em Chaos}, 10(2):344--349.

\bibitem[Bryant et~al., 1990]{Bryant:1990}
Bryant, P., Brown, R., and Abarbanel, H. D.~I. (1990).
\newblock Lyapunov exponents from observed time series.
\newblock {\em Physical Review Letters}, 65(13).

\bibitem[Chen and Chen, 2008]{Chen:2008}
Chen, J.-H. and Chen, W.-C. (2008).
\newblock Chaotic dynamics of the fractionally damped van der pol equation.
\newblock {\em Chaos, Solitons \& Fractals}, 35(1):188--198.

\bibitem[{Crucifix}, 2011]{Crucifix:2011a}
{Crucifix}, M. (2011).
\newblock {Oscillators and relaxation phenomena in Pleistocene climate theory}.
\newblock {\em Transactions of the Philosophical Transactions of the Royal
  Society (Series A, Physical Mathematical and Engineering Sciences)}, In
  Press.

\bibitem[D'Acunto, 2006]{DAcunto:2006}
D'Acunto, M. (2006).
\newblock Determination of limit cycles for a modified van der pol oscillator.
\newblock {\em Mechanics Research Communications}, 33(1):93--98.

\bibitem[De~Saedeleer et~al., 2010]{De-Saedeleer:2010}
De~Saedeleer, B., Crucifix, M., and Wieczorek, S. (2010).
\newblock Is the synchronization of the climatic system on the orbital forcing
  robust?
\newblock Poster presented at the "SYNCLINE 2010: Synchronization in Complex
  Networks" conference, held on 26--29 May 2010 at the Physikzentrum Bad Honnef
  (Germany).

\bibitem[Degli Esposti~Boschi et~al., 2002]{Degli-Esposti-Boschi:2002}
Degli Esposti~Boschi, C., Ortega, G.~J., and Louis, E. (2002).
\newblock Discriminating dynamical from additive noise in the van der pol
  oscillator.
\newblock {\em Physica D: Nonlinear Phenomena}, 171(1-2):8--18.

\bibitem[Dijkstra et~al., 2003]{Dijkstra:2003}
Dijkstra, H.~A., Weijer, W., and Neelin, J.~D. (2003).
\newblock Imperfections of the three-dimensional thermohaline circulation:
  hysteresis and unique-state regimes.
\newblock {\em J. Phys. Oceanogr.}, 33:2796--2814.

\bibitem[Doedel et~al., 2009]{Doedel:2009}
Doedel, E., Champneys, A., Dercole, F., Fairgrieve, T., Kuznetsov, Y., Oldeman,
  B., Paffenroth, R., Sandstede, B., Wang, X., and Zhang, C. (2009).
\newblock Auto: Software for continuation and bifurcation problems in ordinary
  differential equations.
\newblock Technical report, Montreal.

\bibitem[Donges et~al., 2009]{Donges:2009}
Donges, J.~F., Zou, Y., Marwan, N., and Kurths, J. (2009).
\newblock The backbone of the climate network.
\newblock {\em EPL (Europhysics Letters)}, 87(4):48007.

\bibitem[Drazin and King, 1992]{Drazin:1992}
Drazin, P.~G. and King, G.~P., editors (1992).
\newblock {\em Interpretation of Time Series from Nonlinear Systems},
  volume~58.

\bibitem[Eckmann and Ruelle, 1985]{Eckmann:1985}
Eckmann, J.~P. and Ruelle, D. (1985).
\newblock Ergodic-theory of chaos and strange attractors.
\newblock {\em Reviews of Modern Physics}, 57(3):617--656.

\bibitem[Ginoux and Rossetto, 2006]{Ginoux:2006}
Ginoux, J.-M. and Rossetto, B. (2006).
\newblock Differential geometry and mechanics: Applications to chaotic
  dynamical systems.
\newblock {\em I. J. Bifurcation and Chaos}, 16(4):887--910.

\bibitem[{Glass} and {Sun}, 1994]{Glass:1994}
{Glass}, L. and {Sun}, J. (1994).
\newblock {Periodic forcing of a limit-cycle oscillator: Fixed points, Arnold
  tongues, and the global organization of bifurcations}.
\newblock {\em Phys. Rev. E}, 50:5077--5084.

\bibitem[Glendinning and Wiersig, 1999]{Glendinning:1999}
Glendinning, P. and Wiersig, J. (1999).
\newblock Fine structure of mode-locked regions of the quasi-periodically
  forced circle map.
\newblock {\em Physics Letters A}, 257(1-2):65--69.

\bibitem[{Grasman} et~al., 2005]{Grasman:2005}
{Grasman}, J., {Verhulst}, F., and {Shih}, S. (2005).
\newblock {The Lyapunov exponents of the Van der Pol oscillator}.
\newblock {\em Mathematical Methods in the Applied Sciences}, 28:1131--1139.

\bibitem[Grassberger and Procaccia, 1983]{Grassberger:1983}
Grassberger, P. and Procaccia, I. (1983).
\newblock Measuring the strangeness of strange attractors.
\newblock {\em Physica D: Nonlinear Phenomena}, 9(1-2):189--208.

\bibitem[Guckenheimer and Haiduc, 2005]{Guckenheimer:2005}
Guckenheimer, J. and Haiduc, R. (2005).
\newblock Canards at folded node.
\newblock {\em Mosc. Math. J}, 5:91--103.

\bibitem[Guckenheimer et~al., 2000]{Guckenheimer:2000}
Guckenheimer, J., Hoffman, K., and Weckesser, W. (2000).
\newblock Numerical computation of canards.
\newblock {\em International Journal of Bifurcation and Chaos},
  10(12):2269--2687.

\bibitem[Guckenheimer and Holmes, 1983]{Guckenheimer:1983}
Guckenheimer, J. and Holmes, P. (1983).
\newblock {\em Nonlinear oscillations, dynamical systems, and bifurcations of
  vector fields}.
\newblock Springer-Verlag, New York.

\bibitem[Haken, 1983]{Haken:1983}
Haken, H. (1983).
\newblock At least one lyapunov exponent vanishes if the trajectory of an
  attractor does not contain a fixed point.
\newblock {\em Physics Letters A}, 94(2):71--72.

\bibitem[Hays et~al., 1976]{Hays:1976}
Hays, J.~D., Imbrie, J., and Shackleton, N.~J. (1976).
\newblock Variations in the {E}arth's orbit : {P}acemaker of ice ages.
\newblock {\em Science}, 194:1121--1132.

\bibitem[Hilborn, 2000]{Hilborn:2000}
Hilborn, R. (2000).
\newblock {\em Chaos and Nonlinear Dynamics: an Introduction for Scientists and
  Engineers}.
\newblock Oxford University Press.

\bibitem[Huybers, 2007]{Huybers:2007}
Huybers, P. (2007).
\newblock Glacial variability over the last two millions years: an extended
  depth-derived age model, continous obliquity pacing, and the {P}leistocene
  progression.
\newblock {\em Quaternary Sci. Rev.}, 26:37--55.

\bibitem[Imbrie and Imbrie, 1980]{Imbrie:1980}
Imbrie, J. and Imbrie, J.~Z. (1980).
\newblock Modelling the climatic response to orbital variations.
\newblock {\em Science}, 207:943--953.

\bibitem[Kantz and Schreiber, 2004]{Kantz:2004}
Kantz, H. and Schreiber, T. (2004).
\newblock {\em Nonlinear Time Series Analysis}.
\newblock Cambridge University Press, Cambridge, U.K., 2nd edition.

\bibitem[Kaplan and Yorke, 1979]{Kaplan:1979}
Kaplan, J.~L. and Yorke, J.~A. (1979).
\newblock Chaotic behavior of multidimensional difference equations.
\newblock In {\em Functional differential equations and approximation of fixed
  points ({P}roc. {S}ummer {S}chool and {C}onf., {U}niv. {B}onn, {B}onn,
  1978)}, volume 730 of {\em Lecture Notes in Math.}, pages 204--227. Springer,
  Berlin.

\bibitem[Kloeden, 2000]{Kloeden:2000}
Kloeden, P.~E. (2000).
\newblock {A Lyapunov function for pullback attractors of nonautonomous
  differential equations}.
\newblock {\em Electronic J. Diff. Eqns, Conf. 05}, pages 91--102.

\bibitem[Kosmidis and Pakdaman, 2003]{Kosmidis:2003}
Kosmidis, E.~K. and Pakdaman, K. (2003).
\newblock An analysis of the reliability phenomenon in the fitzhugh-nagumo
  model.
\newblock {\em Journal of Computational Neuroscience}, 14(1):5--22.

\bibitem[Langa et~al., 2002]{Langa:2002}
Langa, J.~A., Robinson, J.~C., and Su{\'a}rez, A. (2002).
\newblock Stability, instability, and bifurcation phenomena in non-autonomous
  differential equations.
\newblock {\em Nonlinearity}, 15(3):1--17.

\bibitem[Laskar et~al., 2004]{Laskar:2004}
Laskar, J., Robutel, P., Joutel, F., Boudin, F., Gastineau, M., Correia, A.
  C.~M., and Levrard, B. (2004).
\newblock A long-term numerical solution for the insolation quantities of the
  {E}arth.
\newblock {\em Astrom. Astroph.}, 428:261--285.

\bibitem[Le~Treut and Ghil, 1983]{Le-Treut:1983}
Le~Treut, H. and Ghil, M. (1983).
\newblock Orbital forcing, climatic interactions and glaciation cycles.
\newblock {\em J. Geophys. Res.}, 88(C9):5167--5190.

\bibitem[Lichtenberg and Lieberman, 1983]{Lichtenberg:1983}
Lichtenberg, A.~J. and Lieberman, M.~A. (1983).
\newblock {\em Regular and stochastic motion}.
\newblock Springer-Verlag, New York.

\bibitem[Lisiecki and Raymo, 2005]{Lisiecki:2005}
Lisiecki, L.~E. and Raymo, M.~E. (2005).
\newblock A pliocene-pleistocene stack of 57 globally distributed benthic
  {$\delta^{18}O$} records.
\newblock {\em Paleoceanography}, 20(1).

\bibitem[Lisiecki and Raymo, 2007]{Lisiecki:2007}
Lisiecki, L.~E. and Raymo, M.~E. (2007).
\newblock Plio-pleistocene climate evolution: trends and transitions in glacial
  cycle dynamics.
\newblock {\em Quaternary Science Reviews}, 26(1-2):56--69.

\bibitem[Liu et~al., 2005]{Liu:2005}
Liu, H.-F., Dai, Z.-H., Li, W.-F., Gong, X., and Yu, Z.-H. (2005).
\newblock Noise robust estimates of the largest lyapunov exponent.
\newblock {\em Physics Letters A}, 341(1-4):119--127.

\bibitem[Luethi et~al., 2008]{Luethi:2008}
Luethi, D., Le~Floch, M., Bereiter, B., Blunier, T., Barnola, J.-M.,
  Siegenthaler, U., Raynaud, D., Jouzel, J., Fischer, H., Kawamura, K., and
  Stocker, T.~F. (2008).
\newblock High-resolution carbon dioxide concentration record 650,000-800,000
  years before present.
\newblock {\em Nature}, 453(7193):379--382.

\bibitem[Marwan et~al., 2009]{Marwan:2009}
Marwan, N., Donges, J.~F., Zou, Y., Donner, R.~V., and Kurths, J. (2009).
\newblock Complex network approach for recurrence analysis of time series.
\newblock {\em Physics Letters A}, 373(46):4246--4254.

\bibitem[McCaffrey et~al., 1992]{McCaffrey:1992}
McCaffrey, D.~F., Ellner, S., Gallant, A.~R., and Nychka, D.~W. (Sep., 1992).
\newblock Estimating the lyapunov exponent of a chaotic system with
  nonparametric regression.
\newblock {\em Journal of the American Statistical Association},
  87(419):682--695.

\bibitem[Mettin et~al., 1993]{Mettin:1993}
Mettin, R., Parlitz, U., and Lauterborn, W. (1993).
\newblock Bifurcation structure of the driven van der pol oscillator.
\newblock {\em Int. J. Bifurcation \& Chaos}, (3):1529--1555.

\bibitem[Milankovitch, 1941]{Milankovitch:1941}
Milankovitch, M. (1941).
\newblock {\em Kanon der Erdbestrahlung und Seine Anwendung auf das
  Eiszeitenproblem (Canon of insolation and the ice-age problem)}.
\newblock K{\"o}nlishe Serbische Akademie, Belgrad.

\bibitem[M{\"u}ller, 1995]{Muller:1995}
M{\"u}ller, P.~C. (1995).
\newblock Calculation of lyapunov exponents for dynamic systems with
  discontinuities.
\newblock {\em Chaos, Solitons \& Fractals}, 5(9):1671--1681.

\bibitem[Oseledec, 1968]{Oseledec:1968}
Oseledec, V. (1968).
\newblock A multiplicative ergodic theorem: Ljapunov characteristic numbers for
  dynamical systems.
\newblock {\em Transactions of the Moscow Mathematical Society}, 19:197--231.

\bibitem[{Osinga} et~al., 2000]{Osinga:2000}
{Osinga}, H., {Wiersig}, J., {Glendinning}, P., and {Feudel}, U. (2000).
\newblock {Multistability and nonsmooth bifurcations in the quasiperiodically
  forced circle map}.
\newblock {\em ArXiv Nonlinear Sciences e-prints}.

\bibitem[Ott, 2002]{Ott:2002}
Ott, E. (2002).
\newblock {\em Chaos in Dynamical Systems}.
\newblock Cambridge University Press.

\bibitem[Paillard, 1998]{Paillard:1998}
Paillard, D. (1998).
\newblock The timing of pleistocene glaciations from a simple multiple-state
  climate model.
\newblock {\em Nature}, 391:378--381.

\bibitem[Paillard and Parrenin, 2004]{Paillard:2004}
Paillard, D. and Parrenin, F. (2004).
\newblock The {A}ntarctic ice sheet and the triggering of deglaciations.
\newblock {\em Earth Planet. Sc. Lett.}, 227:263--271.

\bibitem[Parlitz and Lauterborn, 1987]{Parlitz:1987}
Parlitz, U. and Lauterborn, W. (1987).
\newblock Period-doubling cascades and devil's staircases of the driven van der
  pol oscillator.
\newblock {\em Physical Review A}, 36(3).

\bibitem[Pikovsky et~al., 2001]{Pikovsky:2001}
Pikovsky, A., Rosenblum, M., and Kurths, J. (2001).
\newblock {\em Synchronization A Universal Concept in Nonlinear Sciences}.
\newblock Cambridge University Press, New York.

\bibitem[Rahmstorf et~al., 2005]{Rahmstorf:2005}
Rahmstorf, S., Crucifix, M., Ganopolski, A., Goosse, H., Kamenkovich, I.,
  Knutti, R., Lohmann, G., Marsh, R., Mysak, L.~A., Wang, Z., and Weaver, A.~J.
  (2005).
\newblock Thermohaline circulation hysteresis: A model intercomparison.
\newblock {\em Geophys. Res. Lett.}, 32(23).

\bibitem[Ramasubramanian and Sriram, 2000]{Ramasubramanian:2000}
Ramasubramanian, K. and Sriram, M.~S. (2000).
\newblock A comparative study of computation of lyapunov spectra with different
  algorithms.
\newblock {\em Physica D: Nonlinear Phenomena}, 139(1-2):72--86.

\bibitem[Rosenstein et~al., 1993]{Rosenstein:1993}
Rosenstein, M.~T., Collins, J.~J., and Luca, C. J.~D. (1993).
\newblock A practical method for calculating largest lyapunov exponents from
  small datasets.
\newblock {\em Physica D}, (124).

\bibitem[Rossler, 1979]{Rossler:1979}
Rossler, O.~E. (1979).
\newblock An equation for hyperchaos.
\newblock {\em Physics Letters A}, 71(2-3):155--157.

\bibitem[Ruelle, 1990]{Ruelle:1990}
Ruelle, D. (1990).
\newblock Deterministic chaos: the science and the fiction.
\newblock {\em Proceedings of the Royal Society A, London}, (427):241--248.

\bibitem[Ruihong et~al., 2008]{Ruihong:2008}
Ruihong, L., Wei, X., and Shuang, L. (2008).
\newblock Chaos control and synchronization of the $\phi^6$-van der pol system
  driven by external and parametric excitations.
\newblock {\em Nonlinear Dynamics}, 53(3):261--271.

\bibitem[Rulkov et~al., 1995]{Rulkov:1995}
Rulkov, N.~F., Sushchik, M.~M., Tsimring, L.~S., and Abarbanel, H. D.~I.
  (1995).
\newblock Generalized synchronization of chaos in directionally coupled chaotic
  systems.
\newblock {\em Phys. Rev. E}, 51(2):980--994.

\bibitem[Saltzman, 2002]{Saltzman:2002}
Saltzman, B. (2002).
\newblock {\em Dynamical Paleoclimatology: Generalized Theory of Global Climate
  Change (International Geophysics)}.
\newblock {Academic Press}.

\bibitem[Saltzman et~al., 1984]{Saltzman:1984}
Saltzman, B., Hansen, A.~R., and Maasch, K.~A. (1984).
\newblock The late {Q}uaternary glaciations as the response of a 3-component
  feedback-system to {E}arth-orbital forcing.
\newblock {\em Journal of the Atmospheric Sciences}, 41(23):3380--3389.

\bibitem[Saltzman and Maasch, 1990]{Saltzman:1990}
Saltzman, B. and Maasch, K.~A. (1990).
\newblock A first-order global model of late {C}enozoic climate.
\newblock {\em Trans. R. Soc. Edinburgh Earth Sci}, 81:315--325.

\bibitem[Saltzman and Maasch, 1991]{Saltzman:1991}
Saltzman, B. and Maasch, K.~A. (1991).
\newblock A first-order global model of late {C}enozoic climate. {II} further
  analysis based on a simplification of the {CO}$_2$ dynamics.
\newblock {\em Clim. Dyn.}, 5:201--210.

\bibitem[Savi, 2005]{Savi:2005}
Savi, M.~A. (2005).
\newblock Chaos and order in biomedical rhythms.
\newblock {\em Journal of the Brazilian Society of Mechanical Sciences and
  Engineering}, 27(2):157--169.

\bibitem[Shimada and Nagashima, 1979]{Shimada:1979}
Shimada, I. and Nagashima, T. (1979).
\newblock A numerical approach to ergodic problem of dissipative dynamical
  systems.
\newblock {\em Progr. Theoret. Phys.}, 61(6):1605--1616.

\bibitem[Strogatz, 1994]{Strogatz:1994}
Strogatz, S.~H. (1994).
\newblock {\em Nonlinear Dynamics And Chaos: With Applications To Physics,
  Biology, Chemistry, And Engineering (Studies in Nonlinearity)}.
\newblock Studies in nonlinearity. Perseus Books Group, 1 edition.

\bibitem[Svensson and Coombes, 2009]{Svensson:2009}
Svensson, C.-M. and Coombes, S. (2009).
\newblock Mode locking in a spatially extended neuron model: active soma and
  compartmental tree.
\newblock {\em International Journal of Bifurcation and Chaos},
  19(8):2597--2607.

\bibitem[Theiler, 1990]{Theiler:1990}
Theiler, J. (1990).
\newblock Estimating fractal dimension.
\newblock {\em J. Opt. Soc. Am. A}, 7(6):1055--1073.

\bibitem[Tsiganis et~al., 1999]{Tsiganis:1999}
Tsiganis, K., Anastasiadis, A., and Varvoglis, H. (1999).
\newblock Effective lyapunov numbers and correlation dimensions in a {3-D}
  hamiltonian system.
\newblock In Henrard, J. and Ferraz-Mello, S., editors, {\em IAU-Colloquium No
  172 - The impact of modern dynamics in Astronomy}, pages 447--448.

\bibitem[Tziperman and Gildor, 2003]{Tziperman:2003}
Tziperman, E. and Gildor, H. (2003).
\newblock On the mid-{P}leistocene transtion to 100-kyr glacial cycles and the
  asymmetry between glaciation and deglaciation times.
\newblock {\em Paleoceanography}, 18(1):1001.

\bibitem[Tziperman et~al., 2006]{Tziperman:2006}
Tziperman, E., Raymo, M.~E., Huybers, P., and Wunsch, C. (2006).
\newblock Consequences of pacing the {P}leistocene 100 kyr ice ages by
  nonlinear phase locking to {M}ilankovitch forcing.
\newblock {\em Paleoceanography}, 21:PA4206.

\bibitem[van~der Pol, 1926]{vanderPol:1926}
van~der Pol, B. (1926).
\newblock On relaxation oscillations.
\newblock {\em Phil. Mag.}, 2(11):978--992.

\bibitem[Wieczorek, 2009]{Wieczorek:2009}
Wieczorek, S. (2009).
\newblock Stochastic bifurcation in noise-driven lasers and hopf oscillators.
\newblock {\em Phys. Rev. E}, 79(3):036209.

\bibitem[{Wieczorek}, 2011]{Wieczorek:2011}
{Wieczorek}, S.~M. (2011).
\newblock {Noise synchronisation and stochastic bifurcations in lasers}.
\newblock {\em http://arxiv.org/abs/1104.4052}.

\bibitem[Wiggins, 2003]{Wiggins:2003}
Wiggins, S. (2003).
\newblock {\em {Introduction to Applied Nonlinear Dynamical Systems and
  Chaos}}.
\newblock Texts in Applied Mathematics. Springer, 2nd edition.

\bibitem[Wolf et~al., 1985]{Wolf:1985}
Wolf, A., Swift, J.~B., Swinney, H.~L., and Vastano, J.~A. (1985).
\newblock Determining lyapunov exponents from a time series.
\newblock {\em Physica D: Nonlinear Phenomena}, 16(3):285--317.

\bibitem[Wu et~al., 2006]{Wu:2006}
Wu, L., Zhu, S., and Li, J. (2006).
\newblock Synchronization on fast and slow dynamics in drive-response systems.
\newblock {\em Physica D: Nonlinear Phenomena}, 223(2):208--213.

\end{thebibliography}


\end{document}